\documentclass[11pt]{amsart}

\newcommand{\A}{\mbox{${\mathcal A}$}}

\newcommand{\F}{\mbox{${\mathcal F}$}}

\theoremstyle{remark}

\theoremstyle{definition}

\numberwithin{equation}{section}

\begin{document}

\title[Polynomial Ergodicity and asymptotic behaviour of
unbounded solutions]{{Polynomial ergodicity and asymptotic
behaviour of unbounded solutions of abstract evolution equations}}
\author{Bolis Basit}
       \address {School of Math. Sciences, P.O. Box  28M,
       Monash University, Victoria 3800, Australia}
        \email{ bolis.basit@sci.monash.edu.au and alan.pryde)@sci.monash.edu.au}
\author{A. J. Pryde}

\begin{abstract}
{ In this paper we develop the notion of ergodicity to include
functions dominated by a weight $w$. Such functions have
polynomial means and include, amongst many others, the $w$-almost
periodic functions. This enables us to describe the asymptotic
behaviour of unbounded solutions of linear evolution, recurrence
and convolution equations. To unify the treatment and allow for
further applications, we consider solutions $\phi \,: G\rightarrow
X$ of generalized evolution equations of the form $\,\,
(*)\,\,(B\phi )(t)=A\phi (t)+\psi (t)$ for $t\in G$ where $G$\ is
a locally compact abelian group with a closed subsemigroup $J$,
$A$ is a closed linear operator on a Banach space $X$, $\psi
:G\rightarrow X$ is continuous and $B$ is a linear operator with
characteristic function $\theta _{B}:\widehat{G}\rightarrow
\mathbf{C}$. We introduce the resonance set $\theta
_{B}^{-1}(\sigma (A))$ which contains the Beurling spectra of all
solutions of the homogeneous equation $B\phi
=A\circ \phi$. For certain classes ${\F} $ of functions from $J$ to $%
X,$ the spectrum $sp_{{\F}}(\phi )$ of $\phi $ relative to $%
{\F} $ is used to determine membership of ${\F}.$ Our main result
gives general conditions under which $sp_{{\F} }(\phi )$\ is a
subset of the resonance set. As a simple consequence we obtain
conditions under which $\psi |_{J}\in \F$ implies $\phi |_{J}\in
{\F}.$ An important tool is our generalization to unbounded
functions of a theorem of Loomis. As applications we obtain
generalizations or new proofs of many known results, including
theorems of Gelfand, Hille, Katznelson-Tzafriri, Esterle et al.,
Ph\'{o}ng, Ruess and Arendt-Batty.

\smallskip

\noindent \textit {1991 Mathematics subject classification}:
primary {46J20, 43A60}; secondary {47A35, 34K25, 28B05.}

\smallskip

\noindent\textit {Keywords and phrases}: Weighted ergodicity,
orbits of unbounded semigroup representation, non-quasianalytic
weights, stability, Beurling spectrum. }
 \end{abstract}
\maketitle

\bigskip

\textbf{1. Introduction}

In this paper we develop techniques to describe the asymptotic
behaviour of unbounded solutions of linear evolution, recurrence
and convolution equations. Our initial motivation was the study of
resonance phenomena for solutions of equations of the form $\phi
^{\prime
}(t)=A\phi (t)+\psi (t)$ where $A$ is the generator of a $C_{0}$%
-semigroup on a Banach space $X$. Many authors \cite{AFR},
\cite{ABj}, \cite{ABk}, \cite{ESZ}, \cite{ESZJ}, \cite{BBA},
\cite{BBAS}, \cite{BG}, \cite{BPA}, \cite{BPB}, \cite{BPJ},
\cite{BPH}, \cite{BC}, \cite{BNR}, \cite{DP}, [43, pp. 92-96],
\cite{RS}, \cite{RP} have studied the bounded uniformly continuous
solutions of such equations. However, for unbounded solutions
\cite{AR},
   \cite{BPB},  \cite{BPJ},  \cite{BY},
   \cite{PV} much less is known. In  \cite{BG},  \cite{BBG}, \cite{BGH} unbounded solutions with bounded
   n-th mean  are investigated. Here
we obtain results for more general operators $A$ and with $\phi
^{\prime }$ replaced by $B\phi =\sum_{j=0}^{m}b_{j}\phi ^{(j)}.$

\smallskip

An important tool in the bounded case is a theorem of Loomis
\cite{LH} and its subsequent extensions  \cite{BBA}, \cite{BAG},
[43, p. 92]. We prove a generalization for unbounded functions
(Corollary 6.8). To do this we introduce $w$-bounded ergodic
functions (section 4) and $w$-bounded almost periodic functions
(section 5) for certain weights $w.$ In particular, we have found
a connection between the spectral synthesis of functions with
respect to closed subsets $\Delta $ of characters and
$w$-ergodicity (Corollary 6.13). Such a connection does not appear
to be well-known even for the bounded case $w=1$ (see [15, Theorem
5.5]). However, it allows us to reduce results relating spectral
synthesis and asymptotic stability of $C_{0}$-semigroups to the
generalized Loomis theorem. Moreover, we obtain general conditions
under
which solutions of the equation $B\phi =A\circ \phi +\psi $ are $w$%
-almost periodic.

\smallskip

Secondly, we wished to explore parallels between the problems
mentioned above and convolution equations of the form $k\ast \phi
=A\circ \phi +\psi $ on $\mathbf {R}^{d}$ or $\mathbf{Z}^{d}.$
Using the generalized Loomis theorem, we give conditions under
which $w$-bounded solutions are $w$-almost periodic.

\smallskip

Thirdly, we sought to generalize results of Gelfand and
Katznelson-Tzafriri concerning power-bounded operators (see
\cite{KT}, \cite{AR}, \cite{PQ}). This is achieved for elements of
unital
Banach algebras whose powers are $w$-bounded. For example we show that if $%
\sigma (x)=\{1\}$ and $\left\| x^{n}\right\| \leq cw(n),$ where $w$ has
polynomial growth of order $N$, then $(x-e)^{N+1}=0.$ This is shown to be a
special case of our results on $w$-bounded solutions of recurrence equations
of the form $B\phi =A\circ \phi +\psi $ on $\mathbf{Z}$, where $%
(B\phi )(n)=\sum_{j=0}^{m}b_{j}\phi (n+n_{j}).$

\smallskip

To unify the treatment of these various problems and allow for
further applications, we consider $w$-bounded solutions $\phi
:G\rightarrow X$ of general evolution equations of the form

\bigskip

 (*) $\qquad \qquad(B\phi )(t)=A\phi (t)+\psi (t)$ for
$t\in G.$

\bigskip

Here and throughout this paper $G$ denotes a locally compact abelian
topological group with a fixed Haar measure $\mu $ and dual $\widehat{G};$ $%
J $ is a closed subsemigroup of $G$ with non-empty interior such that $%
G=J-J; $ and $X$ is a complex Banach space. Moreover, $A$ is a closed linear
operator on $X,$ $B$ is a linear operator with characteristic function $%
\theta _{B}:\widehat{G}\rightarrow \mathbf{C},$ and ${\F} $ is a
translation invariant class of $w$-bounded continuous functions from $J$ to $%
X$. The asymptotic behaviour of a solution $\phi $ of (*) is
described by properties of the form: $\psi |_{J}\in {\F} $\
implies $\phi |_{J}\in {\F}$.\ A set of characters $sp_{{\F}
}(\phi ),$ the spectrum of
$\phi $ relative to ${\F} ,$ is used to determine membership of $%
{\F}$. Indeed, under certain conditions, $sp_{{\F}}(\phi
)=\emptyset $ if and only if $\phi |_{J}\in {\F}$. Of course it is
also of interest to study solutions of equations of the form
$B\phi =A\circ \phi +\psi $ on a semigroup $J.$ In many
applications this is achieved by extending solutions on $J$ to
solutions on $G$. See for example Corollary 12.6 and Theorem 11.1.

\smallskip

Our main theorem gives general conditions under which $\psi
|_{J}\in {\F} $ implies $sp_{{\F} }(\phi )\subseteq \theta
_{B}^{-1}(\sigma (A)).$ Simple conditions on the spectra then imply that $%
\phi |_{J}\in {\F}$. This theorem is used, sometimes in
conjunction with the generalized Loomis theorem, in each of the
applications mentioned above. We refer to $\theta _{B}^{-1}(\sigma
(A))$ as the \textit{characteristic or resonance set }of (*). It
contains the Beurling spectra of all solutions on $G$\ of the
homogeneous equation $B\phi =A\circ \phi$.

\smallskip

Though other authors use different characterizations of ergodicity
for bounded functions (see \cite{BBA}, \cite{L}, \cite{RP},
\cite{ABj}), we use that of Maak \cite{M}, \cite{MW} (see also
\cite{E},  \cite{I}, \cite{J}, \cite{DM}, \cite{BPA}) because of
its simplicity and wide applicability. A $w$-bounded continuous
function $\phi :J\rightarrow X$ is called $w$\textit{-ergodic}
with mean $p$ if $(\phi -p)/w$ is bounded and ergodic with mean
$0$ in the sense of Maak for some polynomial $p.$
Many properties and applications of ergodic functions are also valid for $w$%
-ergodic ones. A study of polynomials $p:G\rightarrow X$ was
commenced in \cite{BP} and is continued here. In particular, for
certain $w,$ polynomials are the $w$-bounded functions with
Beurling spectrum $\{1\}$. Functions with finite Beurling spectra
are sums of products of characters and polynomials. As a
consequence we obtain a characterization of the minimal primary
ideals in the associated Beurling algebras. In turn, this is used
to characterize functions for which $sp_{ {\F} }(\phi )$ is a
singleton. Thence we prove the generalized Loomis theorem. In
particular, $w$-uniformly continuous totally $w$-ergodic functions
with residual Beurling spectra are $w$-almost periodic.

\smallskip

Some of the notation we use is as follows.

\smallskip

The \textit{translate} $\phi
_{h}$ and \textit{difference} $\Delta _{h}\phi $ by $h\in J$ of \vspace{%
1pt}a function $\phi :J\rightarrow X$ are given by $\phi
_{h}(t)=\phi (t+h)$ and $\Delta _{h}\phi =\phi _{h}-\phi$.

$\mathbf{R}_- =(-\infty, 0]$, $\,\, \mathbf{R}_+ =[0, \infty)$.

For $J\in \{\mathbf{R_+,R}\}$, we define the \textit{indefinite
integral} or \textit{primitive }$P\phi $ of a function $\phi \in
L^1_{loc}(J,X)$ by $P\phi (t)=\int_{0}^{t}\phi (s)ds$.

Similarly, for $J\in \{\mathbf{Z_+,Z}\}$ and
  $\phi \in X^{J}$, we
define  the $sum\,\, function$  $S\phi$ by $S\phi (t)=
\sum_{k=0}^{t-1} \phi (k)$, $S\phi (-t)= -\sum_{k=-t}^{-1} \phi
(k)$ if $t\ge 1$ and $S\phi (0)=0$.

Iterated primitives  and and iterated sums are defined by
$P^m\phi=P(P^{m-1}\phi)$  and $S^m\phi=S(S^{m-1}\phi)$, $m\in
\mathbf{N}$.

\smallskip

Weights are functions $w:G\rightarrow \mathbf{R}$ always assumed
to satisfy the following three conditions:

\bigskip

(1.1) \qquad $w$ is continuous, $w(t)\geq 1$ and $w(s+t)\leq w(s)w(t)$ for all $%
s,t\in G;$

\smallskip

(1.2) \qquad $w(-t)=w(t)$ for every $t\in G;$

\smallskip

(1.3) \qquad $\sum_{n=1}^{\infty }\frac{1}{n^{2}}\log w(nt)<\infty $ for every $%
t\in G.$

\smallskip

The symmetry condition (1.2) is only used to simplify the
exposition. Condition (1.3) is the Beurling-Domar condition (see
\cite{D}, [52, p. 132]). In the case that $w$ is bounded we will
assume $w=1,$ as this will cause no loss of generality.

\bigskip

A function $\phi :J\rightarrow X$ is called $w$-\textit{%
bounded} if $\phi /w$ is bounded. The spaces $B_{w}(J,X)$,
$BC_{w}(J,X)$ of all $w$-bounded  and continuous $w$-bounded
functions $\phi :J\rightarrow X$  respectively are  Banach spaces
with norm $\left\| \phi \right\| _{w,\infty }=
\displaystyle{\sup_{t\in J}}\,\,
\frac{%
\left\| \phi (t)\right\| }{w(t)}.$\ For this space and others,
when $w=1$ we will omit the subscript $w.$

\smallskip

Following [52, p.142] , we say that a function $%
\phi :J\rightarrow X$ is $w$-\textit{uniformly continuous} if
$\left\| \Delta _{h}\phi \right\| _{w,\infty }\rightarrow 0$ as
$h\rightarrow 0$ in $J$. The space of all $w$-uniformly continuous
functions is denoted $UC_{w}(J,X)$ and the closed subspace of
$BC_{w}(J,X)$ consisting of all $w$-uniformly continuous functions
is denoted $BUC_{w}(J,X)$. Of course if $J$ is discrete
$BUC_{w}(J,X)=BC_{w}(J,X)=B_{w}(J,X)$. We use also the notation $%
C_{w,0}(J,X)=\{w\xi :\xi \in C_{0}(J,X)\},$

\noindent  clearly
 a closed subspace of $%
BUC_{w}(J,X)$.

\smallskip

The space $L_{w}^{1}(G,X)$ consists of the strongly Haar-measurable
functions $f:G\rightarrow X$ for which $\left\| f\right\|
_{w,1}=\int_{G}\left\| f(t)\right\| w(t)d\mu (t)<\infty $ and the space $%
L_{w}^{\infty }(G,X)$ of the strongly Haar-measurable functions
$\phi :G\rightarrow X$ for which $\left\| \phi \right\| _{w,\infty
}=\,\, ess \displaystyle{\sup _{t\in G}}\,\,\,$ $\frac{\left\|
\phi (t)\right\| }{w(t)}<\infty$.

\smallskip

The space $L_{w}^{1}(G)=L_{w}^{1}(G,\mathbf{C})$ is a \textit{%
Beurling algebra} (see \cite{D}, [52, p. 83]). By (1.1) it is a subalgebra of the convolution algebra $%
L^{1}(G)$ and a Banach algebra under the norm $\left\| .\right\|
_{w,1}$ (see  [52, p.14] ). By (1.3), it is a Weiner algebra (see
[52, p.132]). Its Banach space dual is $L_{w}^{\infty
}(G)=L_{w}^{\infty }(G,\mathbf{C}).$

\smallskip

The Banach space dual of $X$ will be denoted $X^{\ast },$ the
space of bounded linear operators $T:X\rightarrow X$ by $L(X),$
the spectrum of a closed operator $T$\ on $X$ by $\sigma (T),$ its
resolvent set by $\rho (T)$, the spectrum of an element $x$ of a
unital Banach algebra by $\sigma (x)$ and its resolvent set by
$\rho (x).$

\smallskip

We use additive notation for $G$\ and multiplicative for the dual group $%
\widehat{G}.$ The \textit{Fourier transform} of $\ $a function\
$f\in L^{1}(G,X)$ is then given by $\widehat{f}(\gamma
)=\int_{G}\gamma (-t)f(t)d\mu (t),$ where $\gamma \in
\widehat{G}$. For $\phi \in BC_{w}(G,X)$ and $f\in L_{w}^{1}(G),$
or more generally $\phi \in BC_{w}(G,L(X))$ and $f\in
L_{w}^{1}(G,X),$ the \textit{convolution} $\phi
\ast f(t)=\int_{G}\phi (t-s)f(s)d\mu (s)$ is well-defined. Moreover, $%
\phi \ast f\in BUC_{w}(G,X).$

\bigskip

The structure of this paper is as follows. Examples and properties
of weights are given in section 2. In section 3 we develop the
theory of polynomials on groups. This is used in section 4 to
study the spaces $E_{w}(J,X)$ of $w$-ergodic functions. In
particular we obtain conditions on a subspace ${\F} $ of
$BC_{w}(J,X)$ under which a w -ergodic function belongs to ${\F} $
whenever its differences belong to ${\F}$. The space $AP_{w}(G,X)$
of $w$-almost periodic functions is studied in section 5. We
develop the tools of spectral analysis for unbounded functions in
section 6.

 \smallskip

The remainder of the paper deals with applications of the main
theorem. Firstly, in section 7 we study derivatives,  indefinite
integrals  and sums of functions on $\mathbf{R}$ and $\mathbf{Z}$.
In section 8 we prove our main result on general evolution
equations. Section 9 deals with the equation $B\phi =A\circ \phi
+\psi $ for differential operators $B.$\ In particular we give
conditions under which $\psi \in AP_{w}(\mathbf{R},X)$ implies
$\phi \in AP_{w}(\mathbf{R},X)$ thereby generalizing recent
results of Arendt and Batty. Similar results are obtained in
section 10 for convolution equations on $\mathbf{R}^{d}$ or
$\mathbf{Z}^{d}.$ In section 11 we obtain asymptotic stability for
$C_{0}$-semigroups $\{T(t):t\geq 0\}$ of operators on $X.$ In
particular we give conditions under which $\left\| T(t)\right\|
\leq
a(t)w(t) $ implies ${\displaystyle{\lim_{t\rightarrow \infty }}}\frac{\left\| T(t)x\right\| }{%
a(t)w(t)}=0$ for all $x\in X.$ As in Ph\'{o}ng \cite{PV} , the
function $a$ may have exponential growth. Finally, in section 12
we apply our results to recurrence equations on $\mathbf{Z.}$

\bigskip

2. \textbf{Weights}

Here we give some examples and make some general comments about
weights. Unless otherwise stated, {\textbf{in addition to
conditions (1.1) - (1.3) we will also assume our weights satisfy
the following}}:

\bigskip

(2.1) $\qquad \qquad w=1$ or $\frac{1}{w}\in C_{0}(G);$

\smallskip

(2.2) $\qquad \qquad \frac{\Delta _{h}w}{w}\in C_{0}(G)\,\,$ for
every $h\in G;$

\smallskip

(2.3) $\qquad \qquad{\displaystyle{\sup_{t\in G}}}\,\,\frac{\left|
\Delta _{h}w(t)\right| }{w(t)}\rightarrow 0$ as $h\rightarrow 0$
in $G.$

\smallskip

Occasionally we will also assume the existence of $N\in \mathbf{Z%
}_{+}$,  the non-negative integers,  such that

\smallskip

(2.4) $\qquad \qquad {\displaystyle{\lim_{\left| m\right|
\rightarrow \infty }}}\,\,\frac{w(mt)}{1+\left| m\right|
^{N+1}}=0$ for all $t\in G$ $;$ and

\smallskip

(2.5) $\qquad \qquad {\displaystyle{\inf_{m\in \mathbf{Z}}}}\,\,\frac{w(mt)}{\left| m\right| ^{N}}>0$ for some $%
t\in G.$

\bigskip

We will say that a weight $w$ has \textit{polynomial growth} of
order $N\in \mathbf{Z}_{+}$ if it satisfies (2.4), (2.5). The
Beurling-Domar condition (1.3) follows from (2.4). Condition (2.1)
is used in Proposition 5.4 and Lemma 8.8; (2.2) is used in Lemmas
4.4, 4.7, 4.11,  Theorem 4.12 and elsewhere; (2.3) is equivalent
to $w\in BUC_{w}(G,\mathbf{C});$ (2.4) and (2.5) are needed for
the spectral analysis in section 6.

\bigskip

\textbf{Example 2.1}.

\smallskip

 (a) The function $w_N (t)=(1+\left| t\right| )^{N}$, where $N\in \mathbf{N}$,
 is a weight on $\mathbf{R}^{d}$ or $\mathbf{Z}^{d}$ with polynomial growth
of order $N$.

\smallskip

(b) If $\alpha :G\rightarrow \mathbf{C}$ is a continuous additive function
with $\limsup_{t\rightarrow \infty }\left| \alpha (t)\right| =\infty $\ ,
then $w(t)=1+\left| \alpha (t)\right| $ is a weight with polynomial growth
of order $1.$

(c)  The function $\alpha(t_{1},t_{2})= t_{1}$ is a continuous
additive function on $\mathbf{R}^{2}$ or $\mathbf{Z}^{2}$ but
$\limsup_{t_2\rightarrow \infty }  \left| \alpha
(t_{1},t_{2})\right| = |t_1|\not =\infty $. Note $w
(t_1,t_2)=1+|t_1|$ is not a weight on $\mathbf{R}^{2}$ or
$\mathbf{Z}^{2}$ since it does not satisfy (2.1).

\smallskip

(d) Let $w=w_{1}w_{2}$ where $w_{1},w_{2}$ are two weights on $G$
with polynomial growth of orders $N_{1},N_{2}$ respectively. Then
$w$ is a weight with polynomial growth of order $N=N_{1}+ N_{2}$.

\smallskip

(e)  The function $w(t_{1},t_{2})=(1+\left| t_{1}\right|
)(1+\left| t_{2}\right| )$ is a weight on $\mathbf{R}^{2}$ or
$\mathbf{Z}^{2}$ with polynomial growth of order $2.$

\smallskip

(f) The functions $w(t)=\exp (\left| t\right| ^{p})$ and $%
w(t)=\exp ((1+\left| t\right| )^{p}),$ where $0\leq p<1,$ are weights on $%
\mathbf{R}^{d}$ or $\mathbf{Z}^{d}$ which are not of polynomial growth.

\smallskip

 (g) The function $ w(t)=\exp({1+|t|})$ does not satisfy (1.3),(2.2)
but it satisfies (2.3) for $G= \mathbf{Z}$ or $\mathbf{R}$.

\smallskip

 (h) Phong \cite{PV} considers weights on
$\mathbf{R}_{+}$ satisfying the condition $\displaystyle{\liminf_
{t\rightarrow \infty }}\,\, \frac{w(t+s)}{w(t)}\geq 1$ for all
$s\in \mathbf{R}_{+}.$ From this he concludes that
$w_{1}(s)={\displaystyle{\limsup_{t\rightarrow \infty
}}}\,\,\frac{w(t+s)}{w(t)}$ defines a non-decreasing weight on
$\mathbf{R}_{+}.$ The function $w(t)=(1+\left| \sin t\right|
)(1+\left| t\right| )^{N}$, where $N\in \mathbf{Z}_{+}$, is a weight on $%
\mathbf{R}^{d}$ or $\mathbf{Z}^{d}$ with polynomial growth of
order $N,$ not satisfying the condition used by Phong. Indeed,
$w_{1}(s)=1+\left| \sin s\right|$. Moreover,
${\displaystyle{\liminf_{\left| t\right| \rightarrow \infty
}}}\,\, \frac{w(t+\frac{\pi }{2})}{w(t)}=\frac{1}{2}$ and
${\displaystyle{\liminf_{\left| t\right| \rightarrow \infty}
}}\,\,\frac{w(t+\pi )}{w(t)}=1$ which contrasts with a statement
without proof in [49, after (4)].

\bigskip

\textbf{Proposition 2.2}. If $w$ satisfies (1.1) - (1.3) then
${\displaystyle{\lim_{n\rightarrow \infty }}}\,\,\frac{1}{n}\log
w(s+nt)=0$ uniformly with respect to $s,t\in K,$ for each compact
$K\subseteq G.$

\,

\begin{proof} It is enough to prove ${\displaystyle{\lim_{n\rightarrow \infty }}}\,\,\frac{1}{n}%
\log w(nt)=0$ uniformly on $K.$ It is well-known that the limit exists
pointwise. Indeed, for arbitrary positive integers $n,k$ choose $m,l\in
\mathbf{Z}$ with $n=mk+l$ and $0\leq l<k.$ By (1.1),(1.2) $w(nt)^{\frac{1}{n}%
}\leq w(kt)^{\frac{m}{n}}w(lt)^{\frac{1}{n}}\leq w(kt)^{\frac{1}{k}}w(lt)^{%
\frac{1}{n}}.$ Hence, ${\displaystyle{\limsup_{n\rightarrow \infty }}}\,\,w(nt)^{\frac{1}{n}%
}\leq w(kt)^{\frac{1}{k}}.$ Therefore ${\displaystyle{\limsup_{n\rightarrow \infty }}}\,\,w(nt)^{%
\frac{1}{n}}$exists and, by (1.3), ${\displaystyle{\limsup_{n\rightarrow \infty }}}\,\,w(nt)^{%
\frac{1}{n}}=0.$ Next note that the sequence $f_{m}(t)=w(2^{m}t)^{2^{-m}}$
is non-increasing and converges to $1$ for each $t\in G.$ By Dini's theorem
[7, p. 194] $(f_{m})$ converges uniformly to $1$ on $K.$ So for each $%
\varepsilon >0$ there exists $N(\varepsilon ,K)$ such that $%
w(2^{m}t)^{2^{-m}}<1+\varepsilon $ and $w(t)^{2^{-m/2}}<1+\varepsilon $\ for
all $t\in K$ and $m\geq N(\varepsilon ,K)$. For such $m$ set $n=2^{m},$ take
$1\leq k\leq n,$ $t\in K$ and the smallest integer $n^{\ast }\geq \sqrt{n}.$
We show that $w(kt)^{\frac{1}{n+k}}<(1+\varepsilon )^{2}.$ Firstly, if $%
1\leq k<n^{\ast }$ then $w(kt)^{\frac{1}{n+k}}\leq
w(t)^{\frac{k}{n+k}}\leq w(t)^{\frac{1}{n^{\ast }}}\leq
1+\varepsilon$. Secondly, if $n^{\ast }\leq
k<n$ then $k=dn^{\ast }+e$ with $0\leq d,e<n^{\ast }$ and so $w(kt)^{\frac{1%
}{n+k}}\leq w(n^{\ast }t)^{\frac{d}{n+k}}w(et)^{\frac{1}{n+k}}\leq w(n^{\ast
}t)^{\frac{1}{n^{\ast }}}w(et)^{\frac{1}{n+e}}\leq (1+\varepsilon )^{2}.$
Therefore, for $1\leq k\leq n,$ $w((n+k)t)^{\frac{1}{n+k}}\leq w(nt)^{\frac{1%
}{n+k}}w(kt)^{\frac{1}{n+k}}\leq (1+\varepsilon )^{4}.$ Thus, $w(pt)^{\frac{1%
}{p}}\leq (1+\varepsilon )^{4}$ for all $p\geq 2^{2N(\varepsilon
,K)}.$ \end{proof}

\bigskip

The following condition will be used in Proposition 5.4 and Lemma
8.8. By the next proposition, (2.6)\ is a consequence of (2.1).

\bigskip

(2.6)\qquad\qquad  For each compact neighborhood $H$ of $0$ there is a compact set $%
K\supseteq H$

\qquad \qquad \qquad \qquad such that ${\displaystyle{\sup_{t\in
K}}}\,\, w(t)\leq $ ${\displaystyle{\inf_{t\notin K}}}\,\, w(t).$

\bigskip

\textbf{Proposition 2.3} If $w$ is a continuous positive function on $G,$
then $\frac{1}{w}\in C_{0}(G)$\ if and only if (2.6) holds and $G$ is
compact or $w$ is unbounded.

\,

\begin{proof} Suppose $\frac{1}{w}\in C_{0}(G).$ Then for each $r>0$ the
sublevel set $K_{r}=\{t\in G:w(t)\leq r\}$ is compact. Since each compact
neighbourhood $H$ of $0$ is contained in some $K_{r},$ (2.6) is satisfied.
If $G$ is not compact, then $G$ $\backslash $ $K_{r}$ is non-empty for each $%
r>0$ and so $w$ is unbounded. Conversely, suppose (2.6) holds and $w$ is
unbounded. Given $\varepsilon >0,$ choose $t_{\varepsilon }\in G$ such that $%
w(t_{\varepsilon })>\frac{1}{\varepsilon }$ and a compact neighbourhood $H$
of $0$ and $t_{\varepsilon }.$ By (2.6) there is a compact set $K\supseteq H$
such that ${\displaystyle{\sup_{t\in K}}}\,\,w(t)\leq $
 ${\displaystyle{\inf_{t\notin K}}}\,\,w(t).$ Hence, $\frac{1}{%
w(t)}<\varepsilon $ for $t\notin K,$ proving $\frac{1}{w}\in C_{0}(G).$ If $%
G $ is compact, there is nothing to prove. \end{proof}

\bigskip

Now note that, if $w$ satisfies (1.1) and (1.2), then $\left|
\frac{\Delta _{h}w}{w}\right| \leq (w(h)-1)$ for all $h\in G$ and
hence
(2.3) would hold if $w(0)=1.$ Moreover, $\Delta _{h}(\frac{\phi }{w})=%
\frac{\Delta _{h}\phi }{w}-(\frac{\phi }{w})_{h}\frac{\Delta
_{h}w}{w}$ and therefore

\bigskip

(2.7) \qquad$\phi \in BUC_{w}(J,X)$ if and only if $\frac{\phi
}{w}$ is uniformly continuous and bounded.

\bigskip

Furthermore, $\left\| \phi _{t+h}-\phi _{t}\right\| _{w,\infty
}\leq w(t)\left\| \phi _{h}-\phi \right\| _{w,\infty }$ and so

\bigskip

(2.8) \qquad if $\phi \in BUC_{w}(J,X)$ then the function $t\to
\phi _{t}:J\to BUC_{w}(J,X)$ is continuous.

\bigskip

3. \textbf{Polynomials}

Following \cite{D} (see also \cite{BP}) we say that a function
$p\in C(J,X) $ is a \textit{polynomial} if $\Delta _{t}^{n+1}p=0$
for some $n\in \mathbf{N}$ and all $t\in J.$ Equivalently,
$p(s+mt)$ is a polynomial in $m\in \mathbf{Z}_{+}$\ of degree at
most $n$ for all $s,t\in J. $ For a non-zero polynomial $p,$ the
minimal such $n$ is its degree $\deg
(p) $. The space of all polynomials of degree at most $n$ is denoted $%
P^{n}(J,X). $\ Since $\Delta _{t}$ is a continuous mapping on $BC_{w}(J,X)$,
the polynomials in $BC_{w}(J,X)$\ form a closed subspace which we denote by $%
P_{w}(J,X).$ If $w$ has polynomial growth $N$, then $P_{w}(J,X)\subseteq
P^{N}(J,X).$ We also set $P_{w}^{n}(J,X)=P^{n}(J,X)\cap P_{w}(J,X).$

\bigskip

\textbf{Example 3.1}

 (a) If $J$ is $\mathbf{Z}$ or $\mathbf{R}$ then $%
P^{n}(J,X) $ is the space of\ ordinary polynomials.  Indeed, that
each (ordinary) polynomial is in $P^{n}(J,X)$ is clear.
Conversely, if $\ p\in P^{n}(\mathbf{Z},X)$ then $\Delta
_{1}^{n}p(m)=c$, a constant, and $\Delta
_{t}^{n}(p(m)-cm^{n}/n!)=0$ for all $t\in \mathbf{Z.}$ An
induction argument shows $\ p$ is an (ordinary) polynomial. If $\
p\in P^{n}(\mathbf{R},X)$
then $p|_{\mathbf{Z}}\in P^{n}(\mathbf{Z},X)$ and so $p|_{\mathbf{Z}}=q|_{%
\mathbf{Z}}$ for some (ordinary) polynomial $q:\mathbf{R\rightarrow }X.$ If $%
t=a/b$ where $a,b$ are non-zero integers, then $p(t+m/b)=q(t+m/b)$ for $m\in
-a+b\mathbf{Z.}$ These are polynomials in $m\in \mathbf{Z}_{+}$ and so agree
for all $m\in \mathbf{Z}_{+}$. In particular $p(t)=q(t)$ for all rationals $%
t $ and, by continuity, for all reals $t.$

(b) If $p:G\rightarrow \mathbf{C}$ is a continuous homomorphism,
then $p\in P^{1}(G,\mathbf{C).}$

(c) If $p\in P_{w}(G,X)$ and $f\in L_{w}^{1}(G),$ then $p\ast f\in
P_{w}(G,X).$ Indeed $\Delta _{t}^{n+1}(p\ast f)=(\Delta _{t}^{n+1}p)\ast f$
for all $t\in G.$

\bigskip

Occasionally it will be necessary for us to to assume that $%
P^{n}(J,\mathbf{C})$ is finite dimensional. That this condition is not
always satisfied is shown by the following example.

\bigskip

\textbf{Example 3.2. }Let $G=\{s:\mathbf{N\rightarrow Z}$ ; $\ s$
has finite support$\}$ with the discrete topology. So $G$ is
countable, locally compact,  $\sigma $-compact and not finitely
generated. Moreover the evaluation maps
$p_{n}:G\rightarrow \mathbf{C}$ defined by $p_{n}(s)=s(n)$ for $n\in \mathbf{%
N}$ are polynomials of degree $1$ and so $\dim
(P^{1}(G,\mathbf{C))=\infty } $. Note also that $J= \{s\in G: s(n)
=0 $ for $ n < 0\} $ satisfies the conditions of the following
proposition.

 \bigskip

\textbf{Proposition 3.3.} Assume that $J$ is a subsemigroup of $G$
with non-empty interior such that $G=J-J.$ Then  for each $n$ the restriction map $%
r:P^{n}(G,X)\rightarrow P^{n}(J,X)$ is a linear bijection. In particular,
every polynomial on $J$ has a unique extension to $G$.

\,

\begin{proof}  Firstly, let $p\in $ $P^{n}(G,X)$
be zero on $J$. For any $t\in G$ there are $u,v\in J$ with $t=u-v.$ Now $%
p(t+mv)=p(u+(m-1)v)$ is a polynomial in $m\in \mathbf{Z}_{+}$ which is zero
for all $m\geq 1.$ It is therefore zero for $m=0,$ showing $p(t)=0$ and $r$
is one-to-one. Secondly, let $q\in $ $P^{n}(J,X)$. We define an extension $p$
of $q$ to $G$ as follows. If $t\in G$ then $t=u-v$ for some $u,v\in J$ and
we set $p(t)=\sum_{j=0}^{n}(-1)^{j}\Delta _{v}^{j}q(u).$ If also $t=%
\widetilde{u}-\widetilde{v}$\ where $\widetilde{u},\widetilde{v}\in J$ we
must show

\smallskip
(3.1) $\qquad \qquad  \sum_{j=0}^{n}(-1)^{j}\Delta
_{v}^{j}q(u)=\sum_{j=0}^{n}(-1)^{j}\Delta _{\widetilde{v}}^{j}q(\widetilde{u}%
).$

\smallskip

To do this define a function $L:P^{n}(J,X)\rightarrow X$ by $%
L(q)=q(u)-\Delta _{v}q(\widetilde{u})=q(\widetilde{u})-\Delta _{\widetilde{v}%
}q(u).$ We prove by induction on $s\in \mathbf{Z}_{+}$ that

\smallskip

(3.2)  $ \qquad\sum_{j=0}^{n}(-1)^{j}\Delta
_{v}^{j}q(u)=\sum_{j=0}^{s-1}L(\Delta
_{v}^{j}\Delta _{\widetilde{v}}^{j}q)+\Delta _{v}^{s}\Delta _{\widetilde{v}%
}^{s}q(w)$ \,\,  if \, $q$ is of degree $n=2s,$

 \qquad\qquad \qquad\qquad and

 \qquad\qquad $\sum_{j=0}^{n}(-1)^{j}\Delta
_{v}^{j}q(u)=\sum_{j=0}^{s}L(\Delta _{v}^{j}\Delta
_{\widetilde{v}}^{j}q)$ \,\qquad  if \,\, $q$ is of degree
$n=2s+1$

\smallskip

\noindent where $w$ is an arbitrary element of $J.$ When $n=0,$ $q(u)=q(w)$ and when $%
n=1$, $\Delta _{v}q(u)=\Delta _{v}q(w)$ and so $q(u)-\Delta
_{v}q(u)=q(u)-\Delta _{v}q(\widetilde{u})=L(q).$ Hence, (3.2)
holds for $s=0.$ Assume (3.2) holds for polynomials of degree less
than $2s.$ If $n=2s,$ we can apply (3.2) to $\Delta _{v}q$ and
obtain

\qquad $\sum_{j=0}^{n}(-1)^{j}\Delta
_{v}^{j}q(u)=q(u)-\sum_{j=0}^{s-1}L(\Delta _{v}^{j}\Delta _{\widetilde{v}%
}^{j}\Delta _{v}q)$

\ \ \ \ \ \ \ \ \ \ \ \ \ \ \ \ \ \ \ \ \ \ \ \ \ \ \ \ \ \ \ \ \ $%
=q(u)+\sum_{j=0}^{s-1}\Delta _{\widetilde{v}}\Delta _{v}^{j}\Delta _{%
\widetilde{v}}^{j}\Delta _{v}q(u)-\sum_{j=0}^{s-1}\Delta _{v}^{j}\Delta _{%
\widetilde{v}}^{j}\Delta _{v}q(\widetilde{u})$

\ \ \ \ \ \ \ \ \ \ \ \ \ \ \ \ \ \ \ \ \ \ \ \ \ \ \ \ \ \ \ \ \ $%
=\sum_{j=0}^{s}\Delta _{v}^{j}\Delta _{\widetilde{v}}^{j}q(u)-%
\sum_{j=0}^{s-1}\Delta _{v}^{j}\Delta _{\widetilde{v}}^{j}\Delta _{v}q(%
\widetilde{u})$

\ \ \ \ \ \ \ \ \ \ \ \ \ \ \ \ \ \ \ \ \ \ \ \ \ \ \ \ \ \ \ \ \ $%
=\sum_{j=0}^{s-1}L(\Delta _{v}^{j}\Delta _{\widetilde{v}}^{j}q)+\Delta
_{v}^{s}\Delta _{\widetilde{v}}^{s}q(u).$

If $n=2s+1,$ we can apply this last result to $\Delta _{v}q$ and obtain

\qquad $\sum_{j=0}^{n}(-1)^{j}\Delta
_{v}^{j}q(u)=q(u)-\sum_{j=0}^{s-1}L(\Delta _{v}^{j}\Delta _{\widetilde{v}%
}^{j}\Delta _{v}q)-\Delta _{v}^{s}\Delta _{\widetilde{v}}^{s}\Delta _{v}q(w)$

\ \ \ \ \ \ \ \ \ \ \ \ \ \ \ \ \ \ \ \ \ \ \ \ \ \ \ \ \ \ \ \ \ $%
=q(u)+\sum_{j=0}^{s-1}\Delta _{\widetilde{v}}\Delta _{v}^{j}\Delta _{%
\widetilde{v}}^{j}\Delta _{v}q(u)-\sum_{j=0}^{s-1}\Delta _{v}^{j}\Delta _{%
\widetilde{v}}^{j}\Delta _{v}q(\widetilde{u})-\Delta _{v}^{s}\Delta _{%
\widetilde{v}}^{s}\Delta _{v}q(\widetilde{u})$

\ \ \ \ \ \ \ \ \ \ \ \ \ \ \ \ \ \ \ \ \ \ \ \ \ \ \ \ \ \ \ \ \ $%
=\sum_{j=0}^{s}\Delta _{v}^{j}\Delta _{\widetilde{v}}^{j}q(u)-\sum_{j=0}^{s}%
\Delta _{v}^{j}\Delta _{\widetilde{v}}^{j}\Delta _{v}q(\widetilde{u})$

\ \ \ \ \ \ \ \ \ \ \ \ \ \ \ \ \ \ \ \ \ \ \ \ \ \ \ \ \ \ \ \ \ $%
=\sum_{j=0}^{s}L(\Delta _{v}^{j}\Delta _{\widetilde{v}}^{j}q).$

Hence, (3.2) is proved and (3.1) follows, which means $p$ is
well-defined.

\bigskip

Note that if $t\in J$ then we can take $t=u,$ $v=0$ and so $p(t)=q(t).$ So $%
p $ is an extension of $q$. Next, we show $p$ is continuous. Since
$J$ has an interior point $s_{0}$, there is an open neighborhood
$W$ of $0$ in $G$
such that $s_{0}+W\subset J.$ Moreover, if $t=u-v,$where $u,v\in J$ then $t=%
\widetilde{u}-\widetilde{v},$where $\widetilde{u}=s_{0}+u,\widetilde{v}%
=s_{0}+v$. Let $(t_{\alpha })$ be a net in $G$ converging to $t.$ We may
suppose $t_{\alpha }=t+w_{\alpha }$ where $w_{\alpha }\in W.$ Setting $%
u_{\alpha }=\widetilde{u}+w_{\alpha }$ and $v_{\alpha }=\widetilde{v}$ we
find $u_{\alpha },v_{\alpha }\in J,t_{\alpha }=u_{\alpha }-v_{\alpha }$ and $%
(u_{\alpha })\rightarrow \widetilde{u}.$ So $p(t_{\alpha })\rightarrow p(t).$
Finally, if $t_{j}=u_{j}-v_{j}$ where $u_{j},v_{j}\in J$ and $m\in \mathbf{Z}%
_{+}$ then $p(t_{1}+mt_{2})=\sum_{j=0}^{n}(-1)^{j}\Delta
_{v_{1}+mv_{2}}^{j}q(u_{1}+mu_{2})$ which is a polynomial in $m$
of degree at most $n.$ So $p\in $ $P^{n}(G,X),$ proving $r$ is
onto. \end{proof}

\bigskip

\textbf{Proposition 3.4.} If $P^{n}(J,\mathbf{C})$ is finite dimensional
then $P^{n}(J,X)=P^{n}(J,\mathbf{C})$ $\otimes X.$ Similarly, if $%
P_{w}^{n}(J,\mathbf{C})$ is finite dimensional then $%
P_{w}^{n}(J,X)=P_{w}^{n}(J,\mathbf{C})$ $\otimes X.$

\,

\begin{proof} Clearly, $P^{n}(J,\mathbf{C})$ $\otimes X\subseteq $ $%
P^{n}(J,X).$ For the converse, which is clearly true when $n=0,$\ we use
induction on $n.$\ Let $\{p_{1},...,p_{k}\}$ be a basis of a complement $Q$\
of $P^{n-1}(J,\mathbf{C})$ in$\ P^{n}(J,\mathbf{C}).$ Since $\deg (p_{1})=n$
we can choose $t_{1}\in J$ such that $\Delta _{t_{1}}^{n}p_{1}\neq 0.$ Since
each $\Delta _{t_{1}}^{n}p_{j}$ is a constant we can set $q_{1}=p_{1}/\Delta
_{t_{1}}^{n}p_{1}$ and choose $\lambda _{1}\in \mathbf{C}$\ such that $%
\Delta _{t_{1}}^{n}(p_{2}-\lambda _{1}q_{1})=0.$ Set $q_{2}=p_{2}-\lambda
_{1}q_{1}.$ Then $\deg (q_{2})=n$ for otherwise $q_{2}\in Q\cap P^{n-1}(J,%
\mathbf{C})=\{0\},$ contradicting the linear independence of $p_{1},p_{2}.$
Hence we can choose $t_{2}\in J$ such that $\Delta _{t_{2}}^{n}q_{2}\neq 0.$
Continuing in this way, we obtain a basis $\{q_{1},...,q_{k}\}$ of $Q$ and a
subset $\{t_{1},...,t_{k}\}$ of $J$ such that $\Delta
_{t_{i}}^{n}q_{j}=\delta _{i,j}.$ Now let $p\in $ $P^{n}(J,X).$ For each $%
x^{\ast }\in X^{\ast }$\ we have $x^{\ast }\circ p\in P^{n}(J,\mathbf{C}).$
Hence $x^{\ast }\circ p=\sum_{j=1}^{k}q_{j}c_{j}(x^{\ast })+r (x^*)$ for some $%
c_{j} (x^*)\in \mathbf{C}$ and $r (x^*)\in P^{n-1}(J,\mathbf{C})$.
But $x^{\ast }\circ \Delta _{t_{i}}^{n}\,p=\Delta
_{t_{i}}^{n}(x^{\ast }\circ p)=c_{i}(x^{\ast })$ and so $c_i
(x^*)=  x^* (c_i)$ where $c_i = \Delta^n _{t_{i}} \, p$. Moreover,
$x^{\ast }\circ (p-\sum_{j=1}^{k}q_{j}c_{j})=r (x^*)\in
P^{n-1}(J,\mathbf{C})$  and so by the Hahn-Banach theorem
$\Delta^n_t (p-\sum_{j=1}^{k}q_{j}c_{j})=0$  for all $t\in J $  showing $%
p-\sum_{j=1}^{k}q_{j}c_{j}\in P^{n-1}(J,X).$ By the induction hypothesis $%
P^{n-1}(J,X)=P^{n-1}(J,\mathbf{C})$ $\otimes X$ and hence $p\in P^{n}(J,%
\mathbf{C})$ $\otimes X$ as required. The second assertion is
proved in the same way. \end{proof}

\textbf{Proposition 3.5. }If $P^{n}(J_{1},\mathbf{C})$ is finite
dimensional, then

 $ \qquad \qquad \qquad P^{n}(J_{1}\times J_{2},\mathbf{C})=\,
 \sum_{m=0}^{n}P^{m}(J_{1},\mathbf{C})\otimes
P^{n-m}(J_{2},\mathbf{C}).$

\begin{proof}\ The inclusion $\supseteq $ is clear. For the converse, take
any $p\in P^{n}(J_{1}\times J_{2},\mathbf{C}).$ As in the proof of the
previous proposition, we can find for each $m=1,...,n$ a basis $%
\{q_{1}^{m},...,q_{k_{m}}^{m}\}$ of a complement of $P^{m-1}(J_{1},\mathbf{C}%
)$ in$\ P^{m}(J_{1},\mathbf{C})$ and a subset $\{s_{1}^{m},...,s_{k_{m}}^{m}%
\}$ of $J_{1}$ such that $\ \Delta _{s_{i}^{m}}^{m}q_{j}^{m}=\delta _{i,j}.$
Also let $\left\{ q_{1}^{0}\right\} $ be a basis of $P^{0}(J_{1},\mathbf{C}%
). $ For any $t\in J_{2}$ we have $p(.,t)\in P^{n}(J_{1},\mathbf{C})$ and so
$p(s,t)=\sum_{m=0}^{n}\sum_{j=1}^{k_{m}}q_{j}^{m}(s)r_{j}^{m}(t)$ for some $%
r_{j}^{m}(t)\in \mathbf{C}$. We prove by backward induction on
$h$\ that $r_{j}^{h}\in P^{n-h}(J_{2},\mathbf{C}).$ Now

$\Delta
_{(s_{i}^{n},0)}^{n}p(s,t)=\sum_{m=1}^{n}\sum_{j=1}^{k_{m}}\Delta
_{s_{i}}^{n}q_{j}^{m}(s)r_{j}^{m}(t)=r_{i}^{n}(t)$

\noindent so each $r_{i}^{n}$ is a constant as required. So
suppose each $ r_{j}^{m}\in P^{n-m}(J_{2},\mathbf{C})$ for
 $n\geq m\geq h+1\,\,\,$ and
$\,\,\, 1\leq
j\leq k_{m}$. Then

$p(s,t)-\sum_{m=h+1}^{n}%
\sum_{j=1}^{k_{m}}q_{j}^{m}(s)r_{j}^{m}(t)=\sum_{m=0}^{h}\,
\sum_{j=1}^{k_{m}}q_{j}^{m}(s)r_{j}^{m}(t)\, \in \, P^{n}(J_{1}\times J_{2},%
\mathbf{C})$ and

 $\Delta
_{(s_{i}^{h},0)}^{h}\sum_{m=0}^{h}%
\sum_{j=1}^{k_{m}}q_{j}^{m}(s)r_{j}^{m}(t)=r_{i}^{h}(t).$ So each $%
r_{i}^{h}\in P^{n-m}(J_{2},\mathbf{C})$

\noindent and the proposition is proved. \end{proof}

\bigskip Recall that a group is called a $torsion\,\, group$
if every element has finite order. If the orders of the elements
are bounded the group is said to be of $bounded\,\, order$.

\bigskip\textbf{Proposition 3.6.} Let $H$ be a closed subgroup of $G$
such that $G/H$ is a torsion group.

(a) The restriction map $r:$ $P^{n}(G,X)\rightarrow P^{n}(H,X)$ is
one-to-one and so $\dim $ ($P^{n}(G,X))\leq $ $\dim (P^{n}(H,X)).$

(b) If also $G/H$ is of bounded order and $P^{n}(H,\mathbf{C})$ is finite
dimensional, then

$r:$ $P^{n}(G,X)\rightarrow P^{n}(H,X)$ is a linear isomorphism.

\,

\begin{proof}\ (a) Let $p\in P^{n}(G,X)$ satisfy $p(t)=0$ for all $t\in H.$%
\ Let $t\in G\backslash H$ and let $\pi :G\rightarrow G/H$ be the quotient
map. Since $G/H$ is a torsion group, meaning\ every element has finite
order, $\pi (kt)=0$ or $kt\in H$ for some $k\in \mathbf{N.}$ Hence $p(mkt)=0$
for all $m\in \mathbf{N.}$ But $p(mt)$ is a polynomial in $m\in \mathbf{Z}%
_{+}$ and so is zero. In particular, $p(t)=0$ showing $r$ is one-to-one.

(b) Let $\{p_{1},...,p_{m}\}$ be a basis of $P^{n}(H,\mathbf{C})$
and choose $k\in \mathbf{N}$ such that\ $kt\in H$ for all $t\in G.$ Define $%
q_{j}(t)=p_{j}(kt)$ and suppose $\sum_{j=1}^{m}\alpha _{j}q_{j}=0$ on $G$
for some $\alpha _{j}\in \mathbf{C.}$ Then $\sum_{j=1}^{m}\alpha _{j}p_{j}=0$
on $kG.$ But $kG$ is a closed subgroup of $H$ such that $H/kG$ is a torsion
group. By part (a), $\sum_{j=1}^{m}\alpha _{j}p_{j}=0$ on $H.$ Hence each $%
\alpha _{j}=0,$ showing $\{q_{1},...,q_{m}\}$ is linearly independent and $%
\dim $ ($P^{n}(G,\mathbf{C}))\geq $ $\dim (P^{n}(H,\mathbf{C}))$. Therefore $%
r:$ $P^{n}(G,X)\rightarrow P^{n}(H,X)$ is a linear isomorphism, by (a) when $%
X=\mathbf{C}$ and then by Proposition 3.4 for general $X.$
\end{proof}

\bigskip

\textbf{Proposition 3.7. } Let $G$ be compactly generated.

(a) $G$  has  an open subgroup $G_0 \cong \mathbf{R}^{d}\times
\mathbf{Z}^{m}\times K$ for some integers $d, m \ge 0$ and some
compact group $K$ with $G/G_0$  a finite group. Moreover, there is
a finite subgroup $ E\subset G_0^{\perp}  $ such that the
connected component $\Gamma_0$ containing $1$ of $\widehat{G}$
satisfies
 $ \Gamma_0/E \cong \mathbf{R}^{d}\times \mathbf{T}^{m}$.

 (b) $P^{n}(G,\mathbf{C})$ is finite dimensional.

(c) For each $p\in P^{n}(G,X)$ there exist $p_{j}\in P^{n}(G,X)$ and $%
q_{j}\in P^{n}(G,\mathbf{C})$\ with $q_{j}(0)=0$ such that $\Delta
_{h}p(t)=\sum_{j=1}^{k}p_{j}(t)q_{j}(h)$\ for all $h,t\in G$.

(d) $P_{w}(G,X)\subseteq BUC_{w}(G,X).$

(e) If also $w$ has polynomial growth $N$, then
$P_w(G,\mathbf{C})$ is finite dimensional.

\begin{proof} (a) By the principal structure theorem for locally
compact abelian groups (see [53, Theorem 2.4.1]), each such group
$G$ has an open subgroup $G_1\cong \mathbf{R}^{d}\times K$ for
some $d\geq 0$ and some compact group $K.$ Hence $G/G_{1}$ is
discrete. Since $G$ is compactly generated  so is $G/G_{1}$ and
therefore  $G/G_{1}$ is finitely generated. So $G/G_{1}\cong
\mathbf{Z}^{m}\times F$ for some $m\geq 0$ and some finite group
$F$ (see Theorem 9.3 in \cite{F}). As  $\mathbf{Z}^{m}$ is a free
group,  $G$ has a subgroup $H$ isomorphic to $\mathbf{Z}^{m}$ with
$H \cap G_1 =\{0\}$. Since $G_1$ is open, $G_0=G_1 \times H$ is
open too. Moreover, $G/G_{0}\cong F$. By [53, Theorem 2.1.2,
2.2.2], $\widehat{G}/G_0^{\perp}\cong \widehat{G_0}$ and
$\widehat{G_0}\cong \mathbf{R}^{d}\times \mathbf{T}^{m}\times
\widehat {K}$. Since $\widehat {K}$ is discrete, the connected
component $\Omega_0$ of $\widehat{G_0}$ is isomorphic to
$\mathbf{R}^{d}\times \mathbf{T}^{m}$. As  $G_0^{\perp}$ is
finite, $E= G_0^{\perp}\cap \Omega_0$ is finite too. Since the
natural homomorphism of $\widehat {G}$ onto
$\widehat{G}/G_0^{\perp}$ is a continuous open map it follows
$\Gamma_0/E \cong \Omega_0 \cong \mathbf{R}^{d}\times
\mathbf{T}^{m}$.

(b) By Proposition 3.6, $\dim $ $P^{n}(G,\mathbf{C})=$ $\dim $ $P^{n}(%
\mathbf{R}^{d}\times \mathbf{Z}^{m}\times K,\mathbf{C})$. But $P^{n}(\mathbf{%
R}^{d}\times \mathbf{Z}^{m}\times
K,\mathbf{C})=P^{n}(\mathbf{R}^{d}\times
\mathbf{Z}^{m},\mathbf{C})$ which, by example 3.1(a) and
Proposition 3.5, is finite dimensional. This proves (b).

(c) Since $r: P^{n}(G,X)\rightarrow
P^{n}(G_{0},X)$ is a linear isomorphism and  (c) holds when $G=\mathbf{R%
}^{d}\times \mathbf{Z}^{m}$ it holds for compactly generated $G.$

(d) If $p\in P_{w}(G,X)$ choose $p_{j},q_{j}$ as in (c).\ Since
$\left\| \Delta
_{h}p(t)\right\| \leq cw(t)\sum_{j=1}^{k}\left\| q_{j}(h)\right\| ,$ where $%
c=$ sup$_{j}\left\| p_{j}\right\| _{w,\infty },$ it follows that
$p\in BUC_{w}(G,X),$ proving (d).

(e) The assumptions imply $P_{w}(G,X) \subset P^N(G,X)$ and
therefore the result follows from (b).
 \end{proof}

 \bigskip

\textbf{Remark 3.8.} The above proof also shows that $P^{n}(G,\mathbf{C})$
is finite dimensional whenever $G$ has an open subgroup $G_{1}$ such that $%
G_{1}\cong \mathbf{R}^{d}\times K$ for some $d\geq 0$ and some
compact group $K$ with $G/G_{1}\cong \mathbf{Z}^{m}\times F$ for
some $m\geq 0$ and some torsion group $F.$

\bigskip

\textbf{4}. $w$-\textbf{Ergodic functions}

As in \cite{M}, \cite{MW},  \cite{BP} we say that a function $%
\phi :J\rightarrow X$ is \textit{ergodic} if $\phi \in BC(J,X)$
and there exists $M_{\phi }\in X$ such that for each $\varepsilon
>0$ there
is a finite subset $F\subseteq J$ with $|| \frac{1}{|F|}\,
{\displaystyle{\sum_{t\in F}}}\,\,(\phi _{t}-M_{\phi })|| <
\varepsilon $. Below we will use the abreviation $R_{F}\phi
=\frac{1}{\left| F\right| }{\displaystyle{\sum_{t\in F}}}\,\,\phi
_{t}.$ The element $M_{\phi },$ clearly unique, is the (Maak)
\textit{mean} of $\phi $. Denoting the space of all such ergodic
functions by $E(J,X),$ and by $M:E(J,X)\rightarrow X$ the function
$M(\phi )=M_{\phi },$ it follows that $M$\ is linear and
continuous. Moreover, if $E^{0}(J,X)=\{\phi \in E(J,X):M(\phi
)=0\}$ then $E(J,X)=E^{0}(J,X)\oplus X.$

\bigskip

We shall say a function $\phi :J\rightarrow X$ is $w$-\textit{%
ergodic}\ if $\phi \in BC_{w}(J,X)$\ and there is a polynomial
$p\in P_{w}(J,X)$ such that $(\phi -p)/w$ is ergodic with mean
$0.$\ The
polynomial $p,$\ not necessarily unique, is called a $w$-\textit{mean} of $%
\phi $ and the space of all $w$-ergodic functions is denoted$\
E_{w}(J,X). $ The subspace $E_{w}^{0}(J,X)$\ of functions with
$w$-mean $0$\
is therefore a closed subspace of $BC_{w}(J,X)$\ and $E_{w}(J,X)=$\ $%
E_{w}^{0}(J,X) + P_{w}(J,X).$ A function $\phi $ is \textit{totally }$w$%
\textit{-ergodic} if $\gamma \phi $ is $w$-ergodic for all $\gamma
\in \widehat{G}.$

\bigskip

If $\,\, P_{w}(J,\mathbf{C})\,\,$ is finite \, dimensional, it follows \, that $\,\, E_{w}(J,%
\mathbf{C})\,\,$ is a closed \, subspace of

\noindent $BC_{w}(J,\mathbf{C})$. Moreover, we can choose a
subspace $P_{w}^{M}(J,\mathbf{C})$ of $P_{w}(J,\mathbf{C})$
such that $\,E_{w}(J,\mathbf{C})$ $=E_{w}^{0}(J,\mathbf{C})\oplus P_{w}^{M}(J,%
\mathbf{C})$. The (continuous) projection map $M_{w}:E_{w}(J,\mathbf{C}%
)\rightarrow P_{w}^{M}(J,\mathbf{C})$ then provides a unique $w$-mean $%
M_{w}(\phi )$ for each $\phi \in E_{w}(J,\mathbf{C}).$ Now set $%
P_{w}^{M}(J,X)$\,
$=P_{w}^{M}(J,\mathbf{C})\otimes X$ and define $%
M_{w}:E_{w}(J,X)\rightarrow P_{w}^{M}(J,X)$ by $M_{w}(\phi
)=\sum_{j=1}^{k}M_{w}(p_{j})\otimes x_{j}$ where $\phi \in
E_{w}(J,X)$
has $w$-mean $p$ $=\sum_{j=1}^{k}p_{j}\otimes x_{j}\in P_{w}(J,\mathbf{C}%
)\otimes X.$

\bigskip

\textbf{Proposition 4.1.} Assume $P_{w}(J,\mathbf{C})$ is finite
dimensional. Then the map

\smallskip

\qquad  \qquad \qquad \qquad \qquad $M_{w}:E_{w}(J,X)\to
P_{w}^{M}(J,X)$

\smallskip

\noindent  is
well-defined and continuous. Moreover, for each $\phi \in E_{w}(J,X),$ $%
M_{w}(\phi )$ is a $w$-mean for $\phi $ and for each of its
translates. Finally, $E_{w}(J,X)$ is a closed translation
invariant subspace of $BC_{w}(J,X)$ and $E_{w}(J,X)=$\
$E_{w}^{0}(J,X)\oplus P_{w}^{M}(J,X).$

\,

\begin{proof} Let $\phi \in E_{w}(J,X)$ have means $p$ $%
=\sum_{j=1}^{k}p_{j}\otimes x_{j}$ and $q=\sum_{j=1}^{m}q_{j}\otimes y_{j}.$
Then $p-q\in E_{w}^{0}(J,X)$ and so $x^{\ast }\circ (p-q)\in E_{w}^{0}(J,%
\mathbf{C})$ for all $x^{\ast }\in X^{\ast }.$ Hence $M_{w}(x^{\ast }\circ
(p-q))=0=x^{\ast }\circ (\sum_{j=1}^{k}M_{w}(p_{j})\otimes
x_{j}-\sum_{j=1}^{m}M_{w}(q_{j})\otimes y_{j})$ which gives $%
\sum_{j=1}^{k}M_{w}(p_{j})\otimes x_{j}=\sum_{j=1}^{m}M_{w}(q_{j})\otimes
y_{j}$ showing $M$ is well-defined. Also, $p_{j}-M_{w}(p_{j})\in E_{w}^{0}(J,%
\mathbf{C})$ and so by Lemma 4.2 below $p-M_{w}(p)\in
E_{w}^{0}(J,X).$ Hence $M_{w}(\phi )$ is a mean for $\phi$.
Moreover,  $\left\| x^{\ast }\circ M_{w}(\phi )\right\| =\left\|
M_{w}(x^{\ast }\circ \phi )\right\| \leq c\left\| x^{\ast }\circ
\phi \right\| =c\sup_{t\in J}\left\| x^{\ast }\circ \phi
(t)\right\| /w(t)\leq c\left\| x^{\ast }\right\| \left\| \phi
\right\|$. Hence,\ $|||M_{w}(\phi )|||\leq c\left\| \phi \right\|
$ where $|||\psi |||=\sup \frac{\left\| x^{\ast }\circ \psi
\right\| }{\left\| x^{\ast }\right\| }$ for $\psi \in
P_{w}(J,\mathbf{C})\otimes X.$ By Lemma 4.3 below, $M_{w}$ is
continuous. If
$(\phi _{n})$ is a sequence in $E_{w}(J,X)$ converging to $\phi $ in $%
BC_{w}(J,X),$ let $p_{n}=M_{w}(\phi _{n}).$ Then $(p_{n})$\
converges to some $p\in P_{w}^{M}(J,X)$ and so $(\frac{\phi
_{n}-p_{n}}{w})$ converges to $\frac{\phi -p}{w}$ in $BC(J,X).$\
By the continuity of the Maak mean
function, $M(\frac{\phi -p}{w})=0$ and so $E_{w}(J,X)$ is closed. That $%
E_{w}(J,X)=$\ $E_{w}^{0}(J,X)+P_{w}^{M}(J,X)$ is clear and that the sum is
direct follows from the Hahn-Banach theorem. Finally, for each $t\in J$ we
have $\frac{\phi _{t}-p}{w}=\frac{\Delta _{t}\phi }{w}+\frac{\phi -p%
}{w}$ and so, by Lemma 4.4 below, $p$ is a $w$-mean of $\phi _{t}$ and $%
E_{w}(J,X)$ is translation invariant. \end{proof}

 \bigskip

\textbf{Lemma 4.2.} Assume $P_{w}(J,\mathbf{C})$ is finite dimensional. If $%
p\in P_{w}(J,X)$ and $x^{\ast }\circ p\in E_{w}^{0}(J,\mathbf{C)}$ for each $%
x^{\ast }\in X^{\ast }$ then $p\in E_{w}^{0}(J,X\mathbf{).}$

\,

\begin{proof} By Proposition 3.4 we can choose $q_{1},...,q_{m}\in P_{w}(J,%
\mathbf{C)}$ and linearly independent unit vectors $x_{1},...,x_{m}\in X$
such that $p=\sum_{j=1}^{m}q_{j}\otimes x_{j}.$ Then choose unit vectors $%
x_{j}^{\ast }\in X^{\ast }$ such that $\left\langle x_{j}^{\ast
},x_{i}\right\rangle =\delta _{i,j}.$ Given $\varepsilon >0$ there
are finite subsets $F_{j}$ of $J$ such that $\left\|
R_{Fj}(x_{j}^{\ast }\circ p/w)\right\| <\varepsilon /m.$ Setting
$F=F_{1}+...+F_{m}$ we find

$\left\| R_{F}(p/w)\right\| =\left\|
\sum_{j=1}^{m}R_{F}(q_{j}/w)\otimes x_{j}\right\| \leq
\sum_{j=1}^{m}\left\| R_{F}(q_{j}/w)\right\|$

$=\sum_{j=1}^{m}\left\| R_{F}(x_{j}^{\ast }\circ p/w)\right\| \leq
\sum_{j=1}^{m}\left\| R_{F_{j}}(x_{j}^{\ast }\circ p/w)\right\|
<\varepsilon$.

\noindent This proves that $p\in E_{w}^{0}(J,X\mathbf{).}$
\end{proof}

\bigskip

\textbf{Lemma 4.3.} If $P=P_{w}(J,\mathbf{C})$ is a finite
dimensional Banach space, then on $P\otimes X$ the natural norm
and the norm $|||\phi |||=\sup \frac{\left\| x^{\ast }\circ \phi
\right\| }{\left\| x^{\ast }\right\| }$ are equivalent.

\,

\begin{proof} Let $\{p_{1},...,p_{k}\}$ be a basis of $P$
consisting of unit vectors. If $p=\sum_{j=1}^{k}c_{j}p_{j}$, where
$c_{j}\in \mathbf{C,}$ then $\left\| p\right\| \sim
\sum_{j=1}^{k}\left\| c_{j}\right\| ,$ $\,\, \sim $ denoting
equivalence of norms. Also, every $\phi \in P\otimes X$ has a
unique representation $\phi =\sum_{j=1}^{k}p_{j}\otimes
x_{j},$ where $x_{j}\in X,$ and $\left\| \phi \right\| \sim %
\sum_{j=1}^{k}\left\| x_{j}\right\|$. Hence, $|||\phi |||=\sup \frac{%
\left\| \sum_{j=1}^{k}p_{j}\left\langle x^{\ast },x_{j}\right\rangle
\right\| }{\left\| x^{\ast }\right\| }\leq \sum_{j=1}^{k}\left\|
x_{j}\right\| \sim \left\| \phi \right\|$. Conversely, choose $%
j_{0}$ such that $\left\| x_{j}\right\| \leq \left\|
x_{j_{0}}\right\| $ for each $j.$ Then choose $x^{\ast }\in
X^{\ast }$ such that $\left\langle x^{\ast
},x_{j_{0}}\right\rangle =\left\| x_{j_{0}}\right\| $ and $\left\|
x^{\ast }\right\| =1.$ Hence, $|||\phi |||\geq \left\|
\sum_{j=1}^{k}p_{j}\left\langle x^{\ast },x_{j}\right\rangle
\right\| \sim\sum_{j=1}^{k}\left| \left\langle x^{\ast
},x_{j}\right\rangle \right| \geq \left\| x_{j_{0}}\right\| \geq
\frac{1}{k}\sum_{j=1}^{k}\left\| x_{j}\right\| \sim \left\| \phi
\right\|$. \end{proof}

 \bigskip

\textbf{Lemma 4.4.} If $\phi \in BC_{w}(J,X)$ then $\Delta
_{t}\phi \in E_{w}^{0}(J,X)$ for all $t\in J.$

\,

\begin{proof} We have $\frac{\Delta _{t}\phi }{w}=\Delta _{t}(\frac{%
\phi }{w})+(\frac{\phi }{w})_{t}\frac{\Delta _{t}w}{w}$ where
$\Delta _{t}(\frac{\phi }{w})\in E^{0}(J,X)$ by [16, Proposition
2.2]\ and $(\frac{\phi }{w})_{t}\frac{\Delta _{t}w}{w}\in
E^{0}(J,X) $ by condition (2.2). \end{proof}

\bigskip

Even on $\mathbf{R}$, not every continuous function has bounded
differences. However, uniformly continuous functions on
$\mathbf{R}$ have bounded differences. To study this behavior we
denote by

\smallskip

 $UCD_w(J,X)=\{
\phi\in UC_w(J,X): \Delta_h \phi \in BUC_w(J,X), h\in J\}$.

 \noindent See [15,
(4.1)] where it is shown that if $G$ is connected and $w=1$,
$UCD(G,X) =UC(G,X)$. On the other hand the function $f(x,n) =
x+n^2$ defined on $G= \mathbf {R}\times \mathbf {Z}$ is uniformly
continuous with differences  not all bounded. If $J\in
\{\mathbf{R_+,R}\}$, $UCD(J,X)=UC(J,X)$.  Moreover, $UC(J,X)=\{
\phi\in L^1_{loc}(J,X): \Delta_h \phi \in BUC(J,X), h\in J\}$ (see
[11, Corollary 5.5]).

\bigskip

\textbf{Lemma 4.5.} (a) Let $\phi \in UCD(J,X)$ and $h\in J$. The
following are equivalent

(i) $\Delta _{h}\phi \in E (J,X)$.

(ii)  $\Delta _{nh}\phi \in E (J,X)$ for some $n\in \mathbf{N}$.

(iii) $\Delta _{nh}\phi \in E (J,X)$ for all $n\in \mathbf{N}$.

\smallskip

(b) Let $\phi \in BC(J,X)$ and $h\in J$. The following are
equivalent

(i) $\phi \in E (J,X)$

(ii)  $\sum_{k=0}^{n} \phi_{kh} \in E (J,X)$ for some $n\in
\mathbf{N}$.

(iii)  $\sum_{k=0}^{n} \phi_{kh} \in E (J,X)$ for all $n\in
\mathbf{N}$.

\,

\begin{proof} (a) Assume $\Delta _{h}\phi \in E (J,X)$. Since $ E
(J,X)$ is  linear and translation invariant, the identity

\smallskip

(4.1) \qquad $\Delta_{nh}\phi =\sum_{k=0}^{n-1}
(\Delta_{h}\phi)_{kh}$

\smallskip

\noindent implies $\Delta_{nh}\phi\in E(J,X)$ for all $n\in
\mathbf{N}$. This proves (i) implies (ii) and (iii).

 Assume $\Delta_{nh}\phi\in E(J,X)$ for some $n\in
\mathbf{N}$. Using the definition of Maak ergodicity, for each
$\varepsilon
>0$ there exists $t_1, \cdots ,t_m \in G $ such that

$||(1/m)\sum_{k=1}^m (\Delta _{nh}\phi)_{t_k}-M(\Delta
_{nh}\phi)||\le \varepsilon$.

The identity (4.1) gives

 $||(1/(mn))\sum_{k=1}^m
(\sum_{j=0}^{n-1} \Delta_{h}\phi_{jh}) _{t_k}-M(\Delta
_{nh}\phi)/n||\le \varepsilon/n$.

\noindent  Again applying the definition of Maak ergodicity  to
$\Delta_h \phi$, one gets   $\Delta_{h}\phi \in E(J,X)$ with mean
$M(\Delta _{nh}\phi)/n$. This proves (ii) or (iii) implies (i).

(b) Follows similarly as in (a). \end{proof}

\bigskip

\textbf{Corollary 4.6}. Let $J\in\{\mathbf{R}_+, \mathbf{R}\}$,
$\phi \in UC(J,X)$ and $\Delta _{h}\phi \in E (J,X)$ for some
$h>0$. Then $\Delta _{s}\phi \in E (J,X)$ for all $ s >0$.

\,

\begin{proof} By Lemma 4.5 (a), one gets $\Delta_{hn}\phi$,
 $\Delta_{h/n}\phi \in E(J,X)$ for all $n\in \mathbf{N}$. It follows
$\Delta_{rh}\phi \in E(J,X)$ for all $r\in \mathbf{Q}^+$. Since
$\{rh: r\in \mathbf{Q} \}$ is dense in $\mathbf{R}$, it follows
that for each $s_0 > 0$ there exists a sequence $(r_n) \subset
\mathbf{Q}^+$ such that $r_n h \to s_0$ as $n\to \infty$. Since
$\Delta_s \phi \in BUC (J,X)$ for all $s
>0$, one gets  $\Delta_{r_n h} \phi \to \Delta_{s_0} \phi$
as $n\to \infty$. This implies $\Delta _{s}\phi \in E (J,X)$ for
all $s
>0$. \end{proof}

\bigskip

\textbf{Lemma 4.7}. Let $\phi \in UCD_{w}(J,X)$ and $h\in J$.  The
following are equivalent

(i) $\Delta _{h}\phi \in E_w (J,X)$

(ii)  $\Delta _{nh}\phi \in E_w (J,X)$ for some $n\in \mathbf{N}$.

(iii) $\Delta _{nh}\phi \in E_w (J,X)$ for all $n\in \mathbf{N}$.

\begin {proof} Assume $\Delta _{h}\phi \in E_w (J,X)$. Then $(\Delta _{h}\phi)/w\in E (J,X)$. Since $ E (J,X)$
is linear, translation invariant and $(\Delta_s w)/w \in C_0
(J,X)\subset E(J,X)$ for all $s \in J$, the identity

\smallskip

(4.2) \qquad $\frac{\Delta_{nh}\phi}{w} =\sum_{k=0}^{n-1}
(\frac{\Delta_{h}\phi}{w})_{kh}+ \sum_{k=1}^{n-1}
(\frac{\Delta_{h}\phi}{w})_{kh}\cdot \frac {\Delta_{kh} w}{w}$.

\smallskip

\noindent implies $\Delta_{nh}\phi\in E_w(J,X)$ for all $n\in
\mathbf{N}$.

Conversely, assume $\Delta_{nh}\phi\in E_w(J,X)$ for some $n\in
\mathbf{N}$.  By (2.2), $(\Delta_s w)/w \in C_0 (J,X)$ for all $s
\in J$. Therefore (4.2) implies $\sum_{k=0}^{n-1}
(\frac{\Delta_{h}\phi}{w})_{kh} \in E(J,X)$. By Lemma 4.5 (b), one
gets $\frac{\Delta_{h}\phi}{w} \in E(J,X)$. This implies
$\Delta_{h}\phi \in E_w(J,X)$. \end{proof}

\bigskip

\textbf{Corollary 4.8}. Let $J\in \{\mathbf{R}_+, \mathbf{R}\}$,
$\phi \in UC_{w}(J,X)$ and $\Delta _{h}\phi \in E_{w}(J,X)$ for
some $h>0$. Then $\Delta _{s}\phi \in E_{w}(J,X)$ for all $s
>0$.

\begin{proof} This follows similarly as in Corollary 4.6 using Lemma 4.7
instead of Lemma 4.5. \end{proof}

\bigskip

\textbf{Proposition 4.9}. If $\phi \in E_{w}(G,X)$ has $w$%
-mean $p,$ then $\phi |_{J}\in E_{w}(J,X)$ and $\phi |_{J}$ has $w$%
-mean $p|_{J}.$

\,

\begin{proof} Given $\varepsilon >0$ there is a finite subset $%
F=\{t_{1},...,t_{m}\}\subseteq G$ such that

$\left\| \frac{1}{m}%
\sum_{j=1}^{m}(\frac{\phi -p}{w})(t_{j}+t)\right\| <\varepsilon $
for all $t\in G.$

\noindent  Choose $u_{j},v_{j}\in J$ such that $t_{j}=u_{j}-v_{j}.$ Let $%
v=v_{1}+...+v_{m}$ and set $s_{j}=t_{j}+v.$ So $s_{j}\in J$ and
$\left\| \frac{1}{m}\sum_{j=1}^{m}(\frac{\phi
-p}{w})(s_{j}+t)\right\| <\varepsilon $ for all $t\in J.$
\end{proof}

\bigskip

\textbf{Proposition 4.10.} Let $\phi \in BUC_{w}(G,X),$ $f\in
L_{w}^{1}(G)$
and suppose $\phi |_{J}$ is $w$-ergodic with $w$-mean $p|_{J}$ where $%
p\in P_{w}(G,X).$\ If $p=0$\ or $G$ is compactly generated, then
$(\phi \ast f)|_{J}$ is $w$-ergodic with $w$-mean $(p\ast
f)|_{J}$.

\,

\begin{proof} If $\chi $ is the characteristic function of a compact set $%
K\subseteq G$ then $(\phi -p)\ast \chi (s)=\int_{-K}(\phi
-p)_{t}(s)d\mu (t)$ for each $s\in G.$ But for each $t\in G,$\
$(\phi -p)_{t}=\Delta _{t}(\phi -p)+(\phi -p)$ and so by lemma
4.4, $(\phi -p)_{t}|_{J}\in $ $E_{w}^{0}(J,X).$\ Also, by
Proposition 3.7(d), $\phi -p\in BUC_{w}(G,X).$ By (2.8), the
function $t\mapsto (\phi -p)_{t}|_{J}:G\rightarrow E_{w}^{0}(J,X)$
is continuous and hence weakly measurable and separably-valued on
$-K.$ The integral $\int_{-K}(\phi -p)_{t}|_{J}d\mu (t)$ is
therefore a convergent Haar-Bochner integral and so
belongs to $E_{w}^{0}(J,X).$ As evaluation at $s\in J$ is continuous on $%
E_{w}^{0}(J,X)$ we conclude that $((\phi -p)\ast \chi )|_{J}\in
E_{w}^{0}(J,X)$. Hence also $((\phi -p)\ast \sigma )|_{J}\in
E_{w}^{0}(J,X)$ for any step function $\sigma :G\rightarrow
\mathbf{C.}$ By [52, p.83]\ the step functions are dense in $%
L_{w}^{1}(G)$ and so $((\phi -p)\ast f)|_{J}\in E_{w}^{0}(J,X)$\vspace{1pt%
} for any $f\in L_{w}^{1}(G).$ \end{proof}

 \bigskip

The difference theorem below, included here in order to characterize $%
E_{w}(J,X),$ will also be used later.

\bigskip

\textbf{Lemma 4.11.} $\Delta _{t}^{n}p\in C_{w,0}(G,X)$ for all $%
p\in P_{w}(G,X),$ $t\in G$ and $n\geq 1.$

\,

\begin{proof} Let $p$\ have order $m$\ and note that

\smallskip

(4.3) $\qquad \Delta _{t}^{n}(\frac{p}{w})=\sum_{j=0}^{n}(-1)^{j}\binom {n} {j}\,(\frac{%
p}{w})_{(n-j)t}=\frac{\Delta
_{t}^{n}p}{w}-\sum_{j=0}^{n}(-1)^{j}\binom {n} {j}\,
(\frac{p}{w})_{(n-j)t}\frac{\Delta _{(n-j)t}w}{w}.$

\smallskip

\noindent Using (2.2) and $\Delta _{t}^{m+1}p=0$ we find $\Delta _{t}^{m+1}(\frac{p}{w}%
)\in C_{0}(G,X).$ By repeated application of the difference theorem for $w=1$
[16, Theorem 2.7] , we conclude $\Delta _{t}^{n}(\frac{p}{w}%
)\in $ \ $C_{0}(G,X)$ for $1\leq n\leq m+1$ and hence for all
$n\geq 1.$ By (4.1) again we conclude $\frac{\Delta
_{t}^{n}p}{w}\in C_{0}(G,X).$  \end{proof}

 \smallskip

\textbf{Theorem 4.12}. Let ${\F}$ be any translation
invariant closed subspace of $BC_{w}(J,X).$ If $\phi \in E_{w}(J,X)$ has $%
w$-mean $p$ and $\Delta _{t}\phi \in {\F}$ for each $t\in J,$ then
$\, \phi -p\in {{\F}} +C_{w,0}(J,X).$ If also $w=1,$ then $\phi
-p\in {\F}$

\begin{proof} For any finite subset $F\subseteq J$ we have
$\frac{\phi -p}{w}-R_{F}(\frac{\phi -p}{w})=-\frac{1}{\left|
F\right| }{\displaystyle{\sum_{t\in F}}}\,\,\Delta _{t}(\frac{\phi
-p}{w})$ and so $\phi -p=w \,\ R_{F}(\frac{\phi
-p}{w})-\frac{1}{\left| F\right| }{\displaystyle{\sum_{t\in F}}}\,\,\Delta _{t}\phi +\frac{1}{%
\left| F\right| }\sum_{t\in F}(\frac{\phi -p}{w})_{t}\Delta _{t}w+\frac{1%
}{\left| F\right| }{\displaystyle{\sum_{t\in F}}}\,\,\Delta
_{t}p.$ The first term on the right may be made arbitrarily small
in norm by suitable choice of $F.$ The second
term is in ${\F} $ by assumption, the third and fourth terms are in $%
C_{w,0}(J,X)$ by (2.2) and Lemma 4.11. If $w=1,$ then $\Delta
_{t}w=\Delta _{t}p=0$ which shows $\phi -p\in {\F}.$ \end{proof}

 \smallskip

We now use this theorem  to characterize $w$-ergodic functions.
For this we introduce the space of Maak $w$-ergodic functions
defined by

\smallskip

(4.4) \qquad $ME_{w}(J,X)$ the closed span of $\{\Delta _{t}\phi
:t\in J,\, \, \phi \in BC_{w}(J,X)\}$.

\smallskip

In the case $w=1$ the closed span of $\{\Delta _{t}\phi :t\in G,\,
\, \phi \in BUC(G,X)\}$ is studied in [48]. See also the
introduction to section  2 of  [17].

\smallskip

Since $\Delta_t \gamma = \gamma \Delta_t\gamma(0)$, for each
$\gamma \in \widehat {G}$, one gets $\gamma \, x \in ME_{w}(G,X)$
for each $x\in X$ and $\gamma \in \widehat {G}\setminus\{1\}$.

If $w$ is an unbounded  weight on $G$, then $G$  is  not compact.
This implies $\widehat {G}$ is not discrete. It follows that the
set $\{1\}$ is not open and  hence $\widehat {G}\setminus \{1\}$
is not closed. Therefore
  the trivial character $1$ can be
approximated  in $BC_w (G,\mathbf{C})$ by  characters from $
\widehat{G}\setminus\{1\}$.  One gets $ x \in ME_{w}(G,X)$ for
each $x\in X$. Hence each $X$-valued trigonometric polynomial
$\sum_{j=1}^n x_j \gamma_j \,\in ME_w (G,X)$.
 By Proposition 5.4 below and Proposition 4.9, it follows

\smallskip

(4.5) \qquad $C_{w,0}(J,X) \subset ME_{w}(J,X)$.

\bigskip

If also $w=1,$ then for any $%
\psi \in C_{0}(J,X)$ and any finite subset $F$ of $J$, we have $\psi =-%
\frac{1}{\left| F\right| }{\displaystyle{\sum_{t\in F}}}\,\,\Delta _{t}\psi +R_{F}\psi$. As $%
\left\| R_{F}\psi \right\| _{\infty }$ may be made arbitrarily
small, we conclude that $\psi \in ME(J,X).$ Hence

\smallskip

(4.6) \qquad $C_{0}(J,X) \subset ME (J,X)$.

\bigskip

 \textbf{Corollary 4.13}. If $G$ is compactly generated
$E^0_{w}(J,X)=ME_{w}(J,X).$

\,

\begin{proof} By Lemma 4.4 and the closedness of $E^0_{w}(J,X)$, we have $%
ME_{w}(J,X)\subseteq E^0_{w}(J,X)$. Conversely, let $\phi \in
E^0_{w}(J,X)$. Then $\Delta _{t}\phi \in ME_{w}(J,X)$ for all
$t\in J$ and hence  $\phi \in ME_{w}(J,X)$ by Theorem 4.12 and
(4.5) with $\F=  ME_{w}(J,X)$.  \end{proof}

 \bigskip

\textbf{Theorem 4.14}.  Let  $J\in \{ \mathbf {R_+, R}\} $ or $\{
\mathbf {Z_+, Z}\} $, $w_m (t) =(1 + |t|)^m$, $m\in \mathbf{N}$
and assume $w(s)\le k w(t)$ for $0\le s \le t $ and some positive
constant $k$. If $\phi \in E^0_ {w}(J,X)$ then $P^m\phi \in  C_{
ww_m, 0}(J,X)$ respectively $S^m\phi \in C_{ ww_m, 0}(J,X)$.

\begin{proof} We consider the case $m=1$ since the general case follows by induction.
Take $J\in \{ \mathbf {R_+, R}\} $, $m=1$. Let $\phi\in BUC_w
(J,X)$, $t\in J$. Then $||P\phi (t)||\le k |t| ||\phi
||_{w,\infty} w (t)$ and so $P$ maps $BUC_w (J,X)$ continuously
into $BUC_{ww_1} (J,X)$. Moreover $P(\Delta_t \phi)= \Delta_t (P
\phi) - P \phi (t)$. By  Lemma 4.4, one has $\Delta_t (P \phi)\in
E^0_ {ww_1}(J,X)$.
 Since  the constant $P\phi (t)\in \,E^0_ {ww_1}(J,X)$ we conclude
 $P(\Delta_t \phi)\in E^0_ {ww_1}(J,X)$. By Corollary 4.13 and the continuity of $\,P\,$
we conclude that $\,P\,\,$ maps $\,\,E^0_ {w}(J,X)\,\,\,$ into
$\,\,\,E^0_ {ww_1}(J,X)$.
 Finally $\,\,||\Delta_t P \phi(s)||\,\, \le$
 $ k |t|  w(t)w(s)
||\phi||_{w, \infty} $ for all $s\in J$. Hence $ \Delta_t (P
\phi)\in C_{ww_1, 0}(J,X)$. If $\phi \in E^0_ {w}(J,X)$, we can
apply Theorem 4.12 to $P\phi$ to obtain $P\phi \in   C_{ww_1,
0}(J,X)$.
 The case  $J\in \{ \mathbf {Z_+, Z}\} $ follows similarly.
\end{proof}

\bigskip

\smallskip

5. $w$-\textbf{almost periodic functions }

In this section we introduce a new class called weighted almost
periodic functions. This new class is a natural generalization of
the space of  continuous almost periodic functions $\phi: G\to X$
denoted by $AP(G,X)$.

So, firstly we recall some basics about $AP(G,X)$ (see [3], [43]
for the case $G=\mathbf{R}$ and [8] for any $G$).

A subset $E\subset G$ is said to be $relatively\,\,\, dense$ in
$G$ if  there is a finite set $F\subset G$ such that $G=
\cup_{a\in F} (E+a)$.

 A continuous function $\phi: G\to X$ is
said to be $(Bohr)\,\,\, almost\,\, periodic\,$ if  for each
$\varepsilon
>0$ the set  $E(\phi,\varepsilon)= \{\tau\in G: ||\phi(t+\tau)
-\phi(t))|| \le \varepsilon,\,\,\, t\in G \}$ of its
$\varepsilon$-$periods$ is relatively dense in  $G$.

The space of all  continuous $X$-valued almost periodic functions
defined on $G$ is denoted by $AP(G,X)$. It is well known that
$AP(G,X)$ is a subspace of $BUC(G,X)$. Moreover,
 if $\phi_1, \cdots,\phi_n \in AP(G,X)$, then
$\cap_{k=1}^{n}E(\phi_k,\varepsilon)$ is relatively dense in $G$
for each $n\in \mathbf{N}$.

Secondly, we denote the space of $w$-\textit{trigonometric
polynomials} $\psi
=\sum_{j=1}^{n}\gamma _{j}p_{j},$ where $\gamma _{j}\in \widehat{G}$ and $%
p_{j}\in P_{w}(G,X),$ by $TP_{w}(G,X).$ The closure of $TP_{w}(G,X)$ in $%
BC_{w}(G,X)$ is the space $AP_{w}(G,X)$ of $w$-\textit{almost periodic }$X$%
\textit{-}valued functions on\textbf{\ }$G.$

\bigskip

\textbf{Proposition 5.1.} If $G$ is compactly generated, then $%
AP_{w}(G,X)\subset BUC_{w}(G,X)$. Moreover, for each $\phi \in
AP_{w}(G,X)$  the range of $\frac{\phi}{w}$ is relatively compact.

\begin{proof} If $p\in P_{w}(G,X)$ then $p\in BUC_{w}(G,X)$ by
Proposition 3.7(d).   By Proposition 3.7 (b) and Proposition 3.4,
we
 conclude that the space  $X_p $ generated by the range of $p$ is
finite dimensional. Moreover, if $\gamma \in \widehat{G}$ then
$\Delta _{h}(\gamma p)(t)=\gamma (t)\gamma (h)\Delta
_{h}p(t)+\gamma (t)p(t)\Delta _{h}\gamma (0)$ and $(\gamma(t)
p(t))/w(t) \in X_p$ for all $t\in J$, and so $\gamma p\in
BUC_{w}(G,X)$ and
 the range of $(\gamma p)/w$ is relatively compact. The
results follow.
\end{proof}

\bigskip

\textbf{Proposition 5.2.} If $P_w(G,\mathbf{C})$ is finite
dimensional, then $AP_{w}(G,X)\subset E_{w}(G,X).$

\,

\begin{proof} By Proposition 4.1 $E_w (G,X)$ is closed in $BC_w(G,X)$. Take any $\gamma \in \widehat{G}$, $\gamma \neq 1$\ \ and
any $p\in P_{w}(G,X),$ $p\neq 0.$ It suffices to prove $\frac{\gamma p}{w}%
\in E^{0}(G,X).$ To do this suppose firstly that $\gamma (s)$ has
finite order $m\neq 1$ for some $s\in G$ and let $F_{1}=\{(j-1)s:$
$1\leq j\leq m\}. $ Then ${\displaystyle{\sum_{t\in
F_{1}}}}\,\,\gamma (t)=\sum_{j=1}^{m}\gamma (s)^{j-1}=0.$
Hence, $R_{F_{1}}(\frac{\gamma p}{w})=\frac{1}{m} {\displaystyle{\sum_{t\in F_{1}}}}\,\,(\frac{%
\gamma p}{w})_{t}=\frac{1}{m}{\displaystyle{\sum_{t\in F_{1}}}}\,\,\gamma _{t}\frac{\Delta _{t}p}{%
w}-\frac{1}{m}{\displaystyle{\sum_{t\in F_{1}}}}\,\,\frac{\gamma
_{t}p_{t}}{w_{t}}\frac{\Delta _{t}w}{w}+\frac{\gamma p}{mw}
\,\sum_{t\in F_{1}}\gamma (t).$ The last term on
the right is $0$ and so, by (2.2) and Lemma 4.11, $R_{F_{1}}(\frac{\gamma p}{w%
})\in C_{0}(G,X).$ Since $C_{0}(G,X)\subseteq E^{0}(G,X),$ for each $%
\varepsilon >0$ we can find a finite set $F_{2}\subset G$ such
that $\left\| R_{F_{2}}R_{F_{1}}(\frac{\gamma p}{w})\right\|
_{\infty }<\varepsilon$.
Since $R_{F_{2}}R_{F_{1}}(\frac{\gamma p}{w})=R_{F_{2}+F_{1}}(\frac{\gamma p%
}{w})$ we conclude that $\frac{\gamma p}{w}\in E^{0}(G,X).$ If $\gamma (s)$
has infinite order for every non-zero $s\in G$ then $\gamma $ has dense
range in $\mathbf{T}$. Given $\varepsilon >0$ we can find $s\in G$ such that
$\left| \gamma (s)+1\right| <\frac{\varepsilon }{\left\| p\right\|
_{w,\infty }}.$ If $F_{1}=\{0,s\},$ then $R_{F_{1}}(\frac{\gamma p}{w})=%
\frac{1}{2}\gamma _{s}\frac{\Delta _{s}p}{w}-\frac{1}{2}\frac{\gamma
_{s}p_{s}}{w_{s}}\frac{\Delta _{s}w}{w}+\frac{\gamma p}{2w}(\gamma (s)+1).$
This time choose a finite set $F_{2}\subset G$ such that $\left\| R_{F_{2}}(%
\frac{1}{2}\gamma _{s}\frac{\Delta _{s}p}{w}-\frac{1}{2}\frac{\gamma
_{s}p_{s}}{w_{s}}\frac{\Delta _{s}w}{w})\right\| _{\infty }<\frac{%
\varepsilon }{2}.$ Since $\left\| R_{F_{2}}(\frac{\gamma p}{2w}(\gamma
(s)+1))\right\| _{\infty }\leq \left\| \frac{\gamma p}{2w}(\gamma
(s)+1)\right\| _{\infty }<\frac{\varepsilon }{2}$ we obtain, as before, $%
\left\| R_{F_{2}}R_{F_{1}}(\frac{\gamma p}{w})\right\| _{\infty
}<\varepsilon$.  \end{proof}

\smallskip

\textbf{Lemma 5.3.}
 (a) Let $\phi \in C(G,X)$ have support in a compact neighbourhood $K$\ of $%
0$ in $G$ and suppose $K$ generates $G.$ Then there is a sequence of
trigonometric polynomials $\psi _{n}:G\rightarrow X$ converging uniformly to
$\phi $ on $K$ with ${\displaystyle{\sup_{t\in G}}}\,\,\left\| \psi _{n}(t)\right\| <\frac{1}{n}%
+{\displaystyle{ \sup_{t\in G}}}\,\,\left\| \phi (t)\right\|$.

(b) If also $G=\mathbf{R}$ and $\phi $ is$\ m$-times continuously
differentiable on $\mathbf{R,}$ then $\psi _{n}$ may also be chosen so that $%
\psi _{n}^{(j)}\rightarrow \phi ^{(j)}$ uniformly on $K$ \ and
${\displaystyle{ \sup_{t\in G}}}\,\,\left\| \psi
_{n}^{(j)}(t)\right\| <\frac{1}{n}+{\displaystyle{ \sup_{t\in
G}}}\,\, \left\| \phi ^{(j)}(t)\right\| $ for $0\leq j\leq m.$

\,

\begin{proof} (a) By [53, 2.4.2] there exists a closed
subgroup $G_{1}$\ of $G$\ with $G/G_{1}$ compact, $G_{1}$ isomorphic to $%
Z^{m}$ for some $m\ge 0$ and $K\cap G_{1}=\{0\}.$ Since $(K-K)\cap
G_1 \,$ is finite, we may also assume (replacing $G_{1}$ by one of
its subgroups of finite index if necessary) that $(K-K)\cap
G_1\,=\{0\}$.
 Let $\pi :G\rightarrow G/G_{1}$ be the quotient mapping and define $
\widetilde{\phi }: G/G_{1}\rightarrow X$ by $\widetilde{\phi }(\pi
(t))=\phi (t)$ if $t\in K$ and $\widetilde{\phi }(\pi (t))=0$ if
$\pi (t)\notin \pi (K).$ Since $(K-K)\cap G_{1}=\{0\},$\
$\widetilde{\phi }$\
is well-defined. Moreover, $\widetilde{\phi }\in C(G/G_{1},X)$ and $%
\widetilde{\phi }\circ \pi |_{K}=\phi |_{K}.$ Since $G/G_{1}$ is
compact, there is a sequence of trigonometric polynomials $\eta
_{n}:G/G_{1}\rightarrow X$ converging to $\widetilde{\phi }$ in $%
C(G/G_{1},X).$ We choose $\eta _{n}$ so that $\left\| \eta _{n}-\widetilde{%
\phi }\right\| _{\infty }<\frac{1}{n}$ and let $\psi _{n}=\eta
_{n}\circ \pi ,$ a trigonometric polynomial. Then $(\psi
_{n}|_{K})$ converges
uniformly to $\phi |_{K}$ and $\left\| \psi _{n}\right\| _{\infty }=$ $%
\left\| \eta _{n}\right\| _{\infty }<\frac{1}{n}+\left\| \widetilde{\phi }%
\right\| _{\infty }=\frac{1}{n}+\left\| \phi \right\| _{\infty }.$

(b) In this case, $G/G_{1}\cong \mathbf{T}$ the circle
group, and so $\widetilde{\phi }\in C^{m}(\mathbf{T},X)\ $the space of $m$%
-times continuously differentiable functions from $\mathbf{T}$ to
$X$. With the norm $\left\| u\right\|
_{(m)}=\sum_{j=0}^{m}\frac{1}{j!}\left\| u^{(j)}\right\| _{\infty
},$ $C^{m}(\mathbf{T},X)$ is a homogeneous Banach space in the
sense of [39, 2.10]. By [39, Theorem
2.11]\ we may choose $\eta _{n}$ so that $\left\| \eta _{n}-\widetilde{%
\phi }\right\| _{(m)}<\frac{1}{n(m!)}.$ (This last reference deals with $%
\mathbf{C}$-valued functions, but the same proof is valid for
general $X$.) Setting $\psi _{n}=\eta _{n}\circ \pi $ we obtain
the lemma. \end{proof}

\smallskip

\textbf{Proposition 5.4.} If $w$ is an unbounded weight and $G$ is compactly
generated, then $TP(G,X)$ is dense in $C_{w,0}(G,X).$

\,

\begin{proof} Take $\phi \in C_{w,0}(G,X)$ and $n\in
\mathbf{N.}$ Since
$G$ is compactly generated, there is a compact generating neighbourhood $%
F_{n}$ of $0$ in $G$ such that $\left\| \phi (t)\right\|
<\frac{1}{n}w(t)$
for all $t\notin F_{n}$. As $w$ is unbounded there is a compact set $%
H_{n}\supset F_{n}$ such that ${\displaystyle{\sup_{t\in
F_{n}}}}\,\,\left\| \phi (t)\right\| <\,\, \frac{1}{n}
{\displaystyle{\sup _{t\in H_{n}}}}\,\,w(t).$ Let $K_{n}$ be a
compact set whose interior contains $H_{n}.$ By Proposition 2.3 we
may assume ${\displaystyle{\sup_{t\in K_{n}}}}\,\,w(t)\leq $
${\displaystyle{\inf_{t\notin K_{n}}}}\,\, w(t)$ and so
${\displaystyle{\sup_{t\in K_{n}}}}\left\| \phi (t)\right\| <$
$\frac{1}{n} {\displaystyle{\inf_{t\notin K_{n}}}}\,\,w(t).$ Since
$G$ is
a regular space [41, p. 146]\, there is a continuous function $%
f:G\rightarrow [0,1] $ which is $1$ on $H_{n}$ and $0$ outside $%
K_{n}.$ By Lemma 5.3 there exists $\psi _{n}\in TP(G,X)$ with
${\displaystyle{\sup_{t\in K_{n}}}}\,\,\left\| f(t)\phi (t)-\psi
_{n}(t)\right\| <\frac{1}{n}$ and ${\displaystyle{\sup _{t\notin
K_{n}}}}\,\, \left\| \psi _{n}(t)\right\| <\frac{1}{n}+
{\displaystyle{\sup_{t\in G}}}\,\,\left\| f(t)\phi (t)\right\|$.
Hence

$\qquad\left\| \phi -\psi _{n}\right\| _{w,\infty }\leq
{\displaystyle{\sup_{t\in K_{n}}}}\,\, \frac{\left\| \phi
(t)-f(t)\phi (t)\right\| }{w(t)}+
 {\displaystyle{\sup_{t\in K_{n}}}}\,\,
\frac{\left\| f(t)\phi (t)-\psi _{n}(t)\right\|
}{w(t)}+{\displaystyle{\sup_{t\notin K_{n}}}}\,\, \frac{\left\|
\phi (t)\right\| }{w(t)}$

$ \qquad + {\displaystyle{\sup{_{t\notin K_{n}}}}}\,\,\frac{
\left\| \psi _{n}(t)\right\| }{w(t)}$ $< {\displaystyle{\sup_{t\in
K_{n}\backslash
H_{n}}}}\,\,\frac{\left\| \phi (t)\right\| }{w(t)}+\frac{1}{n}+\frac{1}{n}+%
\frac{\frac{1}{n}+{\displaystyle{\sup_{t\in K_{n}}}}\,\,\left\|
f(t)\phi (t)\right\|} {\displaystyle{\inf_{t\not\in K_{n}}} \,\,
w(t)} < \,\, \frac{5}{n}$.

This completes the proof. \end{proof}

\smallskip

 \textbf{Proposition 5.5. } As usual let   $w_N (t)=(1+|t|)^N$, \,  $N\in \mathbf{N}$.

 (a) If $w$ is an unbounded weight
and $G$ is compactly generated, then $C_{w,0}(G,X)\subseteq
AP_{w}(G,X).$

(b) If $p\in P_w(G,X)$, then  $p\cdot AP(G,X) \subset
AP_{w}(G,X)$.

(c) If $G\in \{ \mathbf{Z,R}\}$, then $AP_{w_N}(G,X)= t^N\cdot
AP(G,X)\oplus C_{w_N,0} (G,X)$. Moreover, if $\phi = t^N \psi+
\xi$  with $ \psi\in AP(\mathbf{R},X)$ and $\xi \in C_{w_N,0}
(\mathbf{R},X)$ then $||\psi||_{\infty}\le ||\phi||_{w_N,\infty}$.

(d) Let    $\phi \in AP_{w_N}(\mathbf{R},X)$ and  $\phi = t^N
\psi+ \xi$  with $ \psi\in AP(\mathbf{R},X)$ and $\xi \in
C_{w_N,0} (\mathbf{R},X)$. Let $\phi' \in
BUC_{w_N}(\mathbf{R},X)$. Then $\phi' \in AP_{w_N}(\mathbf{R},X)$,
$\psi' \in AP(\mathbf{R},X)$ and $\xi'\in C_{w_N,0}
(\mathbf{R},X)$.

(e) If  $\phi \in AP(\mathbf{R},X)$ and $P^2 \phi \in  AP _{w_1}
(\mathbf{R},X)$, then there is $a\in X$ such that $ (P^2
\phi-at)\in C_{w_1,0} (\mathbf{R},X)$.

\begin{proof} (a) This follows directly from Proposition 5.5.

 (b) By [57, p. 330, Corollary] $TP(G,X)$ is dense in $AP(G,X)$. It follows that if $ \phi \in AP(G,X)$, then $p\phi$
is a limit in $BUC_w (G,X)$ of a sequence from $p\cdot TP (G,X)$.

(c) Note firstly that the sum on the right is direct. For suppose
$t^{N}\psi _{1}+\xi _{1}=t^{N}\psi _{2}+\xi _{2}$ for $\psi
_{j}\in AP(G,X)$
and $\xi _{j}\in C_{w,0}(G,X).$ Set $q(t)=(1+t)^{N}-t^{N}$ and $J=%
G_{+}.$ Then $\left( \psi _{1}-\psi _{2}\right) |_{J}=\frac{1}{w}%
(\xi _{2}-\xi _{1}-q\psi _{2}+q\psi _{1})|_{J}\in C_{0}(J,X).$
This is impossible unless $\psi _{1}=\psi _{2}.$ The sum is also
topological. For
suppose $\phi _{n}=t^{N}\psi _{n}+\xi _{n}$ where $\psi _{n}\in AP(%
G,X)$ and $\xi _{n}\in C_{w,0}(G,X)$ and $\left( \phi
_{n}\right) $ converges to $\phi $ in $BC_{w}(G,X).$ Then

\smallskip

(5.1)\qquad  $\frac{%
\phi _{n}}{w}|_{J}=\psi _{n}|_{J}+\frac{\xi _{n}-q\psi
_{n}}{w}|_{J}\in AP(G,X)|_{J}\oplus C_{0}(J,X) =AAP(J,X).$

\smallskip

\noindent But this last sum is a topological
direct sum  (see \cite{RSW}) or [9, Proposition 2.2.2] and references therein) and so $\left( \psi _{n}|_{J}\right) $ converges to $\psi |_{J%
\text{ }}$ for some $\psi \in AP(G,X).$ Hence $\left( \psi
_{n}\right) $ converges to $\psi $ in $AP(G,X)$ and $\left(
t^{N}\psi _{n}\right) $ converges to $t^{N}\psi $ in
$AP_{w}(G,X)$.
 It follows that $\left( \xi _{n}\right) $ converges to\ some $\xi $ in $%
C_{w,0}(G,X)$ and that $\phi =$ $t^{N}\psi +\xi$.

Next, given $\phi \in AP_{w}(G,X)$ we may choose a sequence $%
\left( \pi _{n}\right) \subset TP_{w}(G,X)$ converging to $\phi $
in $AP_{w}(G,X).$ But $\pi _{n}=t^{N}\psi _{n}+\xi _{n}$ where $%
\psi _{n}\in TP(G,X)$ and $\xi _{n}\in C_{w,0}(G,X).$ It follows
from the previous paragraph that $\phi =$ $t^{N}\psi +\xi $ for
some $\psi \in AP(G,X)$ and $\xi \in C_{w,0}(G,X).$ The converse
follows from (a) and (b). The inequality $||\psi||_{\infty}\le
||\phi||_{w_N,\infty}= ||\phi/w_N||_{\infty}$ follows from [9,
Proposition 2.2.2].

(d) One has $n[\phi(\cdot+ 1/n)]-\phi(\cdot)]-\phi'(\cdot)= n
\int_{0}^{1/n}\, [\phi'(\cdot+s)]-\phi'(\cdot)]\,ds\to 0$  in
$BUC_{w_N} (\mathbf{R},X)$ as $n\to \infty$. This implies $\phi'
\in AP_{w_N}(\mathbf{R},X)$.  By part (c), $\phi' = t^N \eta +\mu$
where $\eta \in AP(\mathbf{R},X)$  and $\mu\in
C_{w_N,0}(\mathbf{R},X)$. Now $\phi/w_N |\, \mathbf{R}_+ = \psi +
(\xi+ \psi(t^N-w_N))/w_N |\, \mathbf{R}_+ \in
AAP(\mathbf{R}_+,X)$.  By (5.1), $\psi'$, $\xi'$ exist on
$\mathbf{R}_+$ and
   $(n[\psi (\cdot+
 1/n)]-\psi(\cdot)] |\, \mathbf{R}_+)$, $(n[\xi (\cdot+
 1/n)/w_N (\cdot+
 1/n)]-\xi(\cdot)/w_N(\cdot)] |\, \mathbf{R}_+)$ are  Cauchy  sequences in $BUC
 (\mathbf{R}_+,X)$. It follows  $ \psi'|\,
 \mathbf{R}_+ \in AP(\mathbf{R},X)|\,
 \mathbf{R}_+$  and $ \xi'/w_N|\,
 \mathbf{R}_+ \in C_0 (\mathbf{R}_+,X)$. Similarly   $ \psi'|\,
 \mathbf{R}_- \in  \in AP(\mathbf{R},X)|\,
 \mathbf{R}_-$  and $ \xi'/w_N|\,
 \mathbf{R}_- \in C_0 (\mathbf{R}_-,X)$. Now $\phi' (t)= N t^{N-1}\psi (t) + t^N \psi' (t)
  + \xi' (t)$
 for all $t\not = 0$. From this we get $\xi' $ is continuous at $t=0$.
 This implies $ \xi'/w_N
  \in C_0 (\mathbf{R},X)$ and so $\psi' =\eta$ and $N t^{N-1} \psi +\xi' =\eta$.

(e) By Theorem 4.14, $P^2 (\phi -M\phi)\in C _{w_2,0}
(\mathbf{R},X)$.  Since $P^2 \phi \in AP _{w_1} (\mathbf{R},X)
\subset C _{w_2,0}
 (\mathbf{R},X)$ and  $P^2 \phi = P^2 (\phi -M\phi) +( M\phi/2)\,
 t^2$ we conclude $M\phi =0 $. By Theorem 4.14, $P\phi \in C_{w_1,0}
(\mathbf{R},X)$. By (c), $P^2 \phi =t\, \psi +\xi$ with $ \psi\in
AP  (\mathbf{R},X)$, and $ \xi\in C _{w_1,0} (\mathbf{R},X)$. By
(d), one gets $P\phi =t\, \psi' + \psi+\xi'$ where $\psi'\in AP
(\mathbf{R},X)$. This implies $\psi'=0$ and hence $\psi=a \in X$.
This shows $ P^2 \phi -t\,a\in C _{w_1,0} (\mathbf{R},X)$.
 \end{proof}

\smallskip

\textbf{Proposition 5.6}. \ If $\phi \in AP_{w}(G,X)$ then $%
\phi \ast f\in AP_{w}(G,X)$ for all $f\in L_{w}^{1}(G)$. If also
$G$ is compactly generated, the operator $T_{\phi
}:L_{w}^{1}(G)\rightarrow AP_{w}(G,X)$ defined by $T_{\phi
}(f)=\phi \ast f$\ is compact.

\,

\begin{proof} Let $\pi =\gamma p$ where $\gamma \in \widehat{G},$ $p$\ $%
\in P_{w}(G,X)$ and let $f\in L_{w}^{1}(G).$ Then $\pi \ast f=\gamma (p\ast
(\gamma ^{-1}f))\in TP_{w}(G,X).$ Choose a sequence $(\pi _{n})$ in $%
TP_{w}(G,X)$ converging to $\phi $ in $BC_{w}(G,X).$ It follows
that $\pi _{n}\ast f\rightarrow \phi \ast f$ and so $\phi \ast
f\in AP_{w}(G,X).$ Since $\left\| T_{\phi }(f)\right\| \leq
\left\| \phi \right\| _{w,\infty }\left\| f\right\| _{w,1},$ each
$T_{\phi }$ is bounded. If $G$
is compactly generated, then by Proposition 3.7(b), $p(t-s)=%
\sum_{j=1}^{k}p_{j}(t)q_{j}(s)$ for $s,t\in G$\ where $p_{j}\in
P_{w}(G,X)$ and $q_{j}\in P_{w}(G,\mathbf{C}).$ So $\pi \ast
f=\sum_{j=0}^{k}\gamma p_{j}c_{j}$ where $c_{j}=\int_{G}\gamma
^{-1}(s)q_{j}(s)f(s)d\mu (s).$ Hence $\pi \ast f\in $
span$\{\gamma p_{j}:1\leq j\leq k\}$ and so $T_{\pi }$ has finite
rank. As $T_{\pi _{n}}\rightarrow T_{\phi }$ in the operator norm,
$T_{\phi }$ is compact. \end{proof}

\smallskip

Following [9] where the case $w=1$ is considered, we say that
${\F} $ is a $\Lambda _{w}^{0}$\textit{-class} if it is a
translation invariant subspace of $BC_{w}(J,X)$ satisfying

\bigskip

(5.2) \qquad\qquad if $\phi \in BUC_{w}(G,X)$ and $\phi |_{J}\in {\F}$ then $%
\phi _{t}|_{J}\in {\F} $ for all $t\in G;$

(5.3)\qquad \qquad ${\F} $ is a closed linear subspace of
$BUC_{w}(J,X);$

(5.4)\qquad\qquad ${\F} $ is closed under multiplication by
characters;

(5.5)\qquad\qquad if $w$ is unbounded, ${\F} $ $\supseteq
C_{w,0}(J,X).$

\bigskip

A $\Lambda _{w}^{0}$-class ${\F} $ is called a $\Lambda _{w}$\textit{%
-class} if ${\F} $ $\supseteq P_{w}(J,X).$

\bigskip

Note that if ${\F}$ is a subspace of  $BC_{w}(J,X)$ and $\phi \in
{\F} $, then $\phi_t = \Delta_t \phi +\phi\in  E^0_w + \F  $ by
Lemma 4.4. This means that $E^0_w \subset {\F} $ implies (5.2).

Note that for the case $J=G,$ (5.2) just states that ${\F} $ is
translation invariant. In particular, $C_{w,0}(G,X)$ is a $\Lambda _{w}^{0}$%
-class. So too is $C_{w,0}(J,X)\ $ when $G$ is $\mathbf{R}$ or $%
\mathbf{Z}$ and $J$ is $\mathbf{R}_{+}$ or $\mathbf{Z}_{+}$. Many
examples of $\Lambda _{w}$-classes for the case $w=1$ are given in
\cite{BBA}. These include almost periodic, almost automorphic and
absolutely recurrent functions. The class $0_{J},$ consisting of
just the zero function from $J$ to $X,$ satisfies (5.2) only when
$J=G.$ From Proposition 5.1 and Proposition  5.5 (a) we obtain:

\bigskip

\textbf{Proposition 5.7}. If $G$ is compactly generated or $w$ is bounded,
then $AP_{w}(G,X)$ is a $\Lambda _{w}$-class.

\bigskip

We mention one other family of $\Lambda _{w}^{0}$-classes. The
partial ordering $\leq $ , defined by $s\leq t$ whenever $t\in s+J,$ makes $%
J $ a directed set. So we may define ${{\F}} _{w}(J,X)=\{\phi \in
BUC_{w}(J,X): {\displaystyle{\lim_{t\in J}}}\,\, \left\|
\frac{\phi (t)}{w(t)}\right\| =0\}.$ If $G=\mathbf{R}$ or
$\mathbf{Z}$ and $J=\mathbf{R}_{+}$ or $\mathbf{Z}_{+},$ then
${{\F}} _w (J,X)\,= C_ {w,0}(J,X)$ but in general this is not the
case.
If $G=\mathbf{R}^{d}$ or $\mathbf{Z}^{d}$ and $J=(\mathbf{R}_{+})^{d}$ or ($%
\mathbf{Z}_{+})^{d},$ then ${{\F}} _{w}(J,X)$ is a $\Lambda _{w}^{0}$%
-class\ contained in every $\Lambda _{w}^{0}$-class of functions from $J$ to
$X$. \ These spaces will be developed and used in a subsequent paper.

\bigskip
\textbf{Example 5.8}. If $G\in \{\mathbf{R, Z}\}$  we define

 (5.6) \qquad \qquad $\, \,\, \A=P_N AP(G,X) =\sum_{k=0}^N t^k \cdot AP (G,X)$

\qquad \qquad\qquad \qquad $=\{\phi: \phi(t)= \sum_{k=0}^N t^k
\psi_k (t) $ where $ \psi_k \in AP (G,X),  t\in G\}$
\smallskip

\noindent Then   the following properties are clear

\smallskip

(5.7)\qquad\qquad$\A$  is a linear subspace of $AP_{w_N} (G,X)$,
translation invariant,  closed under

\qquad\qquad\qquad convolution with elements from
$L^1_{w_N}(G,\mathbf {C})$, closed under multiplication

\qquad\qquad\qquad  by characters from $\widehat{G}$ and dense in
$AP_{w_N} (G,X)$.

\noindent  In Proposition 6.3  we investigate harmonic analysis
with respect to classes $\A=P_N AP$ defined by  (5.6).

\bigskip

 6. \textbf{Spectral analysis}

Throughout this section we will assume that $G$ is compactly generated and
that ${\F} $ is a translation invariant closed subspace of $%
BC_{w}(J,X)$ satisfying (5.2).
\bigskip

Recall that for an  ideal $I$ of $L_{w}^{1}(G),$ the
\textit{cospectrum} of $I$ is given by cosp$(I)=\{\gamma \in \widehat{G}:%
\widehat{f}(\gamma )=0$ for all $f\in I\}$. Clearly, cosp $I\, =$
cosp$\overline {I}$, where $\overline{I}$ is the closure of $I$ in
$L_{w}^{1}(G)$. Since $L_{w}^{1}(G)$ is a Wiener
algebra, its maximal ideals are the sets $I_{\alpha }=\{f\in L_{w}^{1}(G):%
\widehat{f}(\alpha )=0\}$ where $\alpha \in \widehat{G}$ and its
primary ideals are those whose cospectrum is a singleton. By
Wiener's tauberian theorem, all (closed) primary ideals in
$L^{1}(G)$ are maximal (\cite{BE}). This is not the case for
general $L_{w}^{1}(G).$ For example, it is proved in [15, Theorem
3.4]\ that, for a weight of polynomial growth $N,$\ the primary
ideals of $L_{w}^{1}(\mathbf{Z})$ are the sets $I_{k}=\{f\in
L_{w}^{1}(\mathbf{Z}):\widehat{f}^{(j)}(1)=0$ for $0\leq j\leq k\}$ where $%
0\leq k\leq N.$ If $g\in L_{w}^{1}(\mathbf{Z})$ and $f=\Delta _{t}^{N+1}g$
for some $t\in \mathbf{Z}$ then $\widehat{f}(\gamma )=$ $(\gamma (t)-1)^{N+1}%
\widehat{g}(\gamma )$ showing $f\in I_{k}.$ The following theorem
generalizes this result. For $t=(t_{1},...,t_{m})\in G^{m}$ we write $\Delta
_{t}g=\Delta _{t_{1}}...\Delta _{t_{m}}g.$

\bigskip

\textbf{Theorem 6.1}. Assume $w$\ has polynomial growth $N$\ and
$I$ is a closed ideal of $L_{w}^{1}(G)$ with cosp$(I)=\{1\}.$ Then
$\Delta _{t}g\in I$ for all $g\in L_{w}^{1}(G)$ and $t\in
G^{N+1}.$

\begin{proof} Consider the annihilator $I^{\bot }=\{\phi \in
L_{w}^{\infty }(G):\phi \ast f=0$ for all $f\in I\},$ a closed
translation invariant subspace of $L_{w}^{\infty }(G).$ If $\phi
\in I^{\bot }$ and $I_{w}(\phi )=\{f\in L_{w}^{1}(G):\phi \ast
f=0\}$ then $I_{w}(\phi )\supseteq I.$ This implies that
cosp$(I_{w}(\phi ))\subseteq $ cosp$(I)=\{1\}$.  By [15, Theorem
3.4]\ or Theorem 6.6 below, $\phi \in P_{w}(G,\mathbf{C})$ and so
$\Delta _{t}\phi =0$ for all $t\in G^{N+1}.$ If $g\in
L_{w}^{1}(G)$ then $\phi \ast \Delta _{t}g=\Delta _{t}\phi \ast
g=0,$ showing $\Delta _{t}g\in I^{\bot \bot }.$ Since $I^{\bot
\bot }=I$ the theorem is proved. \end{proof}

\smallskip

Let $\phi \in BC_{w}(G,X).$ The set $I_{w}(\phi )=\{f\in
L_{w}^{1}(G):\phi \ast f=0\}$ is a closed ideal of $L_{w}^{1}(G)$
and the \textit{Beurling spectrum} of $\phi $ is defined to be
$sp_{w}(\phi )=$ cosp$(I_{w}(\phi ))$. More generally, following
[9, p.20] , set $I_{{\F} }(\phi )=\{f\in L_{w}^{1}(G):(\phi \ast
f)|_{J}\in {\F} \}.$ By condition (5.2), $I_{{\F} }(\phi )$ is a
closed translation invariant subspace of $L_{w}^{1}(G)$\ and is
therefore an
ideal. We define the \textit{spectrum of} $\phi $ \textit{relative to} $%
{\F} $ (respectively $\A$) to be $sp_{{\F} }(\phi )=$ cosp
$(I_{{\F} }(\phi ))$ ($sp_{\A}(\phi )=$ cosp $(I_{{\A} }(\phi
))$).

\bigskip

\textbf{Lemma 6.2}.  For each $\phi \in BUC_{w}(G,X)$ there is
a sequence of approximate units, that is a sequence $(g_{n})$ in $%
L_{w}^{1}(G)$ such that $\phi \ast g_{n}\rightarrow \phi $ in $%
BUC_{w}(G,X).$

\,

\begin{proof}  Since $\phi $ is $w$-uniformly continuous, there
is a compact neighbourhood $V_{n}$ of $0$ in $G$ such that
$\left\| \phi _{-s}-\phi \right\| _{w,\infty }<\frac{1}{n}$ for
all $s\in V_{n}.$ Choose $g_{n}\in C_{c}(G)$ with
supp$(g_{n})\subseteq V_{n}$, $g_{n}\geq 0$
and $\int_{G}g_{n}(s)d\mu (s)=1.$ So $g_{n}\in L_{w}^{1}(G)$ and for each $%
t\in G,$ $\left\| \phi \ast g_{n}(t)-\phi (t)\right\| =\left\|
\int_{V_{n}}[\phi (t-s)-\phi (t)] g_{n}(s)d\mu (s)\right\|
<\frac{1}{n}w(t)$ (see Theorem 33.12 in \cite{HR}). \end{proof}

\smallskip

\textbf{\textbf{Proposition 6.3}}. (a) Let $\phi ,\psi \in
BC_{w}(G,X).$

(i) $sp_{{\F} }(\phi _{t})=sp_{{\F} }(\phi )$ for all $t\in G;$

(ii)$\ sp_{{\F} }(\phi \ast f)\subseteq sp_{{\F} }(\phi )\cap
supp(\widehat{f})$ for all $f\in L_{w}^{1}(G);$

(iii) $sp_{{\F} }(\phi +\psi )\subseteq sp_{{\F} }(\phi )\cup
sp_{{\F} }(\psi );$

(iv)$\ sp_{{\F} }(\gamma \phi )=\gamma +sp_{{\F} }(\phi ),$
provided ${\F} $ is invariant under multiplication by $\gamma \in
\widehat{G};$

(v)\ Let $\phi \in BUC_{w}(G,X).\ $Then $\phi _{t}|_{J}\in {\F} $
for all $t\in G$ if and only if $sp_{{\F} }(\phi )=\emptyset ;$

(vi) Let $\phi \in BUC_{w}(G,X).\ $If $\Delta
_{t}^{N+1}\phi _{s}|_{J}\in {\F} $ for all $s,t\in G$ then $%
sp_{{\F} }(\phi )\subseteq \{1\};$ and conversely if $w$ has
polynomial growth $N.$

(b) Statements (i), (ii), (iii), (iv), (vi) are true  when $\F$ is
replaced by a class $\A$ defined  by (5.6). We do not know whether
or not statement (v)  is  true for classes $\A$ defined by (5.6),
$N\in \mathbf{N}$. Instead, one has

(v$'$)  Let $\phi \in  \A$, then $sp_{{\A} }(\phi )=\emptyset $.

 (v$''$) If $sp_w (\phi)$ is compact and $sp_{{\A} }(\phi
)=\emptyset $, then $\phi \in \A$.

\,

\begin{proof} (a) The arguments are the same as for the Beurling
spectrum. See for example [28, part II, p.988] or \cite{R}. We
present proofs for (v) and (vi).

Firstly assume $\phi _{t}|_{J}\in {\F} $ for all $t\in G.$ The map $%
t\mapsto \phi _{t}|_{J}:G\rightarrow {\F} $ is continuous and so,
for each $f\in L_{w}^{1}(G)$ we have $(\phi \ast
f)|_{J}=\int_{G}f(t)\phi _{-t}|_{J}d\mu (t)\in {\F}$. So $I_{{\F}
}(\phi )=L_{w}^{1}(G)$ and $sp_{{\F} }(\phi )=\emptyset$.
Conversely, if $sp_{{\F} }(\phi )=\emptyset $ then $(\phi _{t}\ast
f)|_{J}\in {\F} $ for all $t\in G$ and $f\in L_{w}^{1}(G).$ By
Lemma 6.2, $\phi _{t}$ has approximate units and so $\phi
_{t}|_{J}\in {\F}$. This proves (v).

Now assume $\Delta _{t}^{N+1}\phi _{s}|_{J}\in {\F} $ for all
$s,t\in G.$ If $g\in L_{w}^{1}(G)$ then

$\qquad \qquad (\phi \ast \Delta
_{t}^{N+1}g)|_{J}=\int_{G}g(s)(\Delta _{t}^{N+1}\phi
_{-s})|_{J}d\mu (s)\in {\F} $

and so $\Delta _{t}^{N+1}g\in I_{{\F} }(\phi ).$
But $(\Delta _{t}^{N+1}g)\widehat{}(\gamma )=(\gamma (t)-1)^{N+1}\widehat{g}%
(\gamma )$ is zero for all $t\in G$ and $g\in L_{w}^{1}(G)$ only when $%
\gamma =1.$ So $sp_{{\F} }(\phi )\subseteq \{1\}.$ Conversely, if $%
sp_{{\F} }(\phi )\subseteq \{1\}$ then, by Theorem 6.1, $\{\Delta
_{t}^{N+1}g:g\in L_{w}^{1}(G),t\in G\}\subseteq I_{{{\F}} }(\phi
_{s})$ and so $(\Delta _{t}^{N+1}\phi _{s}\ast g)|_{J}=(\phi
_{s}\ast \Delta
_{t}^{N+1}g)|_{J}\in {\F}$. Taking approximate units we conclude $%
\Delta _{t}^{N+1}\phi _{s}|_{J}\in {\F}$. This proves (vi).

(b) We  only  prove (v$''$) as the other parts are similar to (a).
Since $sp_{{\A} }(\phi )=\emptyset $, $I_{\A}\phi$ is dense in
$L_{w}^{1}(G)$. This implies that for each $\lambda \in \widehat
{G}$ there is $f\in I_{\A}\phi$ such that $\hat {f} (\lambda)\not
= 0$, $\hat{f} \ge 0$. Since $sp_w(\phi)$ is compact  and
$I_{\A}\phi$ is linear, one can construct $F \in I_{\A}\phi$ such
that $\widehat {F} (\lambda)>0$ for all $\lambda \in W$, where $W$
is a compact neighbourhood of $sp_{w_N}(\phi)$. By  Corollary p.
19 in \cite{RH}, there is $H\in L_{w}^{1}(G)$ such that $\widehat
{F*H}(\lambda)=1$, $\lambda \in W$. This implies $\phi=F*H*\phi\in
\A$.
\end{proof}

\smallskip

\textbf{Corollary 6.4.} Assume $w$ has polynomial growth $N.$ If
$\phi \in BUC_{w}(G,X),$ $sp_{{{\F} }}(\phi )\subseteq \{1\}$ and $%
\phi |_{J}\in E_{w}^{0}(J,X)\ $ then $\phi |_{J}\in {{\F}}
+C_{w,0}(J,X).$ If also $w=1,$ then $\phi |_{J}\in {\F}$.

\begin{proof} By Proposition 6.3(v) and (vi) , $\Delta
_{t}^{N+1}\phi |_{J}\in {\F} $ for all $t\in J.$ By Lemma 4.4,
$\Delta _{t}^{j}\phi |_{J}\in E_{w}^{0}(J,X)$ for $0\leq j\leq N$
and all $t\in J.$ Therefore we can apply the difference Theorem
4.12 repeatedly to obtain the result. \end{proof}

 \smallskip

We denote by $E^{\prime }$ the set of accumulation points of a
subset $E$ of a topological space. Recall that such a set $E$ is
called \textit{perfect} if $E^{\prime }=E$ and \textit{residual}
if it is closed and has no non-empty perfect subsets. Loomis
\cite{LH} proved that if $\phi \in BUC(G,\mathbf{C})$ and $sp(\phi
)$ is residual, then $\phi \in AP(G,\mathbf{C})$.  Basit (see [43,
pp.92, 97]) and Baskakov  (\cite{BAG})
 proved the same result for $X$-valued functions provided $%
c_0\not\subset X$. Basit \cite{BBA} proved that for certain classes $%
{{\F}} \subseteq BUC(J,X),$ where $J\in \{\mathbf{R_+,R}\},$ if $%
sp_{{\F} }(\phi )$ is residual and $\gamma ^{-1}\phi |_{J}\in
E(J,X)$ for all $\gamma \in sp_{{\F} }(\phi )$ then $\phi \in
{\F}$. Some special cases of his results were proved by Ruess and Ph%
\'{o}ng \cite{RP}. Recently, Arendt and Batty \cite{ABj} gave an
operator theoretic proof for the case ${{\F}} =AP(\mathbf{R},X).$
We extend these results to unbounded functions as follows. (It
will be seen in Example 7.13 that for general weights $w$ the
ergodicity condition in Theorem 6.5 cannot be replaced by
$c_0\not\subset X$).

\bigskip

\textbf{Theorem 6.5.} Assume $w$ has polynomial growth. If $\phi
\in BUC_{w}(G,X),$ then $sp_{{{\F}} }(\phi )$ contains no isolated
points whenever one of the following holds:

(i) ${{\F}} $ is a \ $\Lambda _{w}^{0}$-class and $\gamma
^{-1}\phi |_{J}\in E_{w}^{0}(J,X)$ for all $\gamma \in sp_{{{\F}}
}(\phi )$;

(ii) ${{\F}} $ is a \ $\Lambda _{w}$-class and $\gamma ^{-1}\phi
|_{J}\in E_{w}(J,X)$ for all $\gamma \in sp_{{{\F}} }(\phi )$.

\noindent If also $sp_{{{\F} }}(\phi )$ is residual, then $\phi
|_{J}\in {{\F}}$.

\,

\begin{proof} Suppose $\gamma $ is an isolated point of
$sp_{{\F}
}(\phi ).$ Let $\gamma ^{-1}\phi |_{J}$ have $w$-mean $p|_{J}$ where $%
p\in P_{w}(G,X)$. Take an open neighbourhood $U$ of $\gamma $ in
$\widehat{G} $ such that $U\cap sp_{{{\F}} }(\phi )=\{\gamma \}.$
Choose $f\in
L_{w}^{1}(G)$ such that $\widehat{f}(\gamma )\neq 0$ and supp$(\widehat{f}%
)\subseteq U.$ Then $sp_{{{\F} }}(\phi \ast f)\subseteq \{\gamma
\}$ and so $sp_{{{\F}} }(\gamma ^{-1}(\phi \ast f))\subseteq
\{1\}$. By Proposition 4.10, the restriction to $J$\ of $\gamma
^{-1}(\phi \ast
f)=(\gamma ^{-1}\phi )\ast (\gamma ^{-1}f)$ is $w$-ergodic with $w$-mean $%
(p\ast \gamma ^{-1}f)|_{J}$. By Corollary 6.4, $(\gamma ^{-1}(\phi
\ast f)-(p\ast \gamma ^{-1}f))|_{J}\in {\F}$. Hence $(\phi \ast
f)|_{J}\in {\F} $ which means $\gamma \notin sp_{{{\F} }}(\phi ).$
This is a contradiction and so $sp_{{\F} }(\phi )$ contains no
isolated points. A residual perfect set is empty and so
Proposition 6.3(v) gives the final assertion. \end{proof}

\smallskip

The following result was proved in  Theorem
3.4 of \cite{BP} under a slightly stronger assumption than (2.4) and with $X=%
\mathbf{C}$. The same proof is valid under the present
assumptions. See also  Proposition 0.5 of \cite{PJ}.

\bigskip

\textbf{Theorem 6.6}. Assume $w$ has polynomial growth and $\phi
\in BC_{w}(G,X).$ Then $sp_{w}(\phi )=\{\gamma _{1},...,\gamma
_{n}\}$ if and
only if $\phi =\sum_{j=1}^{n}\gamma _{j}p_{j}$ for some non-zero $%
p_{j}\in P_{w}(G,X).$

\bigskip

\textbf{Corollary 6.7}. Assume $w$ has polynomial growth, $%
{{\F}} $ is a $\Lambda _{w}$-class and $\phi \in BUC_{w}(G,X).$ Then $%
sp_{{{\F}} }(\phi )\subseteq sp_{w}(\phi )^{\prime }.$

\,

\begin{proof} Let $\gamma $ be an isolated point of $sp_{w}(\phi )$ or $%
\gamma \in \widehat{G}$ $\backslash \ sp_{w}(\phi ).$ Then there
is a compact neighbourhood $V$ of $\gamma $ such that $V\cap
sp_{w}(\phi
)\subseteq \{\gamma \}.$ Choose $f\in L_{w}^{1}(G)$ such that $\widehat{f}%
(\gamma )\neq 0$ and supp$(\widehat{f})\subseteq V.$ Then
$sp_{w}(\phi \ast f)\subseteq \{\gamma \}$ and so, by Theorem 6.6,
$(\phi \ast f)|_{J}\in TP_{w}(J,X)\subseteq {{\F} }.$ Hence
$\gamma \notin sp_{{\F} }(\phi )$ and so $sp_{{{\F}} }(\phi
)\subseteq sp_{w}(\phi )^{\prime }.$  \end{proof}

\smallskip

\textbf{Corollary 6.8.} Assume $w$ has polynomial growth and $\phi
\in BUC_{w}(G,X)$. If $sp_{w}(\phi )$ is residual and $\gamma
^{-1}\phi
\in E_{w}(G,X)$ for all $\gamma \in sp_{w}(\phi )^{\prime }$, then $%
\phi \in AP_{w}(G,X).$

\,

\begin{proof} By Proposition 5.7, ${{\F}} =AP_{w}(G,X)$ is a
$\Lambda _{w}$-class. By Corollary 6.7, $sp_{{{\F}} }(\phi
)\subseteq
sp_{w}(\phi )^{\prime }$ which is residual. But by Theorem 6.5, $%
sp_{{{\F}} }(\phi )$ has no isolated points and therefore is a
perfect set. A residual perfect set is empty and so, by
Proposition 6.3(v), $\phi \in {\F}$.  \end{proof}

 \smallskip

Since a discrete set is residual without limit points, we obtain
the following result immediately. For the case $w=1$ see \cite
{LH} and references therein.

\bigskip

\textbf{Corollary 6.9}. Assume $w$ has polynomial growth. If $\phi
\in BUC_{w}(G,X)$ and $sp_{w}(\phi )$ is discrete, then $\phi \in
AP_{w}(G,X).$

\bigskip

The following two results will be used to establish $w$-ergodicity from a
knowledge of spectra. We say that a function $f\in L_{w}^{1}(G)$ is of $w$-%
\textit{spectral synthesis} with respect to a closed subset $\Delta $ of $%
\widehat{G}$ if there is a sequence $(f_{n})$ in $L_{w}^{1}(G)$ converging
to $f$ and satisfying $\widehat{f_{n}}=0$ on a neighbourhood $U_{n}$ of $%
\Delta$.

\bigskip

\textbf{Theorem 6.10}. If $\phi \in BUC_{w}(G,X)$ and $1\notin
sp_{w}(\phi )$, then $\phi \in E_{w}^{0}(G,X).$

\,

\begin{proof} Take a compact neighbourhood $V$ of $0$ in $\widehat{G}$
with $V\cap sp_{w}(\phi )=\emptyset$. Take $f\in L_{w}^{1}(G)$ with $%
\widehat{f}(1)=1$ and supp $(\widehat{f})\subseteq V.$ Then
$sp_{w}(\phi \ast f)=\emptyset $ so $\phi \ast f=0.$ Moreover, $f$
is continuous.
Given $\varepsilon >0,$ choose a compact set $K$ in $G$ such that $\int_{G%
 \setminus K}\left| f(t)\right| w(t)d\mu (t)<\varepsilon
/(1+2\left\| \phi \right\| _{\infty }).$ For $s\in G$ define
$g(s)=(\phi -\phi _{-s})f(s).$ Hence $\int_{G}g(s)d\mu (s)=\phi
-\phi \ast f=\phi$. Moreover, by Lemma 4.4, $g(s)\in
E_{w}^{0}(G,X)$ and since $\phi $ is $w$-uniformly continuous, by (2.8), $%
g:G\rightarrow E_{w}^{0}(G,X)$ is continuous. Since $K$ is compact, $g|_{K}$
is separably-valued and hence Bochner integrable. Therefore $%
\int_{K}g(s)d\mu (s)\in E_{w}^{0}(G,X).$ But $\left\| \phi
-\int_{K}g(s)d\mu (s)\right\| _{w,\infty }\leq \left\|
\int_{G\setminus K}g(s)d\mu (s)\right\| _{w,\infty }<\varepsilon $
and so $\phi \in E_{w}^{0}(G,X)$ as claimed.  \end{proof}

 \smallskip

Finally we establish relationships between $w$-spectral synthesis,
minimal prime ideals and ergodicity. See \cite{G} for
$G=\mathbf{R}$. For $\gamma \in \widehat{G},$ let $J_{w}(\gamma )$
denote the closed span of $\{\gamma \Delta _{t}f:f\in
L_{w}^{1}(G),$ $t\in G^{N+1}\}$.

 For  a closed subset $\Delta \subset \widehat{G}$ set

\smallskip

$I_{w}(\Delta )= \{f\in L_{w}^{1}(G):\widehat{f} (\gamma)= 0$  for
all  $\gamma \in \Delta \}$,

\smallskip

$R_{w}(\Delta )=$ the closure of   $\{f\in I_{w}(\Delta ): $
  supp $\widehat{f}\,\cap \Delta\,$ is
residual $\}$.
\smallskip

\smallskip

  $S_{w}(\Delta )=$ the closure
of $\{f\in L_{w}^{1}(G):\widehat{f}$ is $0$ on a neighbourhood of
$\Delta \}$, the space of functions in $L_{w}^{1}(G)$ which are of $w$%
-spectral synthesis with respect to $\Delta$,

\smallskip

\noindent One can   verify

\qquad\qquad\qquad $S_{w}(\Delta ) \subset R_{w}(\Delta ) \subset
I_{w}(\Delta )$.

\bigskip

\textbf{Theorem 6.11}. Suppose $w$\ has polynomial growth of order $N$ and $%
\gamma \in \widehat{G}.$

(a) $J_{w}(\gamma )$ is the minimal closed ideal of $L_{w}^{1}(G)$ with
cospectrum $\{\gamma \}.$

(b) $S_w(\gamma )=J_{w}(\gamma ).$

\,

\begin{proof} (a) Since $(\gamma \Delta _{t}f)\ast g=\gamma
\Delta _{t}(f\ast \gamma ^{-1}g),$ $J_{w}(\gamma )$ is a closed
ideal. Minimality follows from Theorem 6.1.

(b) Since $J_{w}(\gamma )=\gamma J_{w}(1)$ and $S_w(\gamma
)=\gamma
S_w(1),$ it suffices to prove $S_w(1)=J_{w}(1).$ For $k\in \mathbf{Z}%
_{+},$ let $J_{w}^{k}$ denote the closed span of $\{\Delta _{t}f:f\in
L_{w}^{1}(G),$ $t\in G^{k}\}.$ For $g\in L_{w}^{1}(G)$ satisfying $\widehat{g%
}(1)=1,$\ let $L_{g}$ be the operator on $L_{w}^{1}(G)$ defined by $%
L_{g}f=f-f\ast g=-\int_{G}(\Delta _{-s}f)g(s)d\mu (s).$ The integral is an
absolutely convergent Bochner integral and so $L_{g}$ maps $J_{w}^{k}$ into $%
J_{w}^{k+1}.$ Now take any $f\in L_{w}^{1}(G)$ with
$\widehat{f}=0$ on a neighbourhood $U$\ of $1.$ Choose $g\in
L_{w}^{1}(G)$ with $\widehat{g}(1)=1$ and
supp$(\widehat{g})\subseteq U.$ So $f=L_{g}^{N+1}f\in
J_{w}^{N+1}=J_{w}(1).$ Hence, $S_w(1)\subseteq J_{w}(1).$ But
$S_w(1)$ is an ideal with cospectrum $\{1\}$ and so by (a),
$S_w(1)=J_{w}(1)$.  \end{proof}

 \smallskip

\textbf{Corollary 6.12}. If $f\in J_{w}(\gamma )$ for some $\gamma
\in \widehat{G}$ and $\phi \in BUC_{w}(G,X),$\ then $\gamma
^{-1}(\phi \ast f)\in E_{w}^{0}(G,X).$

\,

\begin{proof} Let $h=\gamma \Delta _{t}g$ where $t\in G$ and
$g\in L_{w}^{1}(G).$ By Lemma 4.4, $\gamma ^{-1}(\phi \ast
h)=\Delta _{t}(\gamma ^{-1}\phi \ast g)\in E_{w}^{0}(G,X)$. Since
$f$ is in the closed linear span of such functions $h$ and
$E_{w}^{0}(G,X)$ is complete, the result follows. \end{proof}

 \smallskip

\textbf{Corollary 6.13}. If $f\in L_{w}^{1}(G)$ is of $w$-spectral
synthesis with respect to a closed subset $\Delta $ of $\widehat{G},\ $then $%
\gamma ^{-1}(\phi \ast f)\in E_{w}^{0}(G,X)$ for all $\gamma \in
\Delta $ and $\phi \in BUC_{w}(G,X).$

\begin{proof} Choose a sequence $(f_{n})$ in $L_{w}^{1}(G)$ converging to $%
f$ and satisfying $\widehat{f_{n}}=0$ on a neighbourhood $U_{n}$
of $\Delta . $ For $\gamma \in \Delta $\ we have $\gamma
^{-1}(\phi \ast f_{n})=(\gamma ^{-1}\phi )\ast (\gamma
^{-1}f_{n})$ and so $sp_{w}(\gamma ^{-1}(\phi \ast
f_{n}))\subseteq $ supp$(\widehat{\gamma ^{-1}f_{n}}).$ Also
$(\widehat{\gamma ^{-1}f_{n}})(\tau )=\widehat{f_{n}}(\gamma \tau
)$
and so $1\notin $ supp$(\widehat{\gamma ^{-1}f_{n}}).$ By Theorem 6.10, $%
\gamma ^{-1}(\phi \ast f_{n})\in E_{w}^{0}(G,X).$ Hence $\gamma
^{-1}(\phi \ast f)\in E_{w}^{0}(G,X)$  for all $\gamma \in \Delta
$.
\end{proof}

\bigskip

\textbf{7. Derivatives and indefinite integrals}

Throughout this section we assume that $%
J\in \{\mathbf{R_+,R, Z_+,Z}\}$ and ${{\F}} = \F(J,X) \subset
BUC_{w}(J,X)$  satisfying (5.2), (5.3) and (5.5). Examples of such classes are $%
C_{w,0}(J,X)$,   $AP_{w}(G,X)$ with $G\in \{\mathbf {R, Z}\}$
 and  if $w$ has polynomial growth $BUC_{w}(J,X)\cap E_{w}(J,X)$.

\bigskip

\textbf{Proposition 7.1.} Let $%
J\in \{\mathbf{R_+,R}\}$.

(a) Let $w$\ have polynomial growth. If $\phi \in
BC_{w}(\mathbf{R},X)$
and $sp_{w}(\phi )$ is compact, then $\phi ^{(j)}\in BUC_{w}(\mathbf{R}%
,X)$ for all $j\geq 0.$

(b) If $\phi \in {{\F}} $ and $\phi ^{\prime }$ is $w$%
-uniformly continuous,\ then $\phi ^{\prime }\in {{\F}}$.

(c) If $\phi \in BC_{w}(J,X)$ and $\phi ^{\prime }$ is
$w$-uniformly continuous,\ then $\phi ^{\prime }\in
E_{w}^{0}(J,X)\cap BUC_{w}(J,X).$

(d) If $\phi ,\phi ^{\prime }\in BC_{w}(\mathbf{R},X)$ then
$sp_{w}(\phi ^{\prime })\subseteq sp_{w}(\phi )\subseteq
sp_{w}(\phi ^{\prime })\cup \{1\}.$

\,

\begin{proof} (a) Choose $f\in S(\mathbf{R}),$ the Schwartz
space of rapidly decreasing functions, such that $f$ has compact
support and is $1$
on a neighbourhood of $sp_{w}(\phi ).$ Then $f^{(j)}\in L_{w}^{1}(\mathbf{%
R})$ for all $j\geq 0.$ Moreover, $\phi =\phi \ast f$ and so $\phi
^{(j)}=\phi \ast f^{(j)}$ for all $j\geq 0.$ Hence $\phi ^{(j)}\in
BUC_{w}(\mathbf{R},X).$

(b) If $\psi _{n}=n\Delta _{1/n}\phi $ then $\psi _{n}\in {{\F}}
.$ Moreover, by the $w$-uniform continuity of $\phi ^{\prime },$
given $\varepsilon >0$ there exists $n_{\varepsilon }$ such that
$\left\| \psi _{n}(t)-\phi ^{\prime }(t)\right\| =\left\|
n\int_{0}^{1/n}(\phi
^{\prime }(t+s)-\phi ^{\prime }(t))ds\right\| <\varepsilon w(t)$ for all $%
t\in J$ and $n>n_{\varepsilon }.$ Hence $\phi ^{\prime }\in {{\F}}
.$

(c) With the notation of the proof of (b), $\psi _{n}\in
E_{w}^{0}(J,X)\cap BUC_{w}(J,X)$ by Lemma 4.4. Hence, so does
$\phi ^{\prime }.$

(d) For any $f\in L_{w}^{1}(\mathbf{R})$ we have $(\phi \ast
f)^{\prime }=\phi ^{\prime }\ast f$ and so $I_{w}(\phi ^{\prime
})\supseteq I_{w}(\phi ).$ Hence, $sp_{w}(\phi ^{\prime
})\subseteq sp_{w}(\phi ).$ For the second inclusion, let
$c=\sup_{0\leq s\leq 1}w(s). $ By Proposition 2.2, $\sup_{n\leq
t\leq n+1}w(t)\leq cw(n)<ce^{n}$
for all $n$ sufficiently large. Hence, if $g(t)=\exp (-t^{2})$ then $%
g,g^{\prime }\in $ $L_{w}^{1}(\mathbf{R}).$ Now take $\gamma \in
\mathbf{R} \setminus (sp_{w}(\phi ^{\prime })\cup \{1\}).$ So
$\gamma (t)=e^{ist}$ for some $s\neq 0.$ Choose $f\in
L_{w}^{1}(\mathbf{R})$ such
that $\phi ^{\prime }\ast f=0$ but $\widehat{f}(\gamma )\neq 0.$ Let Let $%
h=f\ast g^{\prime }.$ Then $\phi \ast h=\phi \ast f\ast g^{\prime
}=\phi ^{\prime }\ast f\ast g=0$ whereas $\widehat{h}(\gamma )=$ $is%
\widehat{f}(\gamma )\widehat{g}(\gamma )\neq 0.$ So $\gamma \notin
sp_{w}(\phi ^{\prime })$ showing $sp_{w}(\phi )\subseteq
sp_{w}(\phi ^{\prime })\cup \{1\}.$ \end{proof}

\smallskip

\textbf{Proposition 7.2}. (a) If $J\in \{\mathbf{R_+, R}\}$, $\phi \in {{\F}}(J,X) $ and $%
P\phi $   is $w$-ergodic with $w$-mean $p,$ then $P\phi -p$
 $\in {{\F}}(J,X)$.

(b) If $J\in\{\mathbf{ Z_+, Z}\}$, $\phi \in {{\F}}(J,X) $ and $%
S\phi $   is $w$-ergodic with $w$-mean $p,$ then $S\phi -p \in
{{\F}}(J,X)$.

\begin{proof} (a) We may assume that $\phi \in BUC_w (\mathbf{R},X)$. For $t\in \mathbf{R}$ set $\chi _{t}=\chi _{[
-t,0] }$ if $t\geq 0$ and $\chi _{t}=-\chi _{[ 0,-t] }$ if
$t<0.$ Then $\Delta _{t}P\phi =\phi \ast \chi _{t}=\int_{\mathbf{R}%
}\phi _{-s}\chi _{t}(s)ds.$ Since $\phi \in BUC_{w}(\mathbf{R},X)$
the integral converges as a Lebesgue-Bochner integral and
therefore by (5,2), (5.3) $\Delta _{t}P\phi \in {{\F}}$. The
result follows from the difference Theorem 4.12 and (5.5).

(b)  Note that $\Delta_1 S\phi=\phi \in {{\F}} (J,X)$ and
therefore by (5,2), (5.3) $\Delta_k S\phi \in {{\F}} (J,X)$ for
each $k\in J$. The result follows from the difference Theorem 4.12
and (5.5).
\end{proof}

\smallskip

\textbf{Proposition 7.3}. (a) If  $\phi \in {{\F}}
=AP_{w}(\mathbf{R},X)$ and $P\phi  \in BUC_{w}(\mathbf{R},X),$
then $sp_{{{\F}} }(P\phi ) \subseteq \{1\}.$

(b) If $\phi \in {{\F}} =AP_{w}(\mathbf{Z},X)$ and $S\phi  \in
BUC_{w}(\mathbf{Z},X),$ then $sp_{{{\F}} }(S\phi ) \subseteq
\{1\}.$

\,

\begin{proof} (a) Let $s,t\in \mathbf{R}$. With $\chi _{s}\ $as in
the previous proof, $(\Delta _{s}P\phi )_{t}=\phi _{t}\ast \chi
_{s}$ and by Proposition 5.6, $(\Delta _{s}P\phi )_{t}\in {{\F}}$.
 By Proposition 6.3(vi), $sp_{{{\F}} }(P\phi )\subseteq \{1\}.$

 (b) This case  follows similarly as $\Delta_k S\phi \in
{{\F}} $ for all $k\in \mathbf{ Z}$.
\end{proof}

\smallskip

\textbf{Remark 7.4.} For general weights $w$ the ergodicity
condition in Theorem 7.5 cannot be replaced by $ c_{0}\not\subset
X.$ Indeed, let $G=\mathbf{R,}$ $w(t)=1+\left| t\right| $ and
$\phi (t)=\frac{t\cos \log w(t)}{w(t)}.$ So $\phi \in
C_{w,0}(\mathbf{R},\mathbf{C})\subseteq
{{\F}} =AP_{w}(\mathbf{R},\mathbf{C})$ and $P\phi (t)=\frac{w(t)}{2}%
[\cos \log w(t)+\sin \log w(t)] -\sin \log w(t)-\frac{1}{2}.$ So
$P\phi \in BUC_{w}(\mathbf{R,C).}$ By Proposition 7.3, $sp_{{{\F}}
}(P\phi )\subseteq \{1\}$ and so $sp_{{{\F}} }(P\phi )$\ is
residual. However, $P\phi \notin AP_{w}(\mathbf{R},\mathbf{C}).$ Indeed, $%
P\phi \notin E_{w}(\mathbf{R},\mathbf{C}).$ For if $\frac{P\phi -ct}{w}%
\in E^{0}(\mathbf{R},\mathbf{C})$ for some constant $c$, then because $(%
\frac{P\phi -ct}{w}|_{\mathbf{R}_{+}})^{\prime }=(\frac{\phi }{w}-\frac{%
P\phi }{w^{2}}-\frac{c}{w^{2}})|_{\mathbf{R}_{+}}\in C_{0}(\mathbf{R}_{+},%
\mathbf{C)}$ it would follow from Theorem 4.12 that $%
\frac{P\phi -ct}{w}|_{\mathbf{R}_{+}}\in
C_{0}(\mathbf{R}_{+},\mathbf{C).}$

\noindent See also Example 7.13.
\bigskip

As before we will take $w_m (t)= (1+ |t|)^m$ for $m\in
\mathbf{N}$.

\bigskip

\textbf{Theorem 7.5}.  Assume $\phi\in AP_w(\mathbf{R},X) $
respectively  $\phi\in AP_w(\mathbf{Z},X) $ has $w$-mean $p$. Then
$P^m(\phi-p)\in C_{ww_m,0} (\mathbf{R},X)$ respectively
$P^m(\phi-p)\in C_{ww_m, 0}(\mathbf{Z},X) $ for each $m \in
\mathbf{N}$.

\begin{proof} The assumptions imply $\phi-p \in E_w^0(\mathbf{R},X)$
 respectively $\phi-p \in E_w^0(\mathbf{Z},X)$. The
result follows from Theorem 4.14. \end{proof}

\smallskip

\textbf{Lemma 7.6}. Let $J\in \{\mathbf{Z_+, Z, R_+, R}\}$.  For
natural numbers $m,N$\ and$\ $non-negative integers $j,k$ set
$a(m,j)=(-1)^{j}\binom{N}{j}\binom{m-1+j}{j}j!.$

\begin{enumerate}

\item[(a)] $S^{m}(t^{N}\phi)=\sum_{j=0}^{N}a(m,j)\,t^{N-j}S^{m+j}\phi  +
\sum_{j=1}^{N}\,t^{N-j}\sum_{k=0}^{j-1} c(m,k,j)\,S^{m+k}\phi$ for
some choice of $ c(m,k,j)$, $0\le k <j \le N$, any $\phi \in X^J$,
$J\in \{\mathbf{Z_+, Z}\}$.

\item[(b)] $P^{m}(t^{N}\phi)=\sum_{j=0}^{N}a(m,j)\,t^{N-j}P^{m+j}\phi$
for any $\phi \in L^1_{loc} (J,X)$, $J\in \{\mathbf{R_+, R}\}$.

\item[(c)] $\sum_{j=0}^{N}\,\frac{a(m,j)}{(j+k)!}=\left\{
\begin{array}{cc}
\binom{N+k-m}{N}\frac{N!}{\left( N+k\right) !} & \text{if }m\leq k \\
0 & \text{if }k+1\leq m\leq k+N \\
(-1)^{N}\binom{m-k-1}{N}\frac{N!}{\left( N+k\right) !} & \text{if }m>k+N.%
\end{array}%
\right.$.

\end{enumerate}

\begin{proof}
(a) For $N=0$ the claim is trivial and for $N=1$ one can prove by
induction on $m$ that $S^m (t\phi)= t S^m \phi- m
(S^{m+1}\phi)_1$. The general case is then proved by induction on
$N$ using $S^m (t^{N+1}\phi)= t S^m t^N\phi- m (S^{m+1}t^N
\phi)_1$.

(b) Follows similarly as in (a) using $P$ instead of $S$.

(c) For $m\geq k+1$ we have

$\sum_{j=0}^{N}\,\frac{a(m,j)}{(j+k)!}=\sum_{j=0}^{N}(-1)^{j}\binom{N}{j}%
\binom{m-1+j}{j}\ \frac{j!}{(j+k)!}=\frac{1}{(m-1)!}\sum_{j=0}^{N}(-1)^{j}%
\binom{N}{j}\ \frac{(m-1+j)!}{(k+j)!}$

\ \ \ \ \ \qquad \qquad $=\ \ \frac{1}{(m-1)!}\sum_{j=0}^{N}(-1)^{j}\binom{N%
}{j}\ D^{m-k-1}t^{m+j-1}|_{t=1}$

\ \ \ \ \ \qquad \qquad $=$\ $\frac{1}{(m-1)!}D^{m-k-1}\
t^{m-1}\sum_{j=0}^{N}(-1)^{j}\binom{N}{j}\ t^{j}|_{t=1}$

\ \ \ \ \ \qquad \qquad $=\frac{1}{(m-1)!}D^{m-k-1}\
t^{m-1}(1-t)^{N}|_{t=1}. $

For $m-k-1<N$ this last expression is $0$ and for $m-k-1\geq N$ it
is

$\ \ \ \ \ \ \ \ \ \ \
\frac{1}{(m-1)!}\binom{m-k-1}{N}(D^{m-k-1-N}\
t^{m-1})D^{N}(1-t)^{N}|_{t=1}=\frac{1}{(m-1)!}\binom{m-k-1}{N}\frac{(m-1)!}{%
(N+k)!}(-1)^{N}N!$

\noindent as claimed.

 For $m\leq k$ the claim follows readily by
substituting $\phi (t)=t^{k-m}$ in (b).

\end{proof}

Our main result is the following:

\ \ \ \

\textbf{Theorem 7.7}. Assume $\phi \in AP(\mathbf{Z},X)$ and $\sum_{j=0}^{N}b_{j}%
\,t^{N-j}S^{j+1}\phi \in BUC_{w_{N}}(\mathbf{Z},X)$ for some
$b_{j}\in \mathbf{C}, b_{0}\neq 0$.

\begin{enumerate}
\item[(a)] $S\phi \in BUC(\mathbf{Z},X)$ and if $\sum_{j=0}^{N}\frac{b_{j}%
}{(j+1)!}\neq 0$ then $M\phi =0.$

\item[(b)] If $c_0 \not\subset X$  then \ $S(\phi -M\phi )\in AP(%
\mathbf{Z},X).$
\end{enumerate}

\begin{proof}
Let $a=M\phi $ and $\psi =\sum_{j=0}^{N}b_{j}\,t^{N-j}S^{j+1}\phi
.$\ Then $\psi =\sum_{j=0}^{N}b_{j}\,t^{N-j}S^{j+1}(\phi
-a)+\,t^{N+1}a\sum_{j=0}^{N}\frac{b_{j}}{(j+1)!} + a r$, where $r
\in P_N (\mathbf{Z},\mathbf {R})$. By Theorem 7.5 and the
assumption,
 $\psi
-\sum_{j=0}^{N}b_{j}\,t^{N-j}S^{j+1}(\phi -a) = t^{N+1}a\sum_{j=0}^{N}\frac{b_{j}}{(j+1)!} + a r \in C_{w_{N+1},0}(\mathbf{Z}%
,X)$ and so either $a=0$ or
$\sum_{j=0}^{N}\frac{b_{j}}{(j+1)!}=0$. To prove
the rest of the theorem, we may assume $a=0$. By Theorem 7.5, $%
S^{j}\phi (t)/w_{j}(t)\rightarrow 0$ as $t\rightarrow \infty $, $1\le j\le N$. Since $%
\phi $ is almost periodic we may choose $(t_{n})\subset
\mathbf{Z}$ such that $t_{n}\rightarrow \infty $ and $\phi
_{t_{n}}\rightarrow \phi $
uniformly on $\mathbf{Z}$. Moreover, as $M\phi =0,$ by Theorem 7.5, $%
S^{j}\phi (s+t_{n})/w_{j}(s+t_{n})\rightarrow 0$, $s\in
\mathbf{Z}$, $1\le j\le N$.  Given $x^{\ast }\in X^{\ast }$, it
follows that $x^{\ast }\circ S\phi _{t_{n}}\rightarrow x^{\ast
}\circ S\phi $ locally uniformly. Moreover, by passing to a
subsequence if necessary, we may assume $x^{\ast }\circ \psi
(t_{n})/w_{N}(t_{n})\rightarrow b$ for some $b\in \mathbf{C}$. By
Theorem
7.5 again, we obtain%
\begin{eqnarray*}
\psi (t+t_{n})
&=&\sum_{j=0}^{N}b_{j}\,(t+t_{n})^{N-j}[\sum_{s=0}^{t_n-1}(S^{j}\phi)
(s)+ \sum_{s=0}^{t-1}(S^{j}\phi)
(s+t_{n})] \\
&=&\psi (t_{n})+b_{0}t_{n}^{N}(S\phi_{t_n}) (t)+o(t_{n}^{N}).
\end{eqnarray*}%
Therefore $x^{\ast }\circ \psi (t+t_{n})/w_{N}(t+t_{n})\rightarrow $ $%
b+b_{0}x^{\ast }\circ S\phi (t)$ for each $t\in \mathbf{Z}$.
Hence, since $\psi /w_{N}$ is bounded, so too is $x^{\ast }\circ
S\phi $. Since $x^{\ast }$ is arbitrary, $S\phi $ is weakly
bounded and therefore bounded. Since $\Delta_1 S\phi =\phi \in
AP(\mathbf{Z},X)$, it follows $\Delta_k S\phi  \in
AP(\mathbf{Z},X)$, for $k\in \mathbf{Z}$. So if  $c_0 \not\subset
X$ then by a generalization of Kadet's Theorem \cite{BB} (see also
\cite{K}), $S\phi $ is almost periodic.

\end{proof}

\textbf{Corollary 7.8.} Assume $\phi \in AP(\mathbf{Z},X)$ and $%
S^m (t^{N}\phi )\in BUC_{w_{N}}(\mathbf{Z},X)$ for some $m\in
\mathbf{N}$.

\begin{enumerate}
\item[(a)] $S\phi, S^2\phi, \cdots, S^m\phi \in BUC(\mathbf{Z},X)$
and $M\phi =0, MS\phi =0, \cdots,  MS ^{m-1}\phi =0.$

\item[(b)] If $c_0 \not\subset X$  then \ $S\phi, S^2\phi, \cdots, S^m\phi \in AP(%
\mathbf{Z},X).$
\end{enumerate}

\begin{proof} Since $S^{m-1} t^N \phi =(S^m t^N\phi)_1 -S^m t^N\phi$ we conclude $ S^{m-1}t^N\phi \in
BUC_{w_n}(\mathbf{Z},X)$; in the same way $ S^{j} t^N\phi \in
BUC_{w_N}(\mathbf{Z},X)$ for all $j=1, \cdots, m $. We prove the
statement by induction on $m$.

Case $m=1$. By Lemma 7.6 (a),

$S (t^{N}\phi)=\sum_{j=0}^{N}a(1,j)\,t^{N-j}S^{1+j}\phi  +
\sum_{j=1}^{N}\, t^{N-j}\sum_{k=0}^{j-1} c(1,k,j)\, S^{N+k}\phi$.

\noindent Since $\phi$ is bounded, $\sum_{j=1}^{N}\,
t^{N-j}\sum_{k=0}^{j-1} c(1,k,j)\, S^{1+k}(\phi-M\phi)\in
BUC_{w_N}(\mathbf{Z},X)$. It follows from the assumptions
$\sum_{j=0}^{N}a(1,j)\,t^{N-j}S^{1+j}\phi \in
BUC_{w_{N}}(\mathbf{Z},X)$. So,
  by Theorem 7.7, $S\phi \in BUC (\mathbf{Z},X)$. By Lemma 7.6
(c), $\sum_{j=0}^N a(1,j)/(j+1)! \not =0$ and hence again by
Theorem 7.7,  $M\phi =0$ and $S\phi \in AP (\mathbf{Z},X)$.

Now assume that the statement is true for the case $m-1$ and let  $%
S^m (t^{N}\phi )\in BUC_{w_{N}}(\mathbf{Z},X)$. It follows $
S^{m-1}\phi \in BUC(\mathbf{Z},X)$. By Lemma 7.6 (a),

$S^{m-1}(t^{N}\phi)=\sum_{j=0}^{N}a(m-1,j)\,t^{N-j}S^{m-1+j}\phi +
\sum_{j=1}^{N}\,t^{N-j}\sum_{k=0}^{j-1} c(m-1,k,j)\,S^{m-1+k}\phi$
for some choice of $ c(m-1,k,j)$, $0\le k <j \le N$, any $\phi \in
X^J$, $J\in \{\mathbf{Z_+, Z}\}$.
 Hence

$S^{m}(t^{N}\phi)= S(S^{m-1}(t^{N}\phi)) = a(m-1,0) \, S (t^{N}
[S^{m-1}\phi]) +$

 $\sum_{j=1}^{N}a(m-1,j)\,S ( t^{N-j} S^j
[S^{m-1}\phi]) + \sum_{j=1}^{N}\, S(t^{N-j}\sum_{k=0}^{j-1}
c(m-1,k,j)\,S^k [S^{m-1}\phi])$.

\noindent Since $ S^{m-1}\phi \in BUC(\mathbf{Z},X)$,
 it follows

 $\sum_{j=1}^{N}a(m-1,j)\,S ( t^{N-j} S^j
[S^{m-1}\phi])$,

 $ \sum_{j=1}^{N}\, S(t^{N-j}\sum_{k=0}^{j-1}
c(m-1,k,j)\,S^k [S^{m-1}\phi]) \in BUC_{w_N}(\mathbf{Z},X)$.

 \noindent Hence $ S (t^{N} [S^{m-1}\phi]) \in
BUC_{w_N} (\mathbf{Z},X)$. Applying the case $m=1$ to $\psi =
S^{m-1}\phi$, one gets $S\psi \in BUC(\mathbf{Z},X) $, $M\psi=0$
and if $c_0 \not\subset X$, $S\psi \in AP(\mathbf{Z},X) $. This
completes the proof of the corollary.

\end{proof}

\textbf{Theorem 7.9}. Let  $\psi \in BUC_{w_1}(\mathbf{R},X)$. Set
$h_n=1/2^n$.

\begin{enumerate}
\item[(a)] If $\psi|\mathbf{Z}h_n\in AP_{w_1}(\mathbf{Z}h_n,X)$
for all  $n \in \mathbf{N}$ then $ \psi \in
AP_{w_1}(\mathbf{R},X)$.

\item[(b)] If $\psi_h|\mathbf{Z}\in AP_{w_1}(\mathbf{Z},X)$ for all $0\le h
\le 1$ then $ \psi \in AP_{w_1}(\mathbf{R},X)$.

\end{enumerate}

\,

\begin{proof} (a) Since $\psi \in BUC_{w_1}(\mathbf{R},X)$,
$\frac{\psi}{w_1} \in BUC (\mathbf{R},X)$. So  for each
$\varepsilon
> 0 $ there is
 $n=n(\varepsilon)$ such that  the range of $\frac{\psi}{w_1} |
\mathbf{Z}h_n$ is an $\varepsilon$-net for the range of
$\frac{\psi}{w_1} $. As $\psi|\mathbf{Z}h_n\in
AP_{w_1}(\mathbf{Z}h_n,X)$, the range of $\frac{\psi}{w_1}
|\mathbf{Z}h_n$ is relatively compact by Proposition 5.1. It
follows that the range of $\frac{\psi}{w_1}$ is relatively
compact. By Proposition 5.5 (c), $\psi |\mathbf{Z}h_n= f^n + t\,
g^n$, where $f^n \in C_{w_1,0}(\mathbf{Z}h_n,X)$ and $g^n \in
AP(\mathbf{Z}h_n,X)$. Choose
 $s_n \in \mathbf{Z}h_n$ such that $||g^n
_{s_n}- g^n||_{\infty}\le 1/n$ and $s_n\to \infty$ as $n \to
\infty$. It follows $||g^n _{s_n}- g^n||_{\infty}\to 0$ as $n\to
\infty$. By the Arzela-Ascoli theorem we can assume
$(\frac{\psi}{w_1})_{s_k}$ converges  to some $F \in
BUC(\mathbf{R},X)$  uniformly on each bounded interval of
$\mathbf{R}$. Since $\mathbf{Z}h_n \subset \mathbf{Z}h_{n+1}$ and
by Proposition 5.5 (c), $g^{n+1}$ extends $g^n$ for all $n \in
\mathbf{N}$, it follows $(g^k _{s_k}- g^k) |\, \mathbf{Z}h_n = g^k
_{s_k}|\, \mathbf{Z}h_n - g^n$ for all $k \ge n$. This implies
$||g^k _{s_k}|\, \mathbf{Z}h_n- g^n||_{\infty}\to 0$   as $k\to
\infty$ for all $n \in \mathbf{N}$. Hence
$(\frac{\psi}{w_1})_{s_k}|\, \mathbf{Z}h_n =(\frac{f^k + t\,
g^k}{w_1})_{s_k} \,|\, \mathbf{Z}h_n $ converges pointwise to  $
F| \mathbf{Z}h_n = g^n$ and so $ F\,|\, \mathbf{Z}h_n\in AP
(\mathbf{Z}h_n,X)$ for all $n\in \mathbf{N}$. As $ F\in BUC
(\mathbf{R},X)$, for each $\varepsilon>0$, one can choose $\delta
>0 $ such that $||F-F_h||_{\infty}\le \varepsilon$ for all $|h|\le \delta$.
Choose $h_n\le \delta $. Then each $\varepsilon$-period $\tau$ of
$g^n$ is a $3\varepsilon$-period for $F$. Moreover, $\tau+h$
 is a $4\varepsilon$-period for $F$, $|h|\le
\delta$ . This gives $F \in AP (\mathbf{R},X)$. Since
$(\psi-tF)\,| \mathbf{Z}h_n \in C_{w_1,0} (\mathbf{Z}h_n,X)$ and
$\frac{\psi-tF} {w_1}\in BUC(\mathbf{R},X)$, we conclude $ (\psi
-tF) \, \in C_{w_1,0}(\mathbf{R},X)$. This proves $\psi \in
AP_{w_1}(\mathbf{R},X)$ by  Proposition 5.5 (c).

(b) By  Proposition 5.5 (c), $\psi_{kh_n}| \mathbf{Z}= f^k + t
g^k$,  where $f^k \in C_{w_1,0}(\mathbf{Z},X) $, $g^k \in AP
(\mathbf{Z},X) $ and $k=0,1 ,\cdots, 2^n -1$. It follows

$E(\varepsilon, n)=\{ \tau\in \mathbf{Z}: ||g^k(m+\tau)-g^k
(m)||\le \varepsilon,\,\,\, k=0,1 ,\cdots, 2^n -1, \,\, m\in
\mathbf{Z}\}$

\noindent is relatively dense in  $\mathbf{Z}$ and therefore in
$\mathbf{Z} h_n = \cup_{k=0}^{2^n-1} [\mathbf{Z}+kh_n]$ for each
$\varepsilon > 0$. Define $g: \mathbf{Z}h_n \to X$ by

$g (t)= g^k (t-kh_n)$ if $t\in \mathbf{Z}+kh_n$, $k=0,1 ,\cdots,
2^n -1$.

\noindent Thus  $||g(t+\tau)-g (t)|| =
||g^k(t-kh_n+\tau)-g^k(t-kh_n)||\le \varepsilon$, $\tau \in
E(\varepsilon, n)$, $t\in \mathbf{Z}h_n$. It follows $g\in AP
(\mathbf{Z}h_n,X)$.  But $(\psi  - t\,g)(m+kh_n)= f^k(m)- kh_n \,
g^k(m)$ so  $\psi|, \mathbf{Z}h_n-tg \in C_{w_1,0}
(\mathbf{Z}h_n,X)$. By Proposition 5.5 (c), $\psi|\mathbf{Z}h_n
\in AP_{w_1}(\mathbf{Z}h_n,X)$. Since $n $ is arbitrary, (b)
follows from (a).
\end{proof}

\smallskip

\textbf{Theorem 7.10}. Let $\psi \in BUC_{w_1}(\mathbf{Z},X)$.

(a)If $\Delta_1\psi \in AP(\mathbf{Z},X)$ then $ \psi \in
AP_{w_1}(\mathbf{Z},X)$.

(b) If   $\Delta_1\psi \in AP(\mathbf{Z},X)+ t  AP(\mathbf{Z},X)$
then $ \psi \in AP_{w_1}(\mathbf{Z},X)$ provided  $c_0 \not\subset
X$.

\,

\begin{proof} (a) Let $\Delta_1\psi=f $ with $ f\in AP(\mathbf{Z},X)$.
Then $\psi = \psi(0)+Sf$. By Theorem 7.5, $S(f-Mf)\in
C_{w_1,0}(\mathbf{Z},X) $ and  hence $Sf\in
AP_{w_1}(\mathbf{Z},X)$, by Proposition 5.5.

(b) Let $\Delta_1\psi=f +t\, g$ with $ f, g\in AP(\mathbf{Z},X)$.
Then $\psi = \psi(0)+Sf+S(tg)$. By (a), $Sf\in
AP_{w_1}(\mathbf{Z},X)$. This and the  assumptions imply  $S
(t\,g) \in B_{w_1}(\mathbf{Z},X)$. So, by Corollary 7.8,  $S g\in
AP(\mathbf{Z},X)$. Therefore  $S(t\,g) = t\, Sg- S^2 g\in
AP_{w_1}(\mathbf{Z},X)$, by  Theorem 7.5 and Proposition 5.5 (c).
It follows $\psi\in AP_{w_1}(\mathbf{Z},X)$. \end{proof}

\smallskip

\textbf{Corollary 7.11}. Let $\psi \in BUC_{w_1}(\mathbf{R},X)$.

(a) If $\Delta_1\psi\in AP(\mathbf{R},X)$ then $ \psi \in
AP_{w_1}(\mathbf{R},X)$.

(b) If
 $\Delta_1\psi\in AP(\mathbf{R},X)+ t\, AP(\mathbf{R},X)$ then $ \psi
\in AP_{w_1}(\mathbf{R},X)$ provided $c_0 \not\subset X$.

\,

\begin{proof}
(a) Let $\Delta_1\psi=f$. It follows $\Delta_1 (
\psi_t|\mathbf{Z})= f_t |\, \mathbf{Z})\in AP(\mathbf{Z},X)$, $0
\le t < 1$.
 By Theorem 7.10 (a), $\psi_t|\mathbf{Z} \in AP_{w_1}(\mathbf{Z},X)$, $0 \le t < 1$
 and by Theorem 7.9 (b), $\psi\in AP_{w_1}(\mathbf{R},X)$.

(b) Follows similarly as in (a).
\end{proof}

\smallskip

 \textbf{Corollary 7.12}.

(a) If $\phi \in AP(\mathbf{R},X)$, then $ P\phi\in
AP_{w_1}(\mathbf{R},X)$.

(b) If $\phi \in AP(\mathbf{R},X)$ and $P(t\phi)\in
BUC_{w_1}(\mathbf{R},X)$, then  $ P(t\phi)\in
AP_{w_1}(\mathbf{R},X)$ provided $c_0 \not\subset X$.

(c) If $\phi \in t\cdot AP(\mathbf{R},X)+  AP(\mathbf{R},X)$ and
$P\phi\in BUC_{w_1}(\mathbf{R},X)$, then  $ P\phi\in
AP_{w_1}(\mathbf{R},X)$ provided $c_0 \not\subset X$.

\begin{proof} (a) Clearly $P\phi \in BUC_{w_1}(\mathbf{R},X)$ and $\Delta_ 1 P\phi\in
AP(\mathbf{R},X)$, therefore  $ P\phi\in AP_{w_1}(\mathbf{R},X)$
by Corollary  7.11 (a).

(b) Follows from Corollary 7.11 (b) since $\Delta_1 P(t\phi)=
t\int_0^1 \phi(\cdot+s)\, ds+ \int_0^1 \, s\phi(\cdot +s)\, ds\in
AP(\mathbf{R},X)+ t\,AP(\mathbf{R},X)$.

(c)  Let $\phi =  f+ tg $, where $f, g \in AP(\mathbf{R},X)$. Then
$P \phi =P f+ P( tg)$. By  Corollary 7.11 (a), $Pf \in
AP_{w_1}(\mathbf{R},X)$. The assumptions imply $P(tg)\in
BUC_{w_1}(\mathbf{R},X)$ and therefore (b) implies $ P(tg)\in
AP_{w_1}(\mathbf{R},X)$. It follows $P\phi \in
AP_{w_1}(\mathbf{R},X)$.
\end{proof}

 \smallskip

The following shows that the condition $\phi \in P_1 AP= t\cdot
AP(\mathbf{R},X)+  AP(\mathbf{R},X)$ of Corollary 7.12 (c)  can
not be replaced by
 $\phi \in AP_{w_1} (\mathbf{R},X)$.  Also Theorem 6.5
 does not hold  without the conditions  $\gamma^{-1}\phi|\,J \in E_w (J,X)$, even for the case  $\phi \in BUC_{w_1} (\mathbf{R},\mathbf{R})$,
 $sp_{w_1} (\phi)$ is residual and  $sp_{AP_{w_1}} (\phi)\subset sp_{P_1 AP} (\phi)=\{0\}$. It also provides  counter
 examples showing  that   the assumption $sp_{\A} (\phi)=\emptyset$ of Proposition 6.3
 (b).

\smallskip

\textbf{Example 7.13}.  Let $\psi (t)= \sum_{k=1}^{\infty} \,
k^{-4/3} \cos (t\, k^{-1/3})$, $t\in \mathbf{R}$. Then

(i) $\psi \in AP(\mathbf{R},\mathbf{R})$, $P^2\psi \in
BUC_{w_1}(\mathbf{R},\mathbf{R})$ but $P\psi\not  \in
AP(\mathbf{R},\mathbf{R})$ and   $P^2\psi\not  \in
AP_{w_1}(\mathbf{R},\mathbf{R})$.

(ii) $sp (\psi)= (k^{-1/3})_{k=1}^{\infty}\cup \{0\}=sp_{w_1}
(P\psi)=sp_{w_1}(P^2\psi) $.

(iii) $sp_{P_1 AP} (P\psi) =sp_{P_1 AP} (P^2 \psi)=\{0\} $.

\begin{proof} (i) That $\psi \in AP(\mathbf{R},\mathbf{R})$ is clear
since the series $\sum_{k=1}^{\infty} \, k^{-4/3} \cos (t\,
k^{-1/3})$ is absolutely convergent. One has $ P\psi=
\sum_{k=1}^{\infty} \, k^{-1} \sin (t\, k^{-1/3})$ and $P^2 \psi=
\sum_{k=1}^{\infty} 2\, k^{-2/3} \sin^2 (1/2)(t\, k^{-1/3})$.
Using the inequalities $|\sin y| \le 1$ for $y\in \mathbf {R}$, $
(2/\pi)\, |y| \le |\sin y| \le |y|$ for $|y| \le \pi/2$ and simple
calculations, one has

\smallskip

(7.1)\qquad  $ 6t^2 (\pi(|t|^3 +\pi^3)^{-1/3}) \le (2t^2/\pi^2)
\sum_{k \,\ge |t|^3/\pi^3} \, k^{-4/3} \le   P^2 \psi (t) \le 9
(1+|t|)$.

\smallskip

\noindent  This implies $P^2\psi \in
BUC_{w_1}(\mathbf{R},\mathbf{R})$. If $P\psi  \in
AP(\mathbf{R},\mathbf{R})$, then its Fourier series is as above
and so $MP \psi=0$. Therefore by Theorem 5.5 (c), $P^2 \psi\in
C_{w_1,0}(\mathbf{R},\mathbf{R})$. This contradicts (7.1) and
proves $P\psi\not  \in AP(\mathbf{R},\mathbf{R})$.

If $P^2\psi  \in AP_{w_1}(\mathbf{R},\mathbf{R})$, by Proposition
5.5 (e), there is $a \in \mathbf{R}$ such that $(P^2\psi -ta)/w_1
\, \in C_0 (\mathbf{R},\mathbf{R})$. This also  contradicts (7.1)
proving $P^2\psi\not  \in AP_{w_1}(\mathbf{R},\mathbf{R})$.

(ii) This follows by [15, Proposition 1.1 and Remark 1.2].

(iii) As in [15, Proposition 1.1  (a), (b) and Theorem 3.4] using
Proposition 6.3 (b), one gets $sp_{P_1AP} (P\psi)$, $sp_{P_1 AP}
(P^2\psi) \subset \{0\} $ and therefore (iii) follows from (i).
\end{proof}

 \bigskip

 8. \textbf{A generalized evolution
equation}

In this section we consider a general equation of which some
recurrence, evolution and convolution\ equations will be seen to
be special cases in subsequent sections. We will assume that
${{\F}} $ is a translation invariant closed subspace of
$BUC_{w}(J,X)$ satisfying (5.2). We study the equation

\bigskip

(8.1) \qquad $(B\phi )(t)=A\phi (t)+\psi (t)$ for $t\in G$

\bigskip

\noindent under the assumptions

\smallskip

(8.2) $A:X\rightarrow X$ is a closed linear operator with resolvent set $%
\rho (A)$\ non-empty and domain $D(A)$ dense in $X$.

\smallskip

(8.3) $B:D(B)\rightarrow C(G,X)$ is a linear operator with
$TP(G,X)\subseteq D(B)\subseteq BUC_{w}(G,X).$

\smallskip

(8.4) Let $(\phi _{n})$ be any bounded sequence in $D(B)$ such that $%
(B\phi _{n})$ is a bounded sequence in $BC_{w}(G,X).$ If $\phi
_{n}\rightarrow \phi $ and $B\phi _{n}\rightarrow \psi $ uniformly
on compact sets for some $\phi \in BUC_{w}(G,X)$, then $\phi \in
D(B)$ and $B\phi =\psi$.

\smallskip

(8.5) $B(\gamma x)=\theta _{B}(\gamma )\gamma x$ for each $x\in
X,$ $\gamma \in \widehat{G}$ and some continuous function $\theta _{B}:%
\widehat{G}\rightarrow \mathbf{C}.$

\smallskip

(8.6) $B(\phi \ast f)=(B\phi )\ast f$ for all $\phi \in D(B)$ and
all $f\in L_{w}^{1}(G)$ with supp$(\widehat{f})$ compact.

\smallskip

(8.7) If $\xi \in D(B)$ and $sp_{w}(\xi )\subseteq V_{1}$ where $\overline{V}%
_{1}\subseteq V_{2}$ and $V_{1},$ $V_{2}$ are relatively compact open
subsets of $\widehat{G}$, then there exists $\pi _{n}\in TP(G,X)$ such that $%
\pi _{n}(G)\subseteq D(A),$ $sp_{w}(\pi _{n})\subseteq V_{2},$
$\pi _{n}(t)\rightarrow \xi (t),$ $B\pi _{n}(t)\rightarrow B\xi
(t)$ uniformly on compact subsets of $G,$ and sup$_{n}(\left\| \pi
_{n}\right\| _{w,\infty }+\left\| B\pi _{n}\right\| _{w,\infty
})<\infty$.

\smallskip

(8.8) There is a closed subspace $Y$ of $L(X)$ containing all
the resolvent operators of $A$ such that if $\phi \in {{\F}} $ and $%
C\in Y$ then $C\circ \phi \in {{\F} }.$

\smallskip

(8.9) For each $\gamma _{0}\in \theta _{B}^{-1}(\rho (A))$ there
is a compact neighbourhood $V$ of $\gamma _{0}$ and a function
$h\in
L_{w}^{1}(G,Y)$ such that $V\subseteq $ $\theta _{B}^{-1}(\rho (A))$ and $%
\widehat{h}(\gamma )=(\theta _{B}(\gamma )-A)^{-1}$ for $\gamma \in V.$

\smallskip

(8.10) If $\phi \in D(B),$ $sp_{w}(\phi )$ is compact and $\phi
|_{J}\in {{\F} }$ then $B\phi |_{J}\in {{\F}}$.

\bigskip

We will refer to $\theta _{B}$ as the \textit{characteristic function} of $%
B. $ In most of our applications we will have $Y=L(X).$ However,
for Corollary 12.4 it is necessary to have more general $Y.$

\bigskip

\textbf{Lemma 8.1}

(a) If $A$ satisfies (8.2), then $A\circ (\phi \ast f)=(A\circ
\phi )\ast f$ for all $\phi \in BUC_{w}(G,X)$ with $\phi
(G)\subseteq D(A)$ and all $f\in L_{w}^{1}(G)$ with
supp$(\widehat{f})$ compact.

(b) If $B$ satisfies (8.3) and (8.5), then $B(R\circ \pi )=R\circ
(B\pi )$ for all $R\in L(X)$ and $\pi \in TP(G,X).$

(c) If $B$ satisfies (8.3), (8.4), (8.5) and (8.7), then $B(R\circ
\xi )=R\circ (B\xi )$ for all $R\in L(X)$ and $\xi \in D(B)$ with
$sp_{w}(\xi )$ compact.

(d) If ${{\F} }$ satisfies (8.8), $\phi \in BUC_{w}(G,X),$ $\phi
|_{J}\in {\F} $ and $h\in L_{w}^{1}(G,Y),$ then $(h\ast \phi
)|_{J}\in {\F}.$

(e) If $\phi \in BC_{w}(G,X),$ $\phi |_{J}\in {{\F} }$ and $C\in
Y$,  then $sp_{w}(C\circ \phi )\subseteq sp_{w}(\phi ).$

\,

\begin{proof} (a) Let $\lambda \in \rho (A)$ and set $R=(\lambda
-A)^{-1}.$ If $\phi \in BUC_{w}(G,X)$ with $\phi (G)\subseteq
D(A)$, $t\in G$ and $f\in L_{w}^{1}(G),$ then

$R\circ A\circ (\phi \ast f)(t)=R\circ A\int_{G}\phi (t-s)f(s)d\mu
(s)=$ $\int_{G}R\circ A\phi (t-s)f(s)d\mu (s)$

\ \ \ \ \ \ \ \ \ \ \ \ \ \ \ \ \ \ \ \ \ \ $=R\int_{G}A\phi
(t-s)f(s)d\mu (s)=R\circ ((A\circ \phi )\ast f)(t)$

\noindent from which the result follows.

(b) If $\pi =\gamma x,$ where $\gamma \in \widehat{G}$ and $x\in X,$ then $%
B(R\circ \pi )=B(\gamma Rx)=\theta _{B}(\gamma )Rx=R(\theta _{B}(\gamma
)x)=R\circ (B\pi ).$ The result follows by linearity.

(c) Choose $V_{1},$ $V_{2}$ and $\pi _{n}$ as in (8.7). Then $%
R\circ \pi _{n}\rightarrow R\circ \xi $ and, by (b), $B(R\circ \pi
_{n})=R\circ (B\pi _{n})\rightarrow R\circ (B\xi )$ uniformly on
compact sets. By (8.4), $R\circ \xi \in D(B)$ and $B(R\circ \xi
)=R\circ (B\xi ).$

(d) By a theorem of Bochner ([57, p.133] ) there is a sequence of
step functions $(h_{n})$ converging to $h$ in $L_{w}^{1}(G,Y).$ Hence, $%
h_{n}\ast \phi \rightarrow $ $h\ast \phi $ in $BUC_{w}(G,Y).$ By
(5.2), $\phi _{-s}|_{J}\in {{\F}} $ for each $s\in G$ and so by (8.8), $%
(h_{n}\ast \phi )|_{J}\in {\F} $. Hence, $(h\ast \phi )|_{J}\in
{\F}$.

(e) For each $f\in L_{w}^{1}(G)$ we have $(C\circ \phi )\ast
f=C\circ (\phi \ast f).$ Hence, $I_{w}(C\circ \phi )\supseteq
I_{w}(\phi )$ and so $sp_{w}(C\circ \phi )\subseteq sp_{w}(\phi
).$ \end{proof}

\smallskip

\textbf{Lemma 8.2}. If $k\in L_{w}^{1}(G,L(X))$ and $(\eta _{n})$
is a bounded sequence in $BUC_{w}(G,X)$ converging to $\eta $
uniformly on compact subsets of $G$, then $k\ast \eta
_{n}\rightarrow k\ast \eta $ uniformly on compact subsets.

\,

\begin{proof} Let $U$ be a symmetric compact neighbourhood of $0$ in $G$.
Given $\varepsilon >0$ choose a compact symmetric set $K_{0}\subseteq G$
such that $\int_{G\backslash K_{0}}\left\| k(s)\right\| w(s)d\mu (s)<\frac{%
\varepsilon }{2c_{1}}$ where $c_{1}=(1+\sup_{n}(\left\| \eta _{n}\right\|
_{w,\infty }+\left\| \eta \right\| _{w,\infty }))\sup_{s\in U}w(s).$ Let $%
K=U+K_{0}$ and choose $n_{\varepsilon }$ such that $\left\| \eta
_{n}(s)-\eta (s)\right\| <\frac{\varepsilon }{2c_{2}}$ for all $%
n>n_{\varepsilon }$ and all $s\in K$ where $c_{2}=1+\left\| k\right\|
_{w,1}. $ For $t\in U$ and $n>n_{\varepsilon }$ we have

$\left\| k\ast \eta _{n}(t)-k\ast \eta (t)\right\| \leq \int_{G}\left\|
k(t-s)\right\| \left\| \eta _{n}(s)-\eta (s)\right\| d\mu (s)$

 $\ \ \ \ \ \ \ \ \ \ \ \ \ \ \ \ \ \ \ \ \ \ \ \ \ \ \ \ \ \ \
\ \leq \int_{G\backslash K}\left\| k(t-s)\right\| w(s)c_{1}d\mu
(s)+\int_{K}\left\| k(t-s)\right\| \frac{\varepsilon }{2c_{2}}w(t-s)d\mu (s)$

\ \ \ \ \ \ \ \ \ \ \ \ \ \ \ \ \ \ \ \ \ \ \ \ \ \ \ \ \ \ \ \ \ $\leq
\int_{G\backslash (t+K)}\left\| k(s)\right\| w(t-s)c_{1}d\mu (s)+\left\|
k\right\| _{w,1}\frac{\varepsilon }{2c_{2}}$

\ \ \ \ \ \ \ \ \ \ \ \ \ \ \ \ \ \ \ \ \ \ \ \ \ \ \ \ \ \ \ \ \ $\leq
\int_{G\backslash K_{0}}\left\| k(s)\right\| w(t-s)c_{1}d\mu (s)+\frac{%
\varepsilon }{2}$

\ \ \ \ \ \ \ \ \ \ \ \ \ \ \ \ \ \ \ \ \ \ \ \ \ \ \ \ \ \ \ \ \ $%
<\varepsilon$.

\noindent Hence $k\ast \eta _{n}\rightarrow k\ast \eta $ uniformly
on $U$ and therefore on each compact subset of $G.$ \end{proof}

\smallskip

\textbf{Theorem 8.3}. Suppose $\phi ,\psi \in BUC_{w}(G,X),$
conditions (8.2)-(8.10) hold and $B\phi =A\circ \phi +\psi$. If
$\psi |_{J}\in {{\F} }$ then $sp_{{{\F} }}(\phi )\subseteq \theta
_{B}^{-1}(\sigma (A)).$

\,

\begin{proof} Let $\lambda _{0}\in \rho (A)\cap \theta _{B}(\widehat{G})$
and choose $\gamma _{0}\in \widehat{G}$ with $\theta _{B}(\gamma
_{0})=\lambda _{0}.$ Write $R=(\theta _{B}(\gamma _{0})-A)^{-1},$ let $%
B_{\delta }(\lambda _{0})=\{\lambda \in \mathbf{C}:\left| \lambda
-\lambda _{0}\right| <\delta \}\ $and choose $\delta >0$ such that
$B_{\delta }(\lambda _{0})\cap \sigma (A)=\emptyset$. By (8.9)
there is a relatively compact open neighborhood $V$ of $\gamma
_{0}$ and a function $h\in
L_{w}^{1}(G,Y)$ such that $V\subseteq $ $\theta _{B}^{-1}(\rho (A))$ and $%
\widehat{h}(\gamma )=(\theta _{B}(\gamma )-A)^{-1}$ for $\gamma \in V.$ Set $%
V_{j}=V\cap \theta _{B}^{-1}(B_{j\delta /5}(\lambda _{0}))$ for $j=1,2,3,4.$
Choose $f\in L_{w}^{1}(G)$ with $\widehat{f}(\gamma _{0})=1$ and supp$(%
\widehat{f})\subseteq V_{1}$ and set $\xi =\phi \ast f.$ It
suffices to prove $\xi |_{J}\in {{\F}} ,$ for then $f\in I_{{{\F}}
}(\phi )$ and so $\gamma _{0}\notin sp_{{{\F}} }(\phi ).$

Note that, by Lemma 8.1(a) and (8.6), we have $B\xi -A\circ \xi
=\psi \ast
f. $ Choose $g\in L_{w}^{1}(G)$ with $\widehat{g}=1$ on $V_{3}$ and supp$(%
\widehat{g})\subseteq V_{4}.$ Setting $k=g\ast h$ gives$\
\widehat{k}(\gamma )=\widehat{g}(\gamma )(\theta _{B}(\gamma
)-A)^{-1}.$ Next, $sp_{w}(\xi )\subseteq $
supp$(\widehat{f})\subseteq V_{1}$ and so we can choose $\pi _{n}$
as in (8.8). Set $\eta _{n}=B(R\circ \pi _{n})-A\circ (R\circ \pi
_{n}) $ and note that $(\eta _{n})$ is a bounded sequence in $%
BUC_{w}(G,X).$ Now for $\gamma \in V_{2}$ and $x\in D(A)$ we have $k\ast
(B(\gamma x)-A\circ (\gamma x))=(k\ast \gamma )\circ (\theta _{B}(\gamma
)-A)x=\widehat{k}(\gamma )\circ (\theta _{B}(\gamma )-A)\gamma x=\gamma x.$
Hence, $k\ast \eta _{n}=R\circ \pi _{n}.$

Since $R$ and $A\circ R$ are continuous operators and by Lemma
8.1(b), $\eta _{n}\rightarrow R\circ (B\xi -A\circ \xi )=R\circ
(\psi \ast f)=(R\circ \psi )\ast f$ and $k\ast \eta
_{n}\rightarrow R\circ \xi $ uniformly on compact subsets of $G.$
By Lemma 7.2, $k\ast \eta _{n}\rightarrow k\ast (R\circ (\psi \ast
f))$ pointwise on $G$ and so $R\circ \xi =k\ast (R\circ \psi )\ast
f$ . By Proposition 8.3(v),(ii) $(\psi \ast f)|_{J}\in {{\F} }$
and
so, by (8.8), $((R\circ \psi )\ast f)|_{J}\in {{\F} }.$ By Lemma 8.1(d), $%
(k\ast (R\circ \psi )\ast f)|_{J}\in {{\F}} ,$ that is $(R\circ
\xi )|_{J}\in {\F}$. By Lemma 8.1(e), $sp_{w}(R\circ \xi )$ is
compact. By Lemma 8.1(c) and (8.10), $(R\circ (B\xi
))|_{J}=B(R\circ \xi )|_{J}\in {{\F} }$ , so $(A\circ R\circ \xi
)|_{J}=B(R\circ \xi )|_{J}-(R\circ \psi \ast f)|_{J}\in {{\F}} .\
$Hence, $\xi |_{J}=\theta _{B}(\gamma _{0})(R\circ \xi
)|_{J}-(A\circ R\circ \xi )|_{J}\in {{\F}} ,$ as required.
\end{proof}

\smallskip

\textbf{Example 8.4}. Consider the equation $B\phi =A\circ
\phi +\psi ,$ where $\phi :\mathbf{R}\rightarrow \mathbf{C}^{2},$ $%
B\phi =\phi ^{\prime \prime }+2i\phi ^{\prime }-\phi $ and $A$ has
matrix representation $\left[
\begin{array}{ll}
2 & -6 \\
3 & -7%
\end{array}
\right]$. Let ${{\F}} \in \{AP_{w_{N}}(\mathbf{R},\mathbf{C}^{2}),$ $%
C_{w_{N},0}(\mathbf{R},\mathbf{C}^{2}),$ $0_{\mathbf{R}}\}$ and $%
w_{N}(t)=(1+\left| t\right| )^{N}.$This is a special case of the
equations studied in section 9, where it is shown that (8.1)-(8.9)
hold. The general solution of the homogeneous equation $B\xi
=A\circ \xi $\ is $\xi =(2c_{1}+c_{2}\gamma _{1}+2c_{3}\gamma
_{-2}+c_{4}\gamma
_{-3},c_{1}+c_{2}\gamma _{1}+c_{3}\gamma _{-2}+c_{4}\gamma _{-3})$ where $%
\gamma _{s}(t)=e^{ist}.$ Moreover, $\theta _{B}^{-1}(\sigma
(A))=\{\gamma _{0},\gamma _{1},\gamma _{-2},\gamma _{-3}\}$ and so
$sp_{w}(\xi )\subseteq \theta _{B}^{-1}(\sigma (A)).$ That this
holds more generally is shown in Corollary 8.5. Now let $\phi $ be
any solution of $B\phi =A\circ \phi +\psi ,$ with say $\psi
(t)=(t,0).$ Then  $sp_{w}(\psi )\cap \theta _{B}^{-1}(\sigma
(A))=\{\gamma_0\}$ and $\phi$ contains a quadratic term. Hence
$\phi \not\in C_{w_{2},0}(\mathbf{R},\mathbf{C}^{2})\cup
AP_{w_{1}}(\mathbf{R},\mathbf{C}^{2})$ whereas $\xi ,\psi \in C_{w_{2},0}(%
\mathbf{R},\mathbf{C}^{2})\cap
AP_{w_{1}}(\mathbf{R},\mathbf{C}^{2})$.

\bigskip

Motivated by examples such as this last one, we refer to $\theta
_{B}^{-1}(\sigma (A))$\ as the \textit{resonance set} of (8.1) and
consider classes ${{\F} }$ satisfying

\bigskip

(8.11) $\xi |_{J}\in {{\F}} $ for all solutions $\xi \in
BUC_{w}(G,X)$ of the homogeneous equation $B\xi =A\circ \xi$.

\bigskip

We say that $\phi \in BUC_{w}(G,X)$ is a \textit{resonance
solution} of (8.1) for the class ${{\F}} $ if (8.11) holds, $\psi
|_{J}\in {{\F}} $ and $\phi $\ satisfies (8.1), \ but $\phi
|_{J}\notin
{\F}$. Corollary 8.5(c) shows that if conditions (8.2)-(8.10) hold and $%
sp_{w}(\psi )\cap \theta _{B}^{-1}(\sigma (A))=\emptyset $, then
there are no resonance solutions for any class ${{\F} }.$

\bigskip

\textbf{Corollary 8.5}. Suppose conditions (8.2)-(8.10) hold,$\
\phi ,\psi \in BUC_{w}(G,X)\ $and $B\phi =A\circ \phi +\psi$.

(a) If $\psi =0,$ then $sp_{w}(\phi )\subseteq \theta
_{B}^{-1}(\sigma (A)).$

(b) If (8.11) holds, then $sp_{{{\F}} }(\phi )\subseteq
sp_{w}(\psi ).$

(c) If (8.11) holds, $sp_{w}(\psi )\cap \theta _{B}^{-1}(\sigma
(A))=\emptyset $ and $\psi |_{J}\in {{\F}} ,$\ then $\phi |_{J}\in
{{\F}}$.

(d) If $\theta _{B}^{-1}(\sigma (A))$ is residual, ${{\F}} $ is a
$\Lambda _{w}^{0}$-class, $\gamma ^{-1}\phi |_{J}\in
E_{w}^{0}(J,X)$ for all $\gamma \in sp_{{{\F}} }(\phi )$ and $\psi
|_{J}\in {{\F}} ,$ then $\phi |_{J}\in {{\F}}$.

(e) If $\theta _{B}^{-1}(\sigma (A))$ is residual, ${{\F}} $ is a
$\Lambda _{w}$-class, $\gamma ^{-1}\phi |_{J}\in E_{w}(J,X)$ for
all $\gamma \in sp_{{{\F}} }(\phi )$ and $\psi |_{J}\in {{\F}} ,$
then $\phi |_{J}\in {{\F}}$.

(f) If $\theta _{B}^{-1}(\sigma (A))\cap sp_{w}(\phi )^{\prime }$
is empty, ${{\F}} $ is a $\Lambda _{w}$-class and $\psi |_{J}\in
{{\F}} ,$ then $\phi |_{J}\in {{\F}}$.

\,

\begin{proof} (a) Conditions (8.8) and (8.10) are trivial when $J=G$ and $%
{{\F}} =\{0\}.$ So we can apply Theorem 8.3 to get $sp_{w}(\phi
)=sp_{\{0\}}(\phi )\subseteq \theta _{B}^{-1}(\sigma (A)).$

(b) For $\gamma \in \widehat{G}$ $\backslash $ $sp_{w}(\psi )$ there exists $%
f\in L_{w}^{1}(G)$ such that $\widehat{f}(\gamma )\neq 0,$
$\widehat{f}$\ has compact support and $\psi \ast f=0.$ By Lemma
8.1(a) and (8.6), $\xi =\phi \ast f$ satisfies $B\xi =A\circ \xi
.$ Hence $\xi |_{J}\in {{\F}} $ and so $\gamma \notin $
$sp_{{{\F}} }(\phi ).$ Thus $sp_{{{\F}} }(\phi )\subseteq
sp_{w}(\psi ).$

(c) By Theorem 8.3 and (b), $sp_{{{\F}} }(\phi ))\subseteq \theta
_{B}^{-1}(\sigma (A))\cap sp_{w}(\psi )=\emptyset$. By Proposition
6.3 (a)(v), $\phi |_{J}\in {{\F}}$.

(d),(e) By Theorem 8.3, $sp_{{{\F}} }(\phi )\subseteq \theta
_{B}^{-1}(\sigma (A))$ which is residual. By Theorem 6.5, $\phi
|_{J}\in {{\F}}$.

(f) By Theorem 8.3 and Corollary 6.7, $sp_{{{\F}} }(\phi
)\subseteq \theta _{B}^{-1}(\sigma (A))\cap sp_{w}(\phi )^{\prime
}$ which is empty. By Proposition 6.3 (a)(v), $\phi |_{J}\in
{{\F}}$.
\end{proof}

\smallskip

We record here a useful theorem for studying almost periodic
functions. Cases (a),(c) for $G=\mathbf{R}$ appear in [43, p.92,
Theorem 4], [54, Theorem 5.2]\ and [5, Theorem 4.3]. Case (b) for
$G=\mathbf{R}$ is in [9, Theorem 4.2.6].

\bigskip

\textbf{Proposition 8.6}. Suppose $\phi \in BUC(G,X),$ ${{\F}}
=AP(G,X)$ and $sp_{{{\F}} }(\phi )$ is residual. Then $\phi \in
AP(G,X)$ provided one of the following conditions is satisfied.

\qquad (a) $c_0\not\subset X;$

\qquad (b) $\gamma ^{-1}\phi \in E(G,X)$ for all $\gamma \in
\widehat{G};$

\qquad (c) $\phi (G)$ is relatively weakly compact.

\bigskip

\textbf{Corollary 8.7}. Suppose $w$\ has polynomial growth and
$G$\ is compactly generated. Suppose conditions (8.2)-(8.10) hold
with ${{\F}}
=AP_{w}(G,X).$ Let $\phi \in BUC_{w}(G,X),$ $\psi \in AP_{w}(G,X)$ and $%
B\phi =A\circ \phi +\psi$. If $\theta _{B}^{-1}(\sigma (A))$ is
residual, then $\phi \in AP_{w}(G,X)$ provided one of the
following conditions is satisfied.

\qquad (a) $w$ is bounded and $c_0\not\subset X;$

\qquad (b) $\gamma ^{-1}\phi \in E_{w}(G,X)$ for all $\gamma \in
\theta _{B}^{-1}(\sigma (A)).$

\qquad (c) $w$ is bounded and $\phi (G)$ is relatively weakly
compact.

\,

\begin{proof} By Theorem 8.3, $sp_{{{\F}} }(\phi )\subseteq
\theta _{B}^{-1}(\sigma (A))$ which is residual. Cases (a),(c)
follow from Proposition 8.6, case (b) from Theorem 6.5.
\end{proof}
\smallskip

The following lemmas will be used in subsequent sections to verify
conditions (8.7),(8.8).

\bigskip

\textbf{Lemma 8.8}. Let $G$\ be compactly generated.

(a) Suppose $\xi \in BUC_{w}(G,X)$ and $sp_{w}(\xi )\subseteq V_{1}$ where $%
\overline{V}_{1}\subseteq V_{2}$ and $V_{1},$ $V_{2}$ are relatively compact
open subsets of $\widehat{G}$. Let $D$\ be a dense linear subspace of $X.$\
Then there exists $\pi _{n}\in TP(G,X)$ with $\pi _{n}(G)\subseteq D,\
sp_{w}(\pi _{n})\subseteq V_{2},$ $\sup_{n}\left\| \pi _{n}\right\|
_{w,\infty }<\infty $ and $\pi _{n}\rightarrow \xi $ uniformly on compact
subsets of $G.$

(b) If in addition $G=\mathbf{R}$\ and $\xi \in
C^{m}(\mathbf{R},X),$ then $\pi _{n}$ can be chosen with

\noindent $\sup_{n}\left\| \pi _{n}^{(j)}\right\| _{w,\infty
}<\infty $ and $\pi _{n}^{(j)}\rightarrow \xi ^{(j)}$ uniformly on
compact sets for $0\leq j\leq m.$

\,

\begin{proof} (a) Since $G$\ is compactly generated we have
$G=\cup _{n=1}^{\infty }U_{n}$ where $\overline{U}_{n}\subset
U_{n+1}$ for some relatively compact, open, generating subsets
$U_{n}.$ By Proposition 2.3 we
may assume $\sup_{t\in U_{n}}w(t)\leq c\inf_{t\notin U_{n}}w(t).$ Choose $%
\alpha _{n}\in C(G)$ with $\alpha _{n}=1$ on $U_{n},$ $\alpha _{n}=0$ on $%
G\backslash U_{n+1}$ and $0\leq \alpha _{n}\leq 1.$ Set $\xi _{n}=\alpha
_{n}\xi$. By Lemma 5.3(a), there exists $\psi _{n}\in TP(G,X)$ with $%
\sup_{t\in U_{n+1}}\left\| \psi _{n}(t)-\xi _{n}(t)\right\|
<\frac{1}{n}$ and $\sup_{t\in G}\left\| \psi _{n}(t)\right\|
<\frac{1}{n}+\sup_{t\in U_{n+1}}\left\| \xi _{n}(t)\right\|$.
Moreover, since $D$ is dense in $X$ we may assume that the
coefficients of $\psi _{n}$ are in $D.$\ Take $f\in
L_{w}^{1}(G)$ with $\widehat{f}=1$ on $V_{1}$ and $\widehat{f}=0$ on $%
G\backslash V_{2}.$\ Set $\pi _{n}=f\ast \psi _{n}.$ \ Now $\psi
_{n}\rightarrow \xi $ uniformly on compact sets and $(\psi _{n})$ is bounded
in $BUC_{w}(G,X)$ since $\sup_{t\in G}\frac{\left\| \psi _{n}(t)\right\| }{%
w(t)}<\sup_{t\in U_{n+1}}\frac{\frac{1}{n}+\left\| \xi (t)\right\| }{w(t)}%
<\infty$. By Lemma 8.2 we are finished.

(b) This follows in the same way using Lemma 5.3(b).\end{proof}

\smallskip

\textbf{Lemma 8.9}. If ${{\F}} \in \{AP_{w}(G,X)$, $C_{w,0}(J,X),$ $%
0_{G}\} $ then (8.8) holds with $Y=L(X).$
 \,

\begin{proof} If $C\in L(X)$ and $\pi \in TP_{w}(G,X)$ then
$C\circ \pi \in TP_{w}(G,X).$ By continuity, (8.8) holds for
${{\F}} =AP_{w}(G,X).$ The result for $C_{w,0}(J,X)$ and $0_{G}$
is obvious.
\end{proof}
\bigskip

\textbf{9. Differential equations}

\smallskip

Throughout this section we assume that $A$ satisfies (8.2) and $J\in \{%
\mathbf{R,R}_{+}\}.\ $Also, the operator $B$\ is given by$\
B=\sum_{j=0}^{m}b_{j}\left( \frac{d}{dt}\right) ^{j}$ where $b_{j}\in
\mathbf{C}.$ Assume $b_{m}\neq 0$ so that $B$ has order $m.$ Given \vspace{%
1pt}$\psi \in $ $BUC_{w}(\mathbf{R},X)$ we investigate the
behaviour of solutions $\phi :J\rightarrow D(A)$ of the evolution
equation

\bigskip

(9.1) $(B\phi )(t)=A\phi (t)+\psi (t)$ for $t\in \mathbf{R}.$

\bigskip

We seek conditions on $X,$ $A$ and translation invariant closed subspaces $%
{{\F}} $ of $BUC_{w}(J,X)$ which ensure that if $\psi |_{J}\in
{{\F}} $
then $\phi |_{J}\in {{\F}}$. We apply Theorem 8.3 with $D(B)=BUC_{w}(%
\mathbf{R},X)\cap C^{m}(\mathbf{R,}X),$ assuming throughout that
${{\F}} $ satisfies (5.2).

\bigskip

Let $\alpha :\widehat{\mathbf{R}}\rightarrow \mathbf{R}$ be the natural
isomorphism given by $\alpha (\gamma _{s})=s$ where $\gamma _{s}(t)=e^{ist}$
for $s,t\in \mathbf{R}.$ As usual, the space of $m$-times continuously
differentiable functions $u:J\mathbf{\rightarrow }X$ is $C^{m}(J\mathbf{,}X)$
and the subspace of functions with compact support is $C_{c}^{m}(J\mathbf{,}%
X).$

\bigskip

To find the characteristic function of $B$ let $p_{B}(s)=%
\sum_{j=0}^{m}b_{j}(is)^{j}$ where $s\in \mathbf{R}.$ For $\gamma
\in \widehat{\mathbf{R}}$ and $x\in X$ we have $B(\gamma
x)=p_{B}(\alpha (\gamma ))\gamma x$ and so $\theta _{B}=p_{B}\circ
\alpha$.

\bigskip

\textbf{Lemma 9.1}. If $w$ has polynomial growth of order $N\ $\
and $u\in C_{c}^{N+2}(\mathbf{R,}X),$ then there exists $h\in L_{w}^{1}(%
\mathbf{R,}X)$ such that $\widehat{h}=u\circ \alpha$.

\,

\begin{proof} Let $h(t)=\frac{1}{2\pi }\int_{-\infty }^{\infty
}u(\tau )e^{it\tau }d\tau$. Integration by parts gives

\,

$(-it)^{N+2}h(t)=%
\frac{1}{2\pi }\int_{-\infty }^{\infty }u^{(N+2)}(\tau )e^{it\tau
}d\tau $ and so $h\in L_{w}^{1}(\mathbf{R,}X).$

\smallskip

\noindent By the Fourier inversion Theorem, for
an appropriate normalization of Haar measure on $\mathbf{R}$ we have $%
x^{\ast }\circ \widehat{h}=\widehat{x^{\ast }\circ h}=x^{\ast
}\circ u\circ \alpha $ for each $x^{\ast }\in X^{\ast }.$ Hence
$\widehat{h}=u\circ \alpha$. \end{proof}

\smallskip

\textbf{Lemma 9.2. }Let $C(J,X)$ have the topology of uniform
convergence on compact sets. Set $P^{0}\phi =\phi $ and $P^{j}\phi
=P(P^{j-1}\phi )$ for $j\geq 1.$

(a) The operator $P:C(J,X)\rightarrow C(J,X)$ is continuous.

(b) For $\phi \in C(J,X)$ and $0\leq j\leq m$ we have

$P^{k}B\phi (t)=\sum_{j=0}^{k}b_{j}P^{k-j}\phi
(t)+\sum_{j=k+1}^{m}b_{j}\phi
^{(j-k)}(t)-\sum_{j=0}^{k-1}\sum_{i=k-j}^{m}b_{i}\phi ^{(i+j-k)}(0)\frac{%
t^{j}}{j!}.$

(c) The operator $B:C(J,X)\rightarrow C(J,X),$ with domain $D(B)=C^{m}(J,X),$
is closed.

\,

\begin{proof} (a) Suppose $\phi _{n}\in C(J,X)$ and $\phi
_{n}\rightarrow \phi $ uniformly on compact sets. Let $K$ be a
compact subset of $J$ and let $\widetilde{K}$ be the convex hull
of $K\cup \{0\}.$ Set $c=\max \{\left| t\right| :t\in
\widetilde{K}\}.$ Given $\varepsilon >0$ there exists
$n_{\varepsilon }$ such that $\left\| \phi _{n}(s)-\phi
(s)\right\| <\frac{\varepsilon }{1+c}$ for all $s\in \widetilde{K}$ and $%
n>n_{\varepsilon }.$ Hence $\left\| P\phi _{n}(t)-P\phi
(t)\right\| =\left\| \int_{0}^{t}(\phi _{n}(s)-\phi (s))ds\right\|
<\varepsilon $ for all $t\in K$ and $n>n_{\varepsilon }.$

(b) A simple induction argument gives this result.

(c) Suppose $\phi _{n}\in D(B),$ $\phi _{n}\rightarrow \phi $ and $%
B\phi _{n}\rightarrow \psi $ uniformly on compact sets. By (b) we have $%
P^{m}B\phi _{n}=\sum_{j=0}^{m}b_{j}P^{m-j}\phi _{n}-p_{n}$ where $%
p_{n}(t)=\sum_{j=0}^{m-1}\sum_{i=m-j}^{m}b_{i}\phi _{n}^{(i+j-m)}(0)\frac{%
t^{j}}{j!}$ is a polynomial of degree at most $m-1.$ By (a) we find $%
p_{n}\rightarrow \sum_{j=0}^{m}b_{j}P^{m-j}\phi -P^{m}\psi $
uniformly on compact sets. Hence,
$p=\sum_{j=0}^{m}b_{j}P^{m-j}\phi -P^{m}\psi \in P^{m-1}(J,X).$ So
$b_{m}\phi =p+P^{m}\psi -\sum_{j=0}^{m-1}b_{j}P^{m-j}\phi \in
C^{1}(J,X).$ \ Since $b_{m}\neq 0$ we have $\phi \in C^{1}(J,X).$
Taking derivatives we get $b_{m}\phi ^{\prime }+b_{m-1}\phi
=p^{\prime }+P^{m-1}\psi -\sum_{j=0}^{m-2}b_{j}P^{m-j-1}\phi \in
C^{1}(J,X).$ Proceeding inductively we find $\phi \in D(B)$ and
$B\phi =p^{(m)}+\psi =\psi$. This proves that $B$ is closed.
\end{proof}
\smallskip

\textbf{Lemma 9.3}. Conditions (8.2)-(8.7) all hold. If $%
{{\F}} \in \{AP_{w}(\mathbf{R},X)$, $C_{w,0}(J,X),$
$0_{\mathbf{R}}\}$ and $Y=L(X)$ then (8.8) holds. If $w$ has
polynomial growth, conditions (8.9),(8.10) hold.

\,

\begin{proof} We have assumed (8.2). Now (8.3),(8.5), (8.6) are
clear, (8.4) follows readily from Lemma 9.2(c), and (8.7) follows
from Lemma 8.8. If ${{\F}} \in \{AP_{w}(\mathbf{R},X)$,
$C_{w,0}(J,X),$ $0_{\mathbf{R}}\}$ and $Y=L(X),$ then (8.8) holds
by Lemma 8.9. Assume $w$ has polynomial
growth. For (8.9) let $V$ be a relatively compact open subset of $\widehat{%
\mathbf{R}}$ with $\theta _{B}(\overline{V})\cap \sigma (A)=\emptyset$. So $%
\overline{V}\subseteq W$ where $W$ is also a relatively compact open subset
of $\widehat{\mathbf{R}}$ with $\theta _{B}(\overline{W})\cap \sigma
(A)=\emptyset$. Choose $v\in C_{c}^{N+2}(\mathbf{C,R})$ such that $v=1$ on $%
\theta _{B}\left( V\right) $ and supp$(v)\subseteq \theta _{B}(W).$ Setting $%
r(\zeta )=v(\zeta )(\zeta -A)^{-1}$ and $u=r\circ p_{B}$\ we have $u\in
C_{c}^{N+2}(\mathbf{R,}Y).$ By Lemma 9.1, there exists $h\in L_{w}^{1}(%
\mathbf{R,}Y)$ such that $\widehat{h}=u\circ \alpha $ and so $\widehat{h}%
(\gamma )=(\theta _{B}(\gamma )-A)^{-1}$ for $\gamma \in V,$
proving (8.9).
To prove (8.10) let $\phi \in D(B)$ with $sp_{w}(\phi )$ compact and $%
\phi |_{J}\in {{\F}}$. By Proposition 8.1(a),(b), $\phi ^{(j)}\in
BUC_{w}(\mathbf{R},X)$ for all $j\geq 0$ and so $\phi
^{(j)}|_{J}\in {{\F}}$. Hence $B\phi |_{J}\in {{\F}} $ as
required. \end{proof}

\smallskip

\textbf{Theorem 9.4}. Suppose $w$ has polynomial growth, $\phi
,\psi \in BUC_{w}(\mathbf{R},X),$ condition (8.8) holds and $B\phi
=A\circ \phi +\psi$. If $\psi |_{J}\in {{\F}} $ , then $\theta
_{B}(sp_{{{\F}} }(\phi ))\subseteq \sigma (A)\cap
p_{B}(\mathbf{R}).$

\,

\begin{proof} By Lemma 9.3, the result follows from Theorem 8.3. \end{proof}

\smallskip

\textbf{Corollary 9.5}. Assume $w$\ has polynomial growth,\ $\phi
\in BUC_{w}(\mathbf{R},X),$ $\psi \in AP_{w}(\mathbf{R},X)$ and
$B\phi =A\circ \phi +\psi$. Let $B$ have positive order and
$\sigma (A)\cap p_{B}(\mathbf{R})$ be residual. Then $\phi \in
AP_{w}(\mathbf{R},X)$ provided one of the following conditions is
satisfied.

\qquad (a) $w$ is bounded and $c_0\not\subset X;$

\qquad (b) $\gamma ^{-1}\phi \in E_{w}(\mathbf{R},X)$ for all
$\gamma \in \theta _{B}^{-1}(\sigma (A));$

\qquad (c) $w$ is bounded and $\phi (\mathbf{R})$ is \allowbreak
relatively weakly compact.

\,

\begin{proof}  Let ${{\F}} =AP_{w}(\mathbf{R},X)$ and $%
Y=L(X).$ Since $\sigma (A)\cap p_{B}(\mathbf{R})$ is residual and
$p_{B}$ is a non-constant polynomial, $\theta _{B}^{-1}(\sigma
(A))$\ is residual. By Lemma 9.3, the result follows from
Corollary 8.7. \end{proof}

\smallskip

\textbf{Example 9.7}.

(a) If $B\phi =\phi ^{\prime }$ then
$p_{B}(\mathbf{R})=i\mathbf{R}$ and Corollary 9.5, with $w=1,$
appears in [6] (see also [11, Theorem 3.3]).

(b) If $B\phi =\phi ^{\prime \prime }$ then $p_{B}(\mathbf{R})=\mathbf{%
(-\infty ,}0\mathbf{]}$ and Corollary 9.5, $w=1,$\ appears in [5].

(c) If $B\phi =\phi ^{\prime }+\phi ^{\prime \prime }$ then $p_{B}(%
\mathbf{R})$ is the parabola $Re(z)+Im(z)^{2}=0.$

\bigskip

The following result was proved by Arendt and Batty [5] for the
cases $B\phi =\phi ^{\prime },\phi ^{\prime \prime
}$ under the additional assumptions that $w=1$\ and$\ \sigma (A)\cap p_{B}(%
\mathbf{R})$ consists only of poles of the resolvent operator of $A.$

\bigskip

\textbf{Corollary 9.7}. Assume $w$\ has polynomial growth.\ Let
$\sigma (A)\cap p_{B}(\mathbf{R})$ be discrete and assume $p_{B}$
is not constant. If $\phi \in BUC_{w}(\mathbf{R},X)$ satisfies
$B\phi =A\circ \phi $ then $\phi \in AP_{w}(\mathbf{R},X).$

\,

\begin{proof} Apply Theorem 9.4 with ${{\F}} =0_{\mathbf{R}}$ to get $%
sp_{w}(\phi )=sp_{{{\F}} }(\phi )\subseteq \theta _{B}^{-1}(\sigma
(A)\cap p_{B}(\mathbf{R})).$ So $sp_{w}(\phi )$ is discrete and by
Corollary 8.7 with ${{\F}} =AP(\mathbf{R},X)$ we have $sp_{{{\F}}
}(\phi )\subseteq sp(\phi )^{\prime }=\emptyset$. By Proposition
6.3 (a)(v), $\phi \in {{\F}}.$ \end{proof}

\bigskip

\textbf{10. Convolution equations}

Throughout this section we assume $A$ satisfies (8.2), $%
G\in \{\mathbf{R}^{d},\mathbf{Z}^{d}\}$ and the operator $%
B:BUC_{w}(G,X)\rightarrow BUC_{w}(G,X)$\ is given by$\ B\phi
=k\ast
\phi $ where $k\in L_{w}^{1}(G).$ Given $\psi \in $ $%
BUC_{w}(G,X)$ we investigate the behaviour of solutions $\phi
:G\rightarrow D(A)$ of the convolution equation

\bigskip

(10.1) $(k\ast \phi )(t)=A\phi (t)+\psi (t)$ for $t\in G.$

\bigskip

We seek conditions on $X,$ $A$ and translation invariant closed subspaces $%
{{\F}} $ of $BUC_{w}(G,X)$ which ensure that if $\psi \in {{\F}} $ then $%
\phi \in {{\F}}$. We apply Theorem 8.3, assuming throughout that $%
{{\F}} $ satisfies (5.2). Since $B\gamma =k\ast \gamma
=\widehat{k}(\gamma )\gamma $ for each $\gamma \in \widehat{G},$\
the characteristic function of $B$ is $\theta _{B}=\widehat{k}.$

\bigskip

Let $\alpha :\widehat{\mathbf{R}^{d}}\rightarrow \mathbf{R}^{d}$ be the
isomorphism given by $\alpha (\gamma _{s})=s$ where $\gamma
_{s}(t)=e^{i\left\langle s,t\right\rangle }$ for $s,t\in \mathbf{R}^{d}\ $%
and $\left\langle s,t\right\rangle =\sum_{j=1}^{d}s_{j}t_{j};$ and $\alpha :%
\widehat{\mathbf{Z}^{d}}\rightarrow \mathbf{T}^{d}$ be the isomorphism given
by $\alpha (\gamma _{\zeta })=\zeta $ where $\gamma _{\zeta }(n)=\zeta
^{n}=\prod_{j=1}^{d}\zeta _{j}^{n_{j}}$ for $\zeta \in \mathbf{T}^{d}\ $and $%
n\in \mathbf{Z}^{d}.$

\bigskip

A simple modification of the proof of Lemma 9.1 yields

\textbf{Lemma 10.1}. If $w$ has polynomial growth of order $N\
$and$\ u\in C_{c}^{N+d+1}(\alpha (\widehat{G})\mathbf{,}X),$ then
there exists $h\in L_{w}^{1}(G\mathbf{,}X)$ such that
$\widehat{h}=u\circ \alpha$.

\bigskip

\textbf{Lemma 10.2}. Conditions (8.2)-(8.6),(8.10) all hold. If
${{\F}} \in \{AP_{w}(\mathbf{R},X),$ $C_{w,0}(J,X),$
$0_{\mathbf{R}}\}$ and $Y=L(X),$ then (8.8) holds. If $w$ has
polynomial growth, condition (8.9) holds.

\,

\begin{proof} We have assumed (8.2). Now (8.3),(8.5),(8.6) are
clear,
(8.4) follows from Lemma 8.2, (8.7) from Lemma 8.8 and (8.10) from (5.2). If $%
{{\F}} \in \{AP_{w}(\mathbf{R},X)$, $C_{w,0}(J,X),$
$0_{\mathbf{R}}\}$ and $Y=L(X),$ then (8.8) holds by Lemma 8.9. If
$w$ has polynomial growth, (8.9) follows readily from Lemma 10.1.
\end{proof}

\smallskip

\textbf{Theorem 10.3}. Suppose $w\ $has polynomial growth, $\phi,
 \psi \in BUC_{w}(G,X),$ condition (8.7) holds and $k\ast \phi
=A\circ \phi +\psi$. If $\psi |_{J}\in {{\F}} $ , then
$\widehat{k}(sp_{{{\F}} }(\phi ))\subseteq \sigma (A).$

\,

\begin{proof} Using Lemma 10.2 the result follows from Theorem
8.3.
\end{proof}
\smallskip

\textbf{Corollary 10.4}. Suppose $w\ $has polynomial growth, $\phi
\in BUC_{w}(G,X),$ $\psi \in AP_{w}(G,X)$ and $k\ast \phi =A\circ
\phi +\psi$. Assume $(\widehat{k})^{-1}(\sigma (A))$ is residual.
Then $\phi \in AP_{w}(G,X)$ provided one of the following
conditions is satisfied.

\qquad (a) $w$ is bounded and $c_0\not\subset X;$

\qquad (b) $\gamma ^{-1}\phi \in E_{w}(G,X)$ for all $\gamma \in (%
\widehat{k})^{-1}(\sigma (A));$

\qquad (c) $w$ is bounded and $\phi (G)$ is  relatively weakly
compact.

\,

\begin{proof} The proof is identical to that of Corollary 9.5
with $\mathbf{R}$ replaced by $G$ and Lemma 9.3 by Lemma 10.2.
\end{proof}

\bigskip

\textbf{11. }$C_{0}$\textbf{-Semigroups}

In [50] , Ph\'{o}ng studied $C_{0}$-semigroups $%
\{T(t):t\in \mathbf{R}_{+}\}$ of operators such that $\left\|
T(t)\right\| $ is dominated by a (weight) function $\alpha (t)$
whose reduced (weight) function

$\alpha _{1}(s)=\limsup _{t\rightarrow \infty }\frac{\alpha
(t+s)}{\alpha (t)}$

\noindent has non-quasianalytic growth. He thereby extended
earlier results of Allan and Ransford [2], Arendt and Batty
\cite{ABcw} , Katznelson and Tzafriri \cite{KT}  and Ph\'{o}ng
\cite{PQ}. In this section, as corollaries of our main theorems,
we obtain results analogous to those in \cite{PQ},\cite{PV}. We
use dominating functions $u(t)$ which appear to be neither more
nor less general than those of Ph\'{o}ng.

\bigskip

Throughout we assume that $A$ is the generator of a $C_{0}$-semigroup $%
\{T(t):t\in \mathbf{R}_{+}\}$ of operators $T(t)\in L(X).$\ Moreover we
assume $\left\| T(t)\right\| $ $\leq u(t)$ for all $t\in \mathbf{R}_{+},$\
where $u(t)=v(t)w(t)$ for some weight $w$ of polynomial growth $N$,
satisfying therefore (1.1) - (1.3), (2.1) - (2.5); and $v:\mathbf{R}%
\rightarrow [1,\infty )$ satisfies (1.1), (2.2), (2.3) and \ \

\bigskip

(11.1)\qquad $v$ is differentiable on $\mathbf{R}_{+}$ and $\frac{v^{\prime }}{v}%
\in C_{0}(\mathbf{R}_{+},\mathbf{R}).$

\bigskip

An example of such a function is given by $v(t)=c\exp (1+\left| t\right|
)^{p},$ where $c\geq \frac{1}{e}$ and $0\leq p<1.$ This function also
satisfies the following condition which we will occasionally need.

\bigskip

(11.2)\qquad $v$ is constant or $\frac{1}{v}\in C_{w,0}(\mathbf{R}_{+},\mathbf{R}%
).$\

\bigskip

We will consider equations of the form $\phi ^{\prime }(t)=A\phi
(t)+\psi (t)$ for $t\in \mathbf{R}_{+},$ where $\frac{\phi }{v}\in
BUC_{w}(\mathbf{R}_{+},X)$ and $\frac{\psi }{v}\in C_{w,0}(\mathbf{R}%
_{+},X). $ Making the substitution

\bigskip

(11.3)\qquad $\phi _{1}=\frac{\phi }{v}$ and $\psi _{1}=\frac{\psi }{v}-%
\frac{v^{\prime }}{v}.\frac{\phi }{v}$

\smallskip

\noindent we find $\phi _{1}^{\prime }(t)=A\phi _{1}(t)+\psi _{1}(t)$ where $%
\phi _{1}\in BUC_{w}(\mathbf{R}_{+},X)$ and $\psi _{1}\in C_{w,0}(\mathbf{%
R}_{+},X)$.
 The study of the asymptotic behaviour of $\phi $\ is
reduced
to that of $\phi _{1}$. For $J\in \{\mathbf{R}_{+},\mathbf{R}\}$, let $%
L_{u}^{1}(J)=\{f\in L_{w}^{1}(\mathbf{R}):f=0$ on $\mathbf{R}$ $\backslash $
$J$ and $vf\in L_{w}^{1}(\mathbf{R})\}.$ If $f\in L_{u}^{1}(J),$ the
integral $\widehat{f}(T)=\int_{J}T(s)f(s)ds$ exists as a strongly convergent
Bochner integral. As before, let $\alpha :\widehat{\mathbf{R}}%
\rightarrow \mathbf{R}$ be the natural isomorphism given by $\alpha (\gamma
_{s})=s$ where $\gamma _{s}(t)=$ $e^{ist}$.

\bigskip

With different restrictions on the dominating function $u$, and with the
additional condition that $f$\ is of $w$-spectral synthesis with respect to $%
\alpha ^{-1}(i\sigma (A)\cap \mathbf{R),}$ the following is [50,
Theorem 8] .

\bigskip

\textbf{Theorem 11.1}. Let $f\in L_{u}^{1}(\mathbf{R}_{+}).$ If $\widehat{f}%
=0$ on $\Delta=\alpha ^{-1}(i\sigma (A)\cap \mathbf{R)}$ and supp
$\widehat{f}\cap \Delta$ is residual then $\lim_{t\rightarrow
\infty }\frac{\left\| T(t)\widehat{f}(T)\right\| }{u(t)}=0$.

\begin{proof} Take ${{\F}} =C_{w,0}(\mathbf{R}_{+},X).$\ Given
$x\in X$ set $\xi (t)=T(t)\widehat{f}(T)x$ and $\xi _{1}=\xi /v.$\
For the moment suppose $x\in D(A^{2}).$ Then,\ for $t\geq 0$\
define $\phi (t)=T(t)x$ and $\psi (t)=0.$ Making the substitution
(11.3) on $\mathbf{R}_{+}$ we find
$\frac{\Delta _{h}\phi _{1}(t)}{w(t)}=\frac{T(t)}{u(t)}(T(h)-I)x-\frac{%
\phi _{1}(t+h)}{w(t+h)}\frac{w(t+h)}{w(t)}\frac{\Delta
_{h}v(t)}{v(t)}.$
Since $v$ is $v$-uniformly continuous, it follows that $\phi _{1}$ is $w$%
-uniformly continuous and hence from (2.7) that $\phi _{1}\in BUC_{w}(%
\mathbf{R}_{+},X).$ Extend $\phi _{1},\psi _{1}$ to $\mathbf{R}$\
by defining, for $t\leq 0,$ $\phi _{1}(t)=x\cos t+Ax\sin t$ and
$\psi
_{1}(t)=\psi _{1}(0)-(x+A^{2}x)\sin t.$ So $\phi _{1}\in BUC_{w}(\mathbf{R%
},X),$ $\psi _{1}|_{\mathbf{R}_{+}}\in {{\F}} $\ and $\phi
_{1}^{\prime }=A\circ \phi _{1}+\psi _{1}$ on $\mathbf{R}.$ By
Lemma 9.3 and Theorem
9.4, $i\alpha (sp_{{{\F}} }(\phi _{1}))\subseteq \sigma (A)\cap i%
\mathbf{R}=i\alpha (\Delta^{-1} )$. If $g(t)=f(-t),$ then
$\widehat{g}=0$\ on $\Delta^{-1} \mathbf{.}$ Also, $\xi
(t)=\int_{0}^{\infty }T(t+s)f(s)xds$ and so $\xi _{1}=\xi
_{2}|_{\mathbf{R}_{+}}+\xi _{3}$ where $\xi _{2}=\phi _{1}\ast g$
and $\xi _{3}=\int_{0}^{\infty }\left( \phi _{1}\right)
_{s}|_{\mathbf{R}_{+}}\frac{\Delta _{s}v}{v}f(s)ds.$ Hence
$sp_{{{\F}}
}(\xi _{2})\subseteq \Delta^{-1} \cap $ supp$(\widehat{g}).$\ By Corollary 6.12, $%
\gamma ^{-1}\xi _{2}\in E_{w}^{0}(\mathbf{R},X)$ for all $\gamma
\in \Delta^{-1} . $ From the  assumptions it follows $\Delta^{-1}
\cap $ supp$(\widehat{g})$
 is residual. Since $%
{{\F}} $ is a $\Lambda _{w}^{0}$-class, Theorem 6.5 shows that $\xi _{2}|_{%
\mathbf{R}_{+}}\in {{\F}}$. Note also that, by (2.8), the functions $%
s\longmapsto \left( \phi _{1}\right) _{s}|_{\mathbf{R}_{+}}:$\ $\mathbf{R}%
_{+}\rightarrow BUC_{w}(\mathbf{R}_{+},X)$ and $s\longmapsto \frac{\Delta
_{s}v}{v}:$\ $\mathbf{R}_{+}\rightarrow C_{0}(\mathbf{R}_{+},\mathbf{R})$
are continuous. Hence the integrand for $\xi _{3}$\ is an almost
separably-valued, absolutely integrable function from $\mathbf{R}_{+}$\ to $%
{{\F}}$.\ So $\xi _{3}\in {{\F}}$.\ Hence $\xi _{1}\in {{\F}}$.
But, $D(A^{2})$ is dense in $X$ and so $\xi _{1}\in {{\F}} $ for
all $x\in
X.$ Therefore we have established the strong convergence of $\frac{T(t)%
\widehat{f}(T)}{u(t)}$ to $0$. Now the argument of Ph\'{o}ng [\
49, p.237]\ gives convergence in the operator norm. \end{proof}



\smallskip \

The following lemma indicates some situations in which functions
are necessarily $w$-ergodic. In the theorem, $A^{\prime }$ denotes
the conjugate of $A$ and $P\sigma (A^{\prime })$ its point
spectrum.

\bigskip

\textbf{Lemma 11.2}. (a) If $\frac{\phi }{v},\frac{\phi ^{\prime
}}{v}\in BUC_{w}(J,X),$ then $\frac{\phi ^{\prime }}{v}\in
E_{w}^{0}(J,X).$

(b) If $\gamma _{s}\in \widehat{\mathbf{R}}$ and $x\in $ range$%
(A-is),$\ then $\frac{\gamma _{s}^{-1}T(.)x}{v}\in E_{w}^{0}(J,X).$

(c) If (11.2) holds, $\gamma _{s}\in \widehat{\mathbf{R}}$ and $%
x\in \ker (A-is),$\ then $\frac{\gamma _{s}^{-1}T(.)x}{v}\in E_{w}(J,X).$ If
also $v$ or $w$ is unbounded, then $\frac{\gamma _{s}^{-1}T(.)x}{v}\in
E_{w}^{0}(J,X).$

\begin{proof} (a) Since $(\frac{\phi }{v})^{\prime }=\frac{\phi
^{\prime }}{v}-\frac{\phi }{v}\frac{v^{\prime }}{v}\in
BUC_{w}(J,X)$, it follows from Proposition 6.1 (c) that
$(\frac{\phi }{v})^{\prime }\in E_{w}^{0}(J,X).$ But $\frac{\phi
}{v}\frac{v^{\prime }}{v}\in C_{w,0}(J,X)\subset E_{w}^0(J,X)$, so
$\frac{\phi ^{\prime }}{v}\in E_{w}^{0}(J,X).$

(b) If $x=(A-is)y$ for some $y\in D(A),$ then $\gamma _{s}^{-1}T(.)x=(\gamma
_{s}^{-1}T(.)y)^{\prime }.$ Hence

$\frac{\gamma _{s}^{-1}T(.)y}{v},\frac{%
(\gamma _{s}^{-1}T(.)y)^{\prime }}{v}\in BUC_{w}(J,X).$

\noindent By (a), $\frac{%
\gamma _{s}^{-1}T(.)x}{v}=\frac{(\gamma _{s}^{-1}T(.)y)^{\prime }}{v}\in
E_{w}^{0}(J,X).$

(c) If $x\in $ $\ker (A-is),$\ then $(\gamma
_{s}^{-1}(t)T(t)x)^{\prime }=\gamma _{s}^{-1}(t)T(t)(A-is)x=0.$ So
$\gamma
_{s}^{-1}T(.)x$\ is constant and by (11.2), $\frac{\gamma _{s}^{-1}T(.)x}{v}%
\in E_{w}(J,X).$ If also $v$ or $w$ is unbounded, then
$\frac{\gamma _{s}^{-1}T(.)x}{vw}\in C_{0}(J,X).$ So $\frac{\gamma
_{s}^{-1}T(.)x}{v}\in E_{w}^{0}(J,X).$ \end{proof}

\smallskip

Part (a) of the following theorem, with different dominating
functions $u$, is [50, Theorem 7] .

\bigskip

\textbf{Theorem 11.3}. If $\sigma (A)\cap i\mathbf{R}$ is countable, then $%
\lim_{t\rightarrow \infty }\frac{\left\| T(t)x\right\| }{u(t)}=0$ for each $%
x\in X$ provided one of the following conditions is satisfied.

\qquad (a) $P\sigma (A^{\prime })\cap i\mathbf{R}=\emptyset ;$

\qquad (b) range$(A-is)$ is dense in $X$\ for each $s\in \mathbf{R};$

\qquad (c) either $v$ or $w$ is unbounded, (11.2) holds and $%
\ker (A-is)+$ range$(A-is)$ is dense in $X$\ for each $s\in \mathbf{R}$.

\,

\begin{proof} It is clear that (a) and (b) are equivalent. Let
${{\F}} =C_{w,0}(J,X),$ a $\Lambda _{w}^{0}$-class. Given $x\in
X$\ define $\phi _{1},\psi _{1}$ as in the proof of Theorem 11.1.
As there, $i\alpha (sp_{{{\F}} }(\phi _{1}))\subseteq \sigma
(A)\cap i\mathbf{R.}$ Recall that $E_{w}^{0}(\mathbf{R}_{+},X)$ is
closed in $BC_{w}(\mathbf{R}_{+},X).$
Hence, by Lemma 11.2, $\gamma ^{-1}\phi _{1}|_{\mathbf{R}_{+}}=\frac{%
\gamma ^{-1}T(.)x}{v}\in E_{w}^{0}(\mathbf{R}_{+},X)$ for all
$\gamma \in \widehat{\mathbf{R}}.$ By Theorem 6.5, $\phi
_{1}|_{\mathbf{R}_{+}}\in {{\F}}$. This means $\lim_{t\rightarrow
\infty }\frac{\left\| T(t)x\right\| }{u(t)}=0.$
\end{proof}
\smallskip

As a final application in this section we prove the following
theorem. Parts (a) and (c) are due to Lyubich, Matsaev and
Fel'dman \cite{LMF}.

\bigskip

\textbf{Theorem 11.4}. Assume the $C_{0}$-group $\{T(t):t\in \mathbf{R}\},$\
with generator $A,$ is dominated by a weight $w$ with polynomial growth of
order $N\ $. Let $f\in L_{w}^{1}(\mathbf{R}).$

(a) $\sigma (A)\cap i\mathbf{R}$ is non-empty.

(b) If $\sigma (A)\cap i\mathbf{R=\{}is\mathbf{\}}$, then $(A-is)^{N+1}=0.$

(c) If $\widehat{f}=0$ on $\alpha ^{-1}(i\sigma (A))$ and
$\alpha^{-1} (i\,\sigma (A))\cap $ supp $\widehat{f}$ is residual,
then $\widehat{f}(T)=0.$

(d) If $f$ is $w$-spectral synthesis with respect to
$\Delta=\alpha ^{-1}(i\sigma (A)\cap \mathbf{R)}$, then
$\widehat{f}(T)=0.$

(e) If $\sigma (A)\cap $ $-i\alpha ($supp$\widehat{f})=\{is\},$ then $%
(T(t)-\gamma _{s}(t))^{N+1}\widehat{f}(T)=0$ and

\qquad\qquad\qquad $(A-is)^{N+1}\widehat{f}\, (T)=0$.

\begin{proof} Take ${{\F}} =0_{\mathbf{R}}$ and set $g(t)=f(-t).$ Given $%
x\in X$ set $\phi (t)=T(t)x$ and $\xi (t)=T(t)\widehat{f}(T)x=\phi
\ast g(t).$ If $x\in D(A),$ then $\phi ^{\prime }=A\circ \phi $ on $%
\mathbf{R}.$ By Lemma 9.3 and Theorem 9.4, $i\alpha (sp_{w}(\phi
))\subseteq \sigma (A)\cap i\mathbf{R.}$ If $x\neq 0$ then $\phi
\neq 0$
and so $sp_{w}(\phi )\not = \emptyset ,$ proving (a). If $\sigma (A)\cap i%
\mathbf{R=\{}is\mathbf{\}}$, then $sp_{w}(\phi )=\{\gamma _{s}\}.$
By
Theorem 6.6, $\gamma _{s}^{-1}\phi \in P^{N}(\mathbf{R},X)$ and so $%
\left( \gamma _{s}^{-1}\phi \right)
^{(N+1)}(t)=T(t)(A-is)^{N+1}x=0$ for all $x\in D(A^{N+1}).$ Hence
(b) follows. Next, $sp_{w}(\xi )\subseteq \alpha ^{-1}(-i\sigma
(A))\cap $ supp$(\widehat{g}).$  The assumptions imply
$sp_{w}(\xi )$\ is residual. As $\widehat{g}=0$ on $%
\alpha ^{-1}(-i\sigma (A))$, it follows from Corollary 6.12 that
$\gamma ^{-1}\xi \in E_{w}^{0}(\mathbf{R},X)$ for all $\gamma \in
\alpha ^{-1}(-i\sigma (A)).$ By Theorem 6.5, $\xi \in {{\F}} $ for
all $x\in D(A)$ and therefore $\widehat{f}(T)=0$,  proving (c).
Now let $\Delta=\alpha ^{-1}(i\sigma (A)\cap \mathbf{R)}$ and
choose $(f_n)\subset L^1_w (\mathbf{R}, \mathbf{C})$, supp
$\widehat {f_n} \cap \Delta=\emptyset$  for all $n\in \mathbf{N}$
and $||f_n - f||_{L^1_w}\to 0$ as $n\to \infty$. Then by (c), it
follows $\widehat{f_n}(T)=0$ for all $n\in \mathbf{N}$. This
implies $\widehat{f}(T)=0$, proving (d).
Finally, if $\alpha ^{-1}(-i\sigma (A))\cap $ supp$(\widehat{g}%
)=\{\gamma _{s}\}$ then $sp_{w}(\xi )=\{\gamma _{s}\}.$\ By Theorem 6.6, $%
\gamma _{s}^{-1}\xi \in P^{N}(\mathbf{R},X).$ So $\Delta _{t}^{N+1}(\gamma
_{s}^{-1}\xi )=0$ for all $t\in \mathbf{R}$ and hence $(T(t)-\gamma
_{s}(t))^{N+1}\widehat{f}(T)x=0.$ Also, $(\gamma _{s}^{-1}\xi
)^{(N+1)}(t)=\gamma _{s}^{-1}(t)T(t)(A-is)^{N+1}\widehat{f}(T)x=0$ for all $%
t\in \mathbf{R}$ and $x\in D(A^{N+1})$. So $(A-is)^{N+1}\widehat{f}(T)=0.$
\end{proof}

\bigskip

\textbf{12. Recurrence equations}

Throughout this section we assume that $A$ satisfies (8.2) and $J\in \{%
\mathbf{Z,Z}_{+}\}.\ $Also, the operator
$B:BUC_{w}(\mathbf{Z},X)\rightarrow BUC_{w}(\mathbf{Z},X)$ \ is
given by $(B\phi
)(n)=\sum_{j=0}^{m}b_{j}\phi (n+n_{j})$ for some $b_{j}\in \mathbf{C},$ $%
n_{j}\in \mathbf{Z}.$ Given $\psi \in $ $BUC_{w}(\mathbf{Z},X)$ we
investigate the behaviour of solutions $\phi
:\mathbf{Z}\rightarrow D(A)$ of the recurrence equation

\bigskip

(12.1)\qquad $(B\phi )(t)=A\phi (t)+\psi (t)$ for $t\in
\mathbf{Z}.$

\bigskip

We seek conditions on $X,\,$ $A\,$ and translation invariant closed subspaces $%
\,\,{{\F}} $ of\,\,\,\, $BUC_{w}(J,X)$ which ensure that if $\psi
|_{J}\in {{\F}} $ then $\phi |_{J}\in {{\F}}$. We apply Theorem
8.3, assuming throughout that ${{\F}} $ satisfies (5.2).

\smallskip

Let $\mathbf{T}$ denote the circle group and $\alpha :\widehat{\mathbf{Z}}%
\rightarrow \mathbf{T}$ \ the isomorphism given by $\alpha (\gamma _{\zeta
})=\zeta $ where $\gamma _{\zeta }(n)=\zeta ^{n}$ for $n\in \mathbf{Z},$ $%
\zeta \in \mathbf{T}.$ The space of $m$-times continuously differentiable
functions $u:\mathbf{T\rightarrow }X$ is $C^{m}(\mathbf{T,}X).$

\smallskip

To find the characteristic function of $B$ let $p_{B}(\zeta
)=\sum_{j=0}^{m}b_{j}\zeta ^{n_{j}}$ where $\zeta \in \mathbf{C}.$ For $%
\gamma \in \widehat{\mathbf{Z}}$ and $x\in X$ we have $B(\gamma
x)=p_{B}(\alpha (\gamma ))\gamma x$ and so $\theta _{B}=p_{B}\circ
\alpha$.

\smallskip

A simple modification of the proof of Lemma 9.1 yields

\smallskip

\textbf{Lemma 12.1}. If $w$ has polynomial growth $N$\ and $u\in C^{N+2}(%
\mathbf{T,}X),$ then there exists $h\in L_{w}^{1}(\mathbf{Z,}X)$ such that $%
\widehat{h}=u\circ \alpha$.

\bigskip

\textbf{Lemma 12.2}. Conditions (8.2)-(8.7), (8.10) all hold. If $%
{{\F}} \in \{AP_{w}(\mathbf{Z},X),C_{w,0}(J,X),$ $0_{\mathbf{Z}}\}$ and $%
Y=L(X),$ then (8.8) holds. If $w$ has polynomial growth, condition
(8.9) holds.

\,

\begin{proof} We have assumed (8.2). Now (8.3)-(8.6) are clear
and (8.7)
follows from Lemma 8.8. If ${{\F}} \in \{AP_{w}(\mathbf{Z}%
,X),C_{w,0}(J,X), $ $0_{\mathbf{Z}}\}$ and $Y=L(X),$ then (8.8)
holds by Lemma 8.9. Let $w$ have polynomial
growth. To prove (8.9) let $V$ be a relatively compact open subset of $%
\widehat{\mathbf{Z}}$ with $\theta _{B}(\overline{V})\cap \sigma
(A)=\emptyset$. So $\overline{V}\subseteq W$ where $W$ is also a
relatively
compact open subset of $\widehat{\mathbf{Z}}$ with $\theta _{B}(\overline{W}%
)\cap \sigma (A)=\emptyset$. Choose $v\in C^{N+2}(\mathbf{T,C})$ such that $%
v=1$ on $\theta _{B}\left( V\right) $ and supp$(v)\subseteq \theta _{B}(W).$
Setting $r(\zeta )=v(\zeta )(\zeta -A)^{-1}$ and $u=r\circ p_{B}$\ we have $%
u\in C^{N+2}(\mathbf{T,}Y).$ By Lemma 12.1, there exists $h\in L_{w}^{1}(%
\mathbf{Z,}Y)$ such that $\widehat{h}=u\circ \alpha $ and so $\widehat{h}%
(\gamma )=(\theta _{B}(\gamma )-A)^{-1}$ for $\gamma \in V,$
proving (8.9).
\end{proof}
\smallskip

\textbf{Theorem 12.3}. Suppose $w$\ has polynomial growth, $\phi
,\psi \in BUC_{w}(\mathbf{Z},X),$ condition (8.8) holds and $B\phi
=A\circ \phi +\psi$. If $\psi |_{J}\in {{\F}} $ , then $\theta
_{B}(sp_{{{\F}} }(\phi ))\subseteq \sigma (A)\cap
p_{B}(\mathbf{T}).$

\,

\begin{proof} Using Lemma 12.2 the result follows from Theorem
8.3. \end{proof}

\smallskip

\textbf{Corollary 12.4}. Suppose $w$\ has polynomial growth, $\phi
,\psi \in BUC_{w}(\mathbf{Z},X),$ $X$ is a unital Banach algebra,
${{\F}} $ is closed under multiplication by elements of $X$ and
$B\phi =x\phi +\psi$. If $\psi |_{J}\in {{\F}} $ , then $\theta
_{B}(sp_{{{\F}} }(\phi ))\subseteq \sigma (x)\cap
p_{B}(\mathbf{T}).$

\,

\begin{proof} For $y\in X,$ define $L_{y}:X\rightarrow X$ by
$L_{y}z=yz.$ Apply Theorem 12.3 with $J=\mathbf{Z},$ $A=L_{x}$\
and $Y=\{L_{y}:y\in X\}.$ Since (8.8) is clear and $\sigma
(A)=\sigma (x)$, we are finished.
\end{proof}

\smallskip

We readily obtain the following two results. For the case $w=a=1,$
the first was proved by Gelfand and the second by Katznelson and
Tzafriri (see [40]\ and the references therein, and [16] ).
Corollary 12.5 is due to Hille [35, Theorem 4.10.1]\ and Corollary
12.7 to Allan and Ransford [2]. However, our proof is new. Here
$X$ is a Banach algebra with unit $e$. We shall say that an
element $x\in X$ is
\textit{power dominated} by $w$\ if $\{\frac{x^{n}}{w(n)}:n\in \mathbf{Z}%
_{+}\}$ is bounded, and \textit{doubly power dominated \ }by $w$\ if $\{%
\frac{x^{n}}{w(n)}:n\in \mathbf{Z}\}$ is bounded. Let $L_{w}^{1}(\mathbf{Z}%
_{+})=\{f\in L_{w}^{1}(\mathbf{Z}):f(n)=0$ for $n<0\}.$ If $x$ is
power dominated by $w$ and $f\in L_{w}^{1}(\mathbf{Z}_{+})$ we
define $\widehat{f}\, (x)=\sum_{n=0}^{\infty }f(n)x^{n}$.

\bigskip

\textbf{Corollary 12.5}. Assume $w$\ has polynomial growth of
order $N.$\ If an element $x\in X$ is doubly power dominated by
$w$\ and $\sigma (x)=\{1\},$ then $(x-e)^{N+1}=0.$

\,

\begin{proof} Apply Corollary 12.4 with $J=\mathbf{Z},$ ${{\F}}
=\{0\}, $ $B\phi =\phi _{1}$ and $\phi (n)=x^{n}$.  We conclude
that $\alpha (sp_{{{\F}} }(\phi ))\subseteq \sigma (x)\cap
\mathbf{T}=\{1\}$. By
Proposition 6.3 (a)(vi), $(x-e)^{N+1}=\sum_{j=0}^{N+1}\binom{N+1}{j}%
(-1)^{j}x^{N+1-j}=\Delta _{1}^{N+1}\phi (0)=0$.
\end{proof}

\smallskip

\textbf{Corollary 12.6}. Assume $a:\mathbf{Z}_+  \to (0,\infty]$,
$\lim_{n\rightarrow \infty }\frac{a(n+1)}{%
a(n)}=1$\ and the element $x\in X$ is power dominated by $a$.

(a) If $\sigma (x)\cap \mathbf{T}=\emptyset ,$ then $\left\| \frac{x^{n}}{%
a(n)}\right\| \rightarrow 0$ as $n\rightarrow \infty $.

(b) If $\sigma (x)\cap \mathbf{T}=\{1\},$ then $\left\| \frac{x^{n+1}-x^{n}}{%
a(n)}\right\| \rightarrow 0$ as $n\rightarrow \infty$.

\,

\begin{proof} Apply Theorem 12.3 with $J=\mathbf{Z}_{+}$,  $w=1$,
 ${{\F}}
=C_{0}(J,X)$,  $B\phi =\phi _{1}$ and $\phi$ , $\psi $ as follows.
 Set

 $\phi (n)=\frac{x^{n}}{a(n)}$,\qquad $\psi (n)=\frac{x^{n+1}}{a(n+1)}%
(1-\frac{a(n+1)}{a(n)})$, \qquad $%
n\geq 0$,

$\phi (n)=\frac{e}{a(0)}$,\qquad $\psi (n)=\frac{e-x}{a(0)}$,
\qquad \qquad\qquad $n<0$.

\noindent  So $\alpha (sp_{{{\F}} }(\phi ))\subseteq \sigma
(x)\cap \mathbf{T}$. If $\sigma (x)\cap \mathbf{T}=\emptyset ,$
then by Proposition 6.3 (a) (v), $\phi |_{J}\in {{\F}}$.\ If
$\sigma (x)\cap \mathbf{T}=\{1\},$ then by Proposition 6.3
(a)(vi), $\Delta _{1}\phi \in C_{0}(J,X).$
\end{proof}

\smallskip

\textbf{Corollary 12.7}. Assume $w$\ has polynomial growth$,$\ $%
\phi \in BUC_{w}(\mathbf{Z},X),$ $\psi \in AP_{w}(\mathbf{Z},X)$ and $%
B\phi =A\circ \phi +\psi$. Assume $p_{B}$ is not constant and
$\sigma
(A)\cap p_{B}(\mathbf{Z})$ is residual. Then $\phi \in AP_{w}(\mathbf{Z}%
,X)$ provided one of the following conditions is satisfied.

\qquad (a) $w$ is bounded and $c_0\not\subset X;$

\qquad (b) $\gamma ^{-1}\phi \in E_{w}(\mathbf{Z},X)$ for all
$\gamma \in \theta _{B}^{-1}(\sigma (A));$

\qquad (c) $w$ is bounded and $\phi (\mathbf{Z})$ is relatively
weakly compact.

\,

\begin{proof} Let ${{\F}} =AP_{w}(\mathbf{Z},X)$ and $Y=L(X).$ Since $%
\sigma (A)\cap p_{B}(\mathbf{Z})$ is residual and $p_{B}$ is a
non-constant polynomial, $\theta _{B}^{-1}(\sigma (A))$\ is
residual. By Lemma 12.2, the result follows from Corollary 8.8.
\end{proof}

\end{document}